\let\latexarabic\arabic
\let\latexdocument\document
\let\latexenddocument\enddocument

\RequirePackage[thmmarks]{ntheorem}
\makeatletter
\renewtheoremstyle{plain} 
  {\item[\hskip\labelsep \theorem@headerfont ##1\ \textup{##2}\theorem@separator]} 
  {\item[\hskip\labelsep \theorem@headerfont ##1\ \textup{##2}\ (##3)\theorem@separator]}
\makeatother

\documentclass[lineno]{biometrika}

\let\document\latexdocument
\let\enddocument\latexenddocument
\AtEndDocument{\printhistory}
\let\arabic\latexarabic
\def\rm{}

\usepackage{times}
\usepackage{amsmath,amssymb}
\usepackage{natbib}
\usepackage{bm}
\usepackage{enumerate}
\usepackage{graphicx}
\usepackage{caption}
\usepackage{subcaption}
\usepackage{bbm}
\usepackage{mathrsfs}
\usepackage{color,soul} 
\usepackage{xcolor} 
\usepackage[plain,noend]{algorithm2e}
\usepackage{url}
\usepackage{comment}

\makeatletter
\renewcommand{\algocf@captiontext}[2]{#1\algocf@typo. \AlCapFnt{}#2} 
\def\@algocf@capt@plain{top}
\renewcommand{\algocf@makecaption}[2]{%
  \addtolength{\hsize}{\algomargin}%
  \sbox\@tempboxa{\algocf@captiontext{#1}{#2}}%
  \ifdim\wd\@tempboxa >\hsize
    \hskip .5\algomargin%
    \parbox[t]{\hsize}{\algocf@captiontext{#1}{#2}}
  \else%
    \global\@minipagefalse%
    \hbox to\hsize{\box\@tempboxa}
  \fi%
  \addtolength{\hsize}{-\algomargin}%
}
\makeatother

\def\Tr{{ \mathrm{\scriptscriptstyle T} }}

\newcommand{\var}{\mathrm{Var}}
\newcommand{\cov}{\mathrm{Cov}}
\newcommand{\tr}{\mathrm{tr}}
\newcommand{\mr}{\mathrm}

\newcommand{\DD}{\mathcal{D}\mathcal{D}}
\newcommand{\T}{\mathcal{T}}

\newcommand{\D}{\mathcal{D}}

\newcommand{\Z}{\mathcal{Z}}

\begin{document}
\nolinenumbers
\jname{Submitted to Biometrika}


\markboth{Y. Zhu et~al.}{Radial Neighbors Gaussian Processes}

\title{Radial Neighbors for Provably Accurate Scalable Approximations of Gaussian Processes}

\author{Yichen Zhu}
\affil{Bocconi Institute for Data Science and Analytics, Bocconi University, Milan, MI, Italy.
\email{yichen.zhu@unibocconi.it}}

\author{Michele Peruzzi}
\affil{Department of Biostatistics, University of Michigan, Ann Arbor, Michigan, U.S.A. \email{peruzzi@umich.edu}}

\author{Cheng Li}
\affil{Department of Statistics and Data Science, National University of Singapore, Sinagpore.
\email{stalic@nus.edu.sg}}

\author{\and David B. Dunson}
\affil{Department of Statistical Science \& Mathematics, Duke University, North Carolina, U.S.A. \email{dunson@duke.edu}}

\maketitle

\begin{abstract}
In geostatistical problems with massive sample size, Gaussian processes can be approximated using sparse directed acyclic graphs to achieve scalable $O(n)$ computational complexity. In these models, data at each location are typically assumed conditionally dependent on a small set of parents which usually include a subset of the nearest neighbors. These methodologies often exhibit excellent empirical performance, but the lack of theoretical validation leads to unclear guidance in specifying the underlying graphical model and sensitivity to graph choice. We address these issues by introducing radial neighbors Gaussian processes (RadGP), a class of Gaussian processes based on directed acyclic graphs in which directed edges connect every location to all of its neighbors within a predetermined radius. 
We prove that any radial neighbors Gaussian process can accurately approximate the corresponding unrestricted Gaussian process in Wasserstein-2 distance, with an error rate determined by the approximation radius, the spatial covariance function, and the spatial dispersion of samples. 
We offer further empirical validation of our approach via applications on simulated and real world data showing excellent performance in both prior and posterior approximations to the original Gaussian process.
\end{abstract}


\section{Introduction}
Data indexed by spatial coordinates are routinely collected in massive quantities due to the rapid pace of technological development in imaging sensors, wearables and tracking devices.
Statistical models of spatially correlated data are most commonly built using Gaussian processes  due to their flexibility and their ability to straightforwardly quantify uncertainty. Spatial dependence in the data is typically characterized by a parametric covariance function whose unknown parameters are estimated by repeatedly evaluating a multivariate normal density with a dense covariance matrix. At each density evaluation, the inverse covariance matrix and its determinant must be computed at a complexity that is cubic on the data dimension $n$. This steep cost severely restricts the applicability of Gaussian processes on modern geolocated datasets. This problem is exacerbated in Bayesian hierarchical models computed via Markov-chain Monte Carlo (MCMC) because several thousands of $O(n^3)$ operations must be performed to fully characterize the joint posterior distribution.

A natural idea for retaining the flexibility of Gaussian processes while circumventing their computational bottlenecks is to use approximate methods for evaluating high dimensional multivariate Gaussian densities. If the scalable approximation is ``good'' in some sense, then one can use the approximate method for backing out inferences on the original process parameters. A multitude of such scalable methods has been proposed in the literature. 
Low-rank methods using inducing variables or knots \citep{QuiRas05, CreJoh08,Banetal08,Finetal09,Guhetal11,Sanetal11} are only appropriate for approximating very smooth Gaussian processes because the number of inducing variables must increase polynomially fast with the sample size to avoid oversmoothing of the spatial surface \citep{Stein14,Buretal20}. Scalability can also be achieved by assuming sparsity of the Gaussian covariance matrix via compactly supported covariance functions \citep{Furetal06,Kauetal08,Bevetal19}; these methods operate by assuming marginal independence of pairs of data points that are beyond a certain distance from each other. See also the related method of \cite{GraApl15}. Using a similar intuition, one can achieve scalability by approximating the high dimensional Gaussian density with a product of small dimensional densities that characterize dependence of nearby data points. Composite likelihood methods \citep{Baietal12,Eidetal14,BevGae15} approximate the high dimensional Gaussian density with the product of sub-likelihoods, whereas Gaussian Markov Random Fields \citep[GMRF;][]{cressie93statistics,rue2005gaussian} 
approximate it via a product of conditional densities whose conditioning sets include all the local spatial neighbors.

Markovian conditional independence assumptions based on spatial neighborhoods can be visualized as a sparse undirected graphical model; these assumptions are intuitively appealing and lead to sparse precision matrices and more efficient computations \citep{rue2001fast}. 
However, the normalizing constant of GMRF densities may require additional expensive computations, and additional approximation steps may be required when making predictions at new spatial locations because GMRFs do not necessarily extend to a standalone stochastic process. Both these problems can be resolved by ordering the data and conditioning only using previously ordered data points.
Because one can now use a directed acyclic graph to represent conditional independence in the data, this approximation  \citep[due to][]{Vec88} immediately leads to a valid density which can be extended to a standalone stochastic process.

There is a rich literature on approximations of Gaussian processes using directed acyclic graphs. In addition to extensions of Vecchia's method \citep{Datetal16a, Finetal19, KatGui21}, scalable models can be built on domain partitions \citep{Peretal21a}, and generalizations are available for non-Gaussian and multivariate data \citep{Zilber21vl, PerDun22, perdun22nongauss} as well as for modeling nonstationary processes \citep{Boretal21, kiddkatzfuss2022ba}.
However, we identify two issues with current approaches based on Vecchia's approximation. First, because the graphical model representation of Vecchia approximations is not unique,  the numerical performance of these methods can be sensitive to the choice of graphical model in certain scenarios \citep{guinness2018permutation}. There is no consensus on how to choose these conditioning sets.
Second, although their performance has been demonstrated in practice \citep{Heaetal19} and some theory on parameter estimation is available \citep{zhang2012,zhangetal2023}, there 
is limited theoretical study focused on accuracy of the approximation. The recent work by \cite{schafer2021sparse} established  Kullback-Leibler divergence bounds for Vecchia approximations, but under strong conditions that are difficult to verify for commonly used spatial covariance functions.

In this paper, we address these two issues by introducing a novel method and corresponding theoretical guarantees. Our radial neighbors Gaussian processes scalably approximate an initial Gaussian process by assuming conditional independence at spatial locations as prescribed by a directed acyclic graph. For every location in the graph, directed edges connect it to all of its neighbors within a predetermined radius $\rho$. 

We also introduce--to the best of our knowledge--the first general theoretical result regarding the closeness between a Vecchia-like method and the 
Gaussian process it approximates.  
We show that any process from the RadGP family can accurately approximate the original Gaussian processes in Wasserstein-2 distance, with an error rate determined by the approximation radius, the spatial covariance function, and the spatial dispersion of samples. This approximation accuracy is insensitive to the graphical structure as long as the graphical structure does not alter the radial condition of RadGP. Our theory can be applied to a wide range of covariance functions including the Mat\'{e}rn, Gaussian and generalized Cauchy models. 
Empirical studies demonstrate that our RadGPs have excellent performance in approximating a Gaussian process prior and the posterior from a Gaussian process regression model.

\section{Radial Neighbors Gaussian process} \label{sec:RadGP}

\subsection{Radial Neighbors Graphs}
\label{subsec:dag def}
Let the spatial domain $\Omega$ be a connected subset of $\mathbb{R}^d$. Let $Z = (Z_s: s\in\Omega)$ be a real valued Gaussian process on $\Omega\subset\mathbb{R}^d$ with zero mean and covariance function $K:\Omega\times\Omega \to \mathbb{R}$ such that for all $s_1, s_2\in\Omega$,
$\mathbb{E} (Z_{s_1}) = 0$ and $\cov (Z_{s_1}, Z_{s_2}) = K(s_1, s_2)$. The objective of this paper is to approximate the process $Z$ with another process $\hat{Z}$ such that the distance between $\hat{Z}$ and $Z$ is small in some metric and finite marginal densities of $\hat{Z}$ are easy to compute.

Let $\D = \{w_1,w_2,\ldots,w_n\}$ be a finite collection of locations in $\Omega$. A $\rho$-radial neighbors graph on $\D$ is defined to be a directed acyclic graph such that for any $w_i, w_j\in \D$ with $\|w_i-w_j\|_2<\rho$, there is a directed edge between $w_i$ and $w_j$. Such radial neighbors graphs are not unique. Section \ref{subsec:dag practice}  provides one practical way of constructing a radial neighbors graph that will be used in our experimental studies, while Section \ref{subsec:construction} presents how to build a Gaussian process from any radial neighbors graph. 
Our theory in Section \ref{sec:theory} applies to any Gaussian process built from a radial neighbors graph, not just the specific graph we constructed in Section \ref{subsec:dag practice}.

\subsection{Construction of a Radial Neighbors Graph}\label{subsec:dag practice}
We present a procedure to construct the radial neighbors graph that will be used for experiments in Section \ref{sec:experiments}.
First, we specify an order of all locations in $\D$. This is done by choosing a ``center'' location of $\D$ and sorting all the locations in $\D$ according to their distances to this ``center'' in a non-decreasing order. 
If the locations in $\D$ are roughly evenly distributed in the domain $\Omega$, then the coordinates of this center are set as the arithmetic mean of the coordinates of all locations in $\D$.

Next, once an order of locations in $\D$ is obtained, for each location $w_i\in\D$, we set the parent set of $w_i$ to be all the locations preceeding $w_i$ in this order and within $\rho$ distance to $w_i$. 
As a consequence, for all locations within $\rho$ radius of $w_i$, the ones that precede it are its parents, whereas the ones that follow it are its children.
This construction satisfies the definition of radial neighbors graphs in Section \ref{subsec:dag def} that any pair of locations $(w_i,w_j)$ with $\|w_i-w_j\|_2<\rho$ have a directed edge connecting them.

Depending on the spatial dispersion of locations in $\D$, some locations may not have neighbors within $\rho$ distance. We thus populate the parent set of a location $w_i \in \D$ without $\rho$-neighbors by drawing an edge from $w_i$'s nearest neighbor among all locations that precede $w_i$ in the order. This step ensures that the resulting directed acyclic graph is connected. 

\begin{figure}[h]
    \centering
    \includegraphics[width=\textwidth]{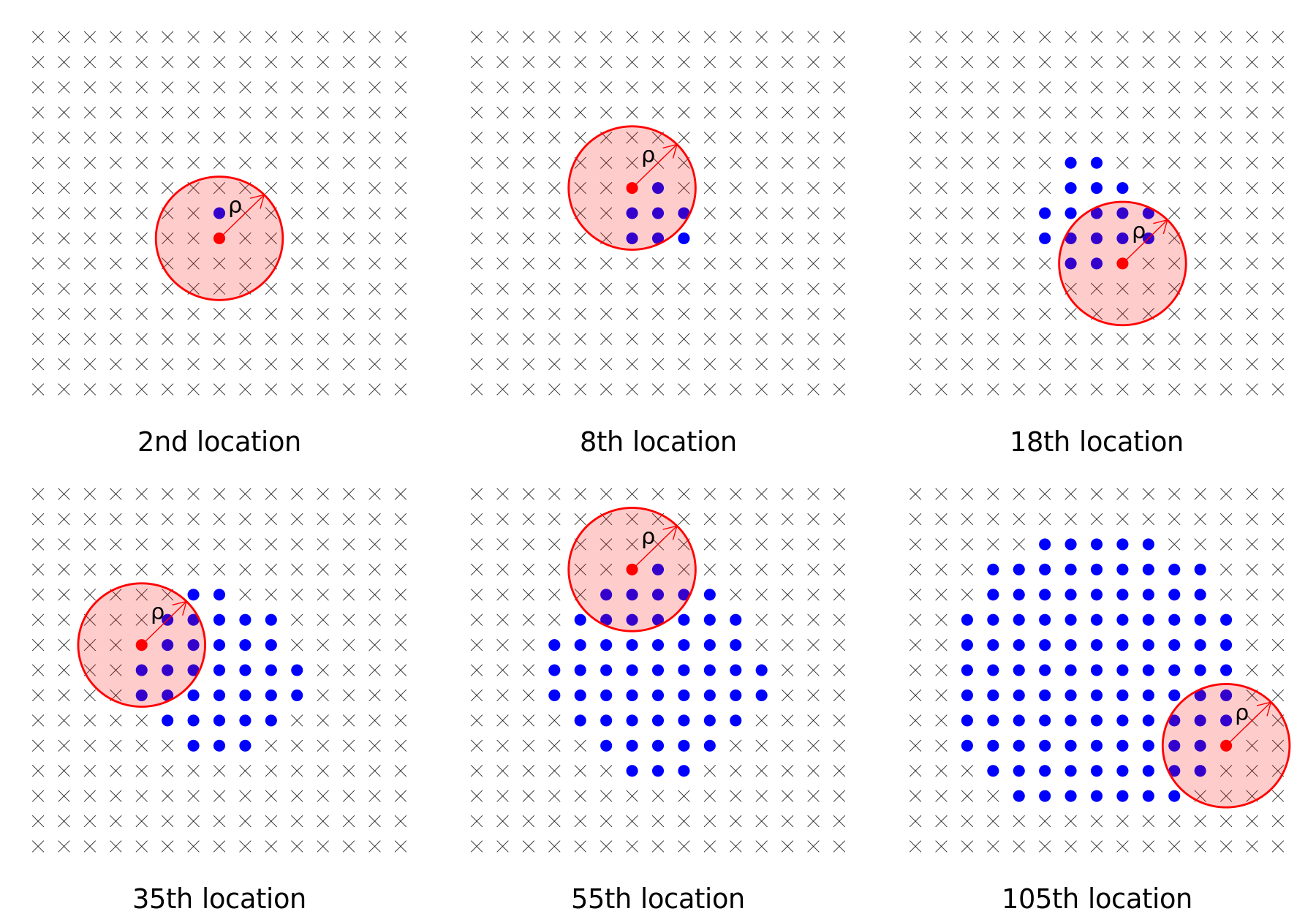}
    \caption{Construction of a radial neighbors graph. In each subpanel, the location in red follows all locations in blue and precedes all locations in gray in the predetermined order. The parent set for the location in red is built by drawing a directed edge to it from each blue location inside the red circle.
    }
    \label{fig:illu}
\end{figure}

If the locations in $\D$ are roughly evenly distributed in the domain $\Omega$, once the first $i$ locations are chosen and the directed graph on them is constructed, our procedure roughly chooses the $(i+1)$th location to be as close to the first $i$ locations as possible. 
Since spatial covariance function $K(s_1, s_2)$ is often relatively large when $s_1$ and $s_2$ are close to each other, our procedure 
aims to minimize the conditional variance of the $(i+1)$th location given the first $i$ locations, or in other words, maximize the spatial information that the $(i+1)$th location can borrow from the first $i$ locations. Figure \ref{fig:illu} illustrates the process of building this radial neighbors graph when $\D$ is a $15\times 15$ grid on $[0,1]^2$.

Finally, if we have a new data set $\tilde\D$, the radial neighbors graph constructed on $\D$ can be extended to a radial neighbors graph on $\D\cup\tilde\D$. Specifically, once an ordering of locations in $\tilde{\D}$ is obtained, we concatenate the ordering of locations in $\D$ and $\tilde{\D}$ to obtain an ordering of locations $\D\cup\tilde{\D}$. We then construct a radial neighbors graph on $\D\cup\tilde{\D}$ via the same procedure as described above. The final graph will contain the graph on $\D$ as a subgraph. This extension framework will be useful when extending a graph from training data to test data.

For some dataset, the spatial dispersion of both training and testing locations can be highly uneven. While a radial neighbor graph ensures all locations are connected in the same directed graph, it can still result in poor predictions on some test locations. We can mitigate this issue by artificially adding more test locations, so that each test locations receives a decent number of parents.

\subsection{Gaussian Processes from Radial Neighbors Graphs} \label{subsec:construction}
Given a $\rho$-radial neighbors graph on $\D=\{w_1,w_2,\ldots\,w_n\}$ defined in Section \ref{subsec:dag def}, we can build a RadGP on $\D$ in which each $w_i$ is conditionally only dependent on its parent set.  Let the symbol $\overset{d}{=}$ denote equality in distribution and define $\Sigma_{ A, B}$ as the covariance matrix between two spatial location sets $A$ and $B$ under the covariance function $K(\cdot,\cdot)$. The conditional distribution of $\hat{Z}_{w_i}$ is defined as
\begin{align}\label{eq alt-dag 2}
&[\hat{Z}_{w_i}\mid \hat{Z}_{w_j}, j<i] 
 \overset{d}{=}  [\hat{Z}_{w_i} \mid \hat{Z}_{\mr{pa}(w_i)}] \nonumber \overset{d}{=} [{Z}_{w_i}\mid {Z}_{\mr{pa}(w_i)}]\\
\sim & N\Big(\Sigma_{\mr{pa}(w_i),w_i}^\Tr \Sigma_{\mr{pa}(w_i),~ \mr{pa}(w_i)}^{-1}\hat{Z}_{\mr{pa}(w_i)},  \;K(w_i,w_i)- \Sigma_{\mr{pa}(w_i),w_i}^\Tr \Sigma_{\mr{pa}(w_i), ~\mr{pa}(w_i)}^{-1} \Sigma_{\mr{pa}(w_i),w_i} \Big),
\end{align}
where for $i=1$ the parent set $\mr{pa}(w_1)$ is empty.

Equation (\ref{eq alt-dag 2}) defines the distribution of $\hat{Z}$ on all finite subsets of $\D$, which  essentially says the conditional distribution of $\hat{Z}_{w_i}$ given the process $\hat{Z}$ on all locations preceding $w_i$ has the same law as the conditional distribution of $Z_{w_i}$ given the process $Z$ on the parent set $\mr{pa}(w_i)$.
We then extend this process $\hat Z$ to arbitrary finite subsets of the whole space $\Omega$. For all $s\in\Omega\backslash \D$, we define the parents of $s$ as
$\mr{pa}(s) = \{s'\in \mathcal{D}: \|s'-s\|_2<\rho\}$.
For any finite set $U\subset\Omega\backslash \D$, we define the conditional distribution of $\hat{Z}_U$ given $\hat{Z}_{\mathcal{T}_1}$ as
\begin{align}\label{eq alt-dag 3}
p(\hat{Z}_U \mid \hat{Z}_{\D}) 
& =  \prod_{s\in U} p(\hat{Z}_s \mid \hat{Z}_{\mr{pa}(s)})
= \prod_{s\in U} p({Z}_s \mid Z_{\mr{pa}(s)}), \nonumber\\
[\hat{Z}_s \mid \hat{Z}_{\mr{pa}(s)}] & \sim  N\Big(\Sigma_{\mr{pa}(s),s}^\Tr \Sigma_{\mr{pa}(s),\mr{pa}(s)}^{-1}\hat{Z}_{\mr{pa}(s)}, \;K(s,s)- \Sigma_{\mr{pa}(s),s}^\Tr \Sigma_{\mr{pa}(s),\mr{pa}(s)}^{-1} \Sigma_{\mr{pa}(s),s} \Big).
\end{align}

Equations (\ref{eq alt-dag 2}) and (\ref{eq alt-dag 3}) together give the joint distribution of $\hat{Z}$ on any finite subset of $\Omega$. Specifically, for a generic finite set $A = U\cup V$ with $U\subset\Omega\backslash\D$ and $V\subset\D$, 
the joint density of $\hat{Z}$ on $A$ implied by Equations (\ref{eq alt-dag 2}) and (\ref{eq alt-dag 3}) is
\begin{equation}\label{eq alt-dag 4}
p(\hat{Z}_{A})
= p(\hat{Z}_{U}|\hat{Z}_{\D}) \int  \prod_{i=1}^n p(\hat{Z}_{w_i}|\hat{Z}_{\mr{pa}(w_i)}) \prod_{w_i\in \D\backslash V} dw_i.
\end{equation}
Equations \eqref{eq alt-dag 2}-\eqref{eq alt-dag 4} complete the definition of our RadGP. 
\begin{lemma}\label{valid process}
The radial neighbors Gaussian process (RadGP) $\hat Z(\cdot)$ 
is a valid Gaussian process on the whole spatial domain $\Omega$.
\end{lemma}
The proof of Lemma \ref{valid process} relies on the Kolmogorov extension theorem and is available at Section S1 of the Supplementary Material.

\section{Theoretical Properties} \label{sec:theory}
\subsection{Wasserstein Distance Between Gaussian Processes} \label{subsec:W2dist}

The objective of this section is to establish theoretical upper bounds for the Wasserstein distance between the original Gaussian process and any RadGP on $\D$ described above. For two generic random vectors $Z_1,Z_2$ defined on the same space $\Z$ following the probability distributions $\mu_1,\mu_2$, their Wasserstein-2 ($W_2$) distance is
$$W_2(Z_1, Z_2) = \left\{\min_{\tilde{\mu}\in \mathcal{M}} \int_{\Z \times \Z} \|z_1-z_2\|_2^2 d\tilde{\mu}(z_1,z_2) \right\}^{1/2}, $$
where $\mathcal{M}$ is the set of all probability distributions on the product space $\Z\times\Z$ whose marginals are $\mu_1$ and $\mu_2$. Let $Z_{\D}$ and $\hat{Z}_{\D}$ be the random vectors from the original Gaussian process and the RadGP on the set $\mathcal{D}$. Denote the covariance matrix of $Z_{\D}$ and $\hat{Z}_{\D}$ as $\Sigma_{\D\D}$ and $\hat{\Sigma}_{\D\D}$, 
respectively. Then the Wasserstein-2 distance between $Z_{\D}$ and $\hat Z_{\D}$ has a closed form \citep{gelbrich1990formula}:
\begin{equation}\label{eq w2 closed}
W_2^2(Z_{\D},\hat{Z}_{\D}) = \tr(\Sigma_{\D\D}) + \tr(\hat{\Sigma}_{\D\D}) - 2\tr\{(\Sigma_{\D\D}^{1/2} \hat{\Sigma}_{\D\D} \Sigma_{\D\D}^{1/2})^{1/2} \},
\end{equation}
where $\tr(A)$ denotes the trace of a square matrix $A$. The right-hand side of equation (\ref{eq w2 closed}) involves various matrix powers that are difficult to analyze. Fortunately, for Gaussian measures, the squared Wasserstein-2 distance can be upper bounded by the trace norm of the difference between their covariance matrices, which greatly simplifies our analysis.

\subsection{Spatial Decaying Families}
We consider the Gaussian process $Z=(Z_s: s\in\Omega)$ with zero mean and an isotropic nonnegative covariance function $K(\cdot,\cdot)$. Thus we can reformulate $K$ as $K(s_1,s_2)=K_0(\|s_1-s_2\|_2)$ for all $s_1,s_2\in \Omega$, where $K_0: [0,+\infty)\to[0,+\infty)$.
For all $r>0$, define the function $v_r(\cdot):\mathbb{R}\to\mathbb{R}_+$ as $v_r(x) = \sum_{k=0}^{+\infty} |x|^k/(k!)^r$. 
Then $1/v_r(x)$ is a monotone decreasing function of $x$, where a larger $r$ results in a slower decay rate in $1/v_r(x)$. Specifically, $1/v_1(x) = \exp(-x)$ and $1/v_r(x)$ with $r>1$ has a decay rate faster than all polynomials but slower than the exponential function.
We also define a series of polynomials $c_r(x)=1+x^r$ for $r>0$. For each $v_r$ with $r>1$ and each $c_r$ with $r>0$, we define the following families of Gaussian processes:
\begin{align}
\mathscr{Z}_{v_r} & = \left\{Z=(Z_s: s\in\Omega): K_0(\|s_1-s_2\|_2) \le \frac{1}{v_r(\|s_1-s_2\|_2) (1+\|s_1-s_2\|_2^{d+1})} \right\}, \label{eq GP family} \\
\mathscr{Z}_{c_r} & = \left\{Z=(Z_s: s\in\Omega): K_0(\|s_1-s_2\|_2) \le \frac{1}{(1+\|s_1-s_2\|_2)^r} \right\}. \label{eq GP polynomial family}
\end{align}
Intuitively, $\mathscr{Z}_{v_r}$ is the family of isotropic Gaussian processes whose covariance function decays no slower than some subexponential rate of the spatial distance, while $\mathscr{Z}_{c_r}$ is the family with covariance functions decaying no slower than an $r$th order polynomial of the spatial distance.

The objective of our theoretical study is to bound the Wasserstein-2 distance between the marginal distribution of RadGP and original Gaussian process on $\D$ when the original process belongs to the two spatial decaying families in \eqref{eq GP family} and \eqref{eq GP polynomial family}:
\begin{equation}\label{eq minimax}
\sup_{Z\in\mathscr{Z}_{v_r}}
W^2_2(Z_{\D},\hat{Z}_{\D})\quad \text{and}\quad \sup_{Z\in\mathscr{Z}_{c_r}}
W^2_2(Z_{\D},\hat{Z}_{\D}).
\end{equation}

\subsection{Rate of Approximation}
For the finite set $\D$ described above, we define the minimal separation distance among all spatial locations in $\D$ as $q = \min_{1\le i < j\le |\mathcal{D}|} \|w_i-w_j\|_2$. 
This minimal separation distance is a parsimonious statistic for describing the spatial dispersion of $\D$, and is useful in bounding various quantities such as the maximal eigenvalue and condition number of $\Sigma_{\D\D}$; see Lemma S2 of the Supplementary Material for details. Let $\hat{K}_0(w)=(2\pi)^{-d/2}\int K_0(x) \exp(-\imath w^{\Tr} x) dx$ be the Fourier transform of $K_0$ where $\imath^2=-1$, and let $\phi_0(x)=\inf_{\|w\|_2\le 2x} \hat{K}_0(w)$ for all $x>0$. For two positive sequences $a_n$ and $b_n$, we use $a_n \lesssim b_n$ to denote the relation that $a_n \le C b_n$ for all $n$ and some finite constant $C$. We have the following theorems on the Wasserstein-2 distance between the RadGP and the original Gaussian process. The constants $c_1, c_2$ are from Lemma S2 while $c_3$ is from the proof of Lemma S5 in the Supplementary Material, all of which only depend on the dimension $d$ of the spatial domain $\Omega$. 

\begin{theorem}\label{W22 rate}
For all radial neighbors Gaussian processes $\hat{Z}$, for the family $\mathscr{Z}_{v_r}$ in \eqref{eq GP family} with $r>1$, if $0<q<1$ and $n^{1/2} \{v_r(\rho d^{-1/2})\}^{-1}\{\phi_0(c_2/q)\}^{-9/2}  v_{r-1}[c_3\{\phi_0(c_2/q)\}^{-1}]  \le c_6$ hold for some constant $c_6$ only dependent on $d$, then 
\begin{equation}\label{eq:W22 rate}
\sup_{Z\in\mathscr{Z}_{v_r}} W_2^2(Z_{\D},\hat{Z}_{\D})  \lesssim  \frac{n}{v_r(\rho d^{-1/2})}\{\phi_0(c_2/q)\}^{-5} q^{-d} v_{r-1}[c_3\{\phi_0(c_2/q)\}^{-1}]. 
\end{equation}
Else if $q\ge 1$ and
$n^{1/2}/v_r(\rho d^{-1/2})\le c'_6$ hold
for some constant $c'_6$ only dependent on $d$, then
$\sup_{Z\in\mathscr{Z}_{v_r}} W^2_2(Z_{\D},\hat{Z}_{\D})  \lesssim n/v_r(\rho d^{-1/2}).$
\end{theorem}

\begin{theorem}\label{W22 rate poly}
For all radial neighbors Gaussian processes $\hat{Z}$, for the family $\mathscr{Z}_{c_r}$ in \eqref{eq GP polynomial family} with $r\ge d+1$, if $0<q<1$ and $n^{1/2}(1+\rho d^{-1/2})^{-(r-d-1)} q^{(r-7)d} \{\phi_0(c_2/q)\}^{-(r+4)}$ $ (c_1c_5d2^{d-1}\pi/\surd{6})^r \le c_7$ hold for some constant $c_7$ only dependent on $d$, then 
\begin{equation}\label{eq:W22 rate poly} \sup_{Z\in\mathscr{Z}_{c_r}} W_2^2(Z_{\D},\hat{Z}_{\D})  \lesssim \frac{n}{ (1+\rho d^{-1/2})^{r-d-1}} q^{(r-8)d} \{\phi_0(c_2/q)\}^{-(r+9/2)} (c_1c_5d2^{d-1}\pi/\surd{6})^r , \end{equation}
Else if $q\ge 1$ and
$ n^{1/2}(1+\rho d^{-1/2})^{-(r-d-1)} \{\phi_0(c_2/q)\}^{-r} (c_1c_5d2^{d-1}\pi/\surd{6})^r \le c'_7$ hold 
for some constant $c'_7$ only dependent on $d$, then  $\sup_{Z\in\mathscr{Z}_{c_r}} W^2_2(Z_{\D},\hat{Z}_{\D})  \lesssim $ 
$n (1+\rho d^{-1/2})^{-(r-d-1)}$.
\end{theorem}

The conditions for the above theorems, namely the equations regarding $c_6, c'_6, c_7, c'_7$, are necessary but not sufficient conditions for the upper bounds of Wasserstein-2 distance to go to zero. Thus, they can be safely ignored if we aim for accurate approximations. We now give a detailed analysis of these upper bounds for Wasserstein-2 distance 
involving three variables: the sample size $n$, the minimal separation distance $q$, and the approximation radius $\rho$ in the RadGP. We will discuss from the perspectives of both finite samples and asymptotics. From a finite sample perspective, we fix two out of three variables and discuss how and why the upper bounds change as the remaining variable changes; from an asymptotic perspective, we describe the relations between $n, q$ and $\rho$ such that as $n$ goes to infinity, the Wasserstein-2 distance goes to zero.

Viewing Theorems \ref{W22 rate} and \ref{W22 rate poly} as finite sample bounds, 
the upper bounds  are linearly dependent on the sample size $n$. This directly follows  
from the definition of Wasserstein distance which measures the difference between two $n$-dimensional Gaussian distributions. Provided both the minimal separation distance $q$ and approximation radius $\rho$ are fixed, as the sample size $n$ increases, the local structure for existing locations will not change and the approximation accuracy for them will not improve. Thus, adding more locations will lead to a linearly increasing approximation error in the squared $W_2$ distance.

Second, the upper bounds on the squared $W_2$ distance decay with respect to the approximation radius $\rho d^{-1/2}$ at the same rate as the spatial covariance function with respect to the spatial distance, up to a $(d+1)$th order polynomial. For Theorem \ref{W22 rate}, the covariance function is required to decay faster than $1/v_r(x)$ by a polynomial $1+x^{d+1}$ as defined in equation (\ref{eq GP family}); for Theorem \ref{W22 rate poly}, the covariance function decays faster than an $r$th order polynomial while the upper bound of $W_2$ distance decays at an $(r-d-1)$th order polynomial.
The difference between the decay rate of the covariance function and the $W_2$ distance in Theorems \ref{W22 rate} and \ref{W22 rate poly} comes from the fact that high dimension reduces the effective decay rate of covariance function. For example, under a fixed minimal separation distance, when $d=1$, as long as the covariance function decays faster than $1/x$ by a polynomial term, the matrix $\Sigma_{\D\D}$ has a bounded $l_1$ norm regardless of sample size; however, for a general $d$ dimensional space, the covariance function must decay faster than $1/x^d$ for  $\Sigma_{\D\D}$ to have a bounded $l_1$ norm.
These upper bounds show that the approximation radius $\rho$ indeed controls the approximation accuracy. Specifically, for fixed $q$ and $n$, we can reach an arbitrarily small $W_2$ error by increasing $\rho$ over some threshold value; see our detailed conditions on $\rho$ in Corollary \ref{corollary} for various types of covariance functions.

Third, for fixed $\rho$ and $n$, both upper bounds in the two theorems increase as the minimal separation distance $q$ decreases.
This is because as the minimal separation distance decreases, spatial locations get closer to each other and the covariance matrix gets close to singular, 
leading to a larger approximation error.

From the asymptotic perspective, we can increase the sample size $n$ while viewing the minimal separation distance $q$ and the approximation radius $\rho$ as some functions of $n$. The objective is to choose an approximation radius $\rho$ dependent on $n$ and $q$ such that the squared Wasserstein distance goes to zero as $n$ goes to infinity. 
We first present a corollary showing how to choose such $\rho$ for three commmonly used covariance functions, where the convergence is understood in the sense that $q$ is also a function of $n$.

\begin{corollary}\label{corollary}
We have $W_2^2(Z_{\D},\hat{Z}_{\D}) \to 0$ as $n\to\infty$ if one of the following conditions is satisfied:\\
\noindent (1) The covariance function is the isotropic Mat\'{e}rn: $K_0(s_1-s_2) = \tau^2 2^{1-\nu} {\Gamma(\nu)}^{-1} \left(\phi\|s_1-s_2\|\right)^{\nu} \mathcal{K}_{\nu}\left(\phi\|s_1-s_2\|_2\right)$ with $q<1/\phi$.
Define the constant $c_{m,1} = \Gamma(\nu) / \{\tau^2 2^d \pi^{d/2}\Gamma(\nu+d/2)\}$. The approximation radius $\rho$ satisfies
$$\rho \ge  \frac{d^{1/2}}{\phi}\left[c_3c_{m,1} \Big(1 + \frac{ 4c_2^2}{\phi^2 q^2}\Big)^{\nu+\frac{d}{2}} + \ln\Big\{ c_{m,1}n q^{-d} \Big(1 + \frac{ 4c_2^2}{\phi^2q^2}\Big)^{5(\nu+\frac{d}{2})} \Big\} \right]^3.$$
\noindent (2) The covariance function is Gaussian $K_0(\|s_1-s_2\|_2) = \tau^2 \exp(-a\|s_1-s_2\|_2^2)$ with $q<a^{-1/2}$. The approximation radius $\rho$ satisfies
$$\rho \ge \left(\frac{d}{a}\right)^{1/2} \left[ c_3  \tau^{-2} 
\exp\big\{c_2^2/(aq^2)\big\} +\ln( n  q^{-d} \tau^{-10}) + \frac{5c_2^2}{aq^2} \right]^3.$$
\noindent (3) The covariance function is the generalized Cauchy $K_0(\|s_1-s_2\|_2) = \tau^2 \big\{1+(\phi\|s_1-s_2\|)^\delta\big\}^{-\lambda/\delta}$ with parameter $\lambda>d+1$ and $q<\phi$. Let $c_9$ be a constant dependent on dimension $d$ and covariance function parameters $\tau^2, \phi, \delta, \lambda$. The approximation radius $\rho$ satisfies
$$\rho \ge c_9 q^{-\frac{\frac{25}{2} d+\delta(\lambda+\frac{9}{2})}{\lambda-(d+1)}}n^{\frac{1}{\lambda-(d+1)}}.$$
\end{corollary}

We then analyze our main theorems from the asymptotic perspective and verify the analysis using the three examples in Corollary \ref{corollary}. 
Depending on whether and how the minimal separation distance $q$ changes with respect to the sample size $n$, there are three commonly used asymptotic frameworks:  fixed domain, mixed domain and increasing domain asymptotics.

Under the increasing domain setting, 
as $n$ goes to infinity, the radius of the domain $\Omega$ increases while the minimal separation distance $q$ is fixed.  
Therefore, as long as $n \gtrsim d^{1/2} v_r^{-1}(n)$ under conditions of Theorem \ref{W22 rate} or $n \gtrsim d^{1/2} n^{1/(r-d-1)}$ under conditions of Theorem \ref{W22 rate poly}, we have $o(1)$ approximation error. When applied to Mat\'{e}rn covariance in case (1) and Gaussian covariance in case (2) of Corollary \ref{corollary}, this yields a slowly increasing rate of $O(\ln^3 n)$ for the approximation radius $\rho$, which is much smaller than the sample size $n$. On the other hand, when the covariance function decays polynomially such as the  generalized Cauchy in case (3) of Corollary \ref{corollary}, the approximation radius $\rho$ needs to increase at a polynomial rate of sample size $n$. This is because for covariance functions that decrease slowly in the spatial distance, one needs more spatial locations and a larger approximation radius to capture the dependence information.

Under the mixed domain setting, as $n$ goes to infinity, the radius of the domain $\Omega$ increases with $n$ while the minimal separation distance $q$ also decreases with $n$.
By solving the inequalities on the right hand side of equations (\ref{eq:W22 rate}) and (\ref{eq:W22 rate poly}) being no greater than $o(1)$, we obtain the lower bound of $\rho$ that guarantees the $o(1)$ approximation error. Corollary \ref{corollary} directly shows such lower bounds for Mat\'{e}rn, Gaussian and generalize Cauchy covariance functions. We can still achieve the $o(1)$ approximation error as long as these lower bounds of $\rho$ in the three cases of Corollary \ref{corollary} increase slower than the radius of the domain $\Omega$. Otherwise, if $\rho$ increases exactly at the same rate as the radius of $\Omega$, RadGP is identical to the original Gaussian process and does not make any approximation at all. Our theory is not applicable when the lower bounds of $\rho$ exceed the radius of domain $\Omega$. In particular, for the fixed domain setting, the radius of $\Omega$ is fixed as a constant and the minimal separation distance $q$ decreases to zero at a rate no slower than $O(n^{-1/d})$. In such case, the lower bounds of $\rho$ in Corollary \ref{corollary} will diverge with $n$, and therefore Theorems \ref{W22 rate} and \ref{W22 rate poly} do not provide a vanishing bound for Wasserstein-2 distances.

Our results fill a major gap in the literature by providing the previously lacking theoretical support for popular local or neighborhood-based Gaussian process approximations. 
Our theory only requires knowledge of the spatial decaying pattern of the covariance functions (e.g., exponential and polynomial in $\mathscr{Z}_{v_r}$ and $\mathscr{Z}_{c_r}$, respectively). These conditions  are easily verifiable, with derivations for the Mat\'ern, Gaussian and the generalized Cauchy covariance functions provided in Corollary \ref{corollary}. Similar techniques can be used for a wider range of stationary covariance functions. Furthermore, our novel theory explicitly links the approximation quality from choosing radius $\rho$ with the minimal separation distance $q$ and the sample size $n$, both of which are readily available from observed data.

\section{Bayesian Regression with RadGP}\label{sec:regression}
Consider a linear regression model with spatial latent effects as:
\begin{equation}\label{eq model}
Y(s) = X(s)\beta + Z(s) + \epsilon(s), \; s\in\Omega,
\end{equation}
where $X(s)$ are covariates at location $s$, $\beta$ are regression coefficients for the covariates, $Z(s)$ is a Gaussian process with zero mean and covariance function $K_\theta: \Omega\times\Omega\to \mathbb{R}$, and  
$\epsilon(s)$ is a white noise process.
Researchers observe $Y(s),X(s)$ at a collection of training locations $\mathcal{T}_1$ with the objective of estimating the regression coefficients $\beta$, the covariance function parameters $\theta$, and the spatial random effects $Z(s)$ at both training locations $\mathcal{T}_1$ and test locations  $\mathcal{T}_2\subset\Omega$. If the spatial random effects $Z(s)$ are not of interest, one can combine $Z(s)$ and $\epsilon(s)$ into one Gaussian process, known as the response model. This is described in Section S4.2 of the Supplementary Material. Hereafter we assume that inference on $Z(s)$ is desired. Let $Z_{\mathcal{T}_1}$ denote the column vector consisting of $(Z(s):s\in\mathcal{T}_1)$. Similarly define other notations with $\mathcal{T}_i$ as subscripts.
Endow the spatial random effects on $\T_1\cup\T_2$ with our RadGP prior defined in Section \ref{subsec:construction}, which is built on the radial neighbors graph constructed in Section \ref{subsec:dag practice} such that the locations in the training set $\T_1$ precede locations in the test set $\T_2$. The full Bayesian model can be outlined as follows:
\begin{align}
& (\beta,\theta,\sigma^2)\sim p(\beta)p(\theta)p(\sigma^2), \quad \epsilon(s) \overset{\text{i.i.d.}}{\sim} N(0, \sigma^2), s\in\Omega,\quad 
\hat{Z}_{\mathcal{T}_1}|\theta \sim 
\prod_{s\in\T_1} p(\hat{Z}_s|\hat{Z}_{\mr{pa}(s)}), \nonumber\\ 
& p(\hat{Z}_{\mathcal{T}_2}|\hat{Z}_{\mathcal{T}_1},\theta) \propto \prod_{s\in\mathcal{T}_2} p(\hat{Z}_s|\hat{Z}_{\mr{pa}(s)}), \quad
Y(s) = X(s)\beta + \hat{Z}(s) + \epsilon(s), \quad \text{for } s\in\Omega\label{eq ZT}.
\end{align}
The priors for $\beta, \theta$ and $\sigma^2$ can be flexible, but it is advised to set a proper prior for either $\sigma^2$ or $\theta$ due to potential non-identifiability.

We employ a Gibbs sampling framework while handling each full conditional with different approaches.
If the prior of $\beta$ is normal $\beta\sim N(\beta_0, \Phi_0^{-1})$, then the posterior of $\beta$ is also normal:
\begin{align}\label{eq beta post} & \beta|{Y}_{\T_1},\hat Z_{\T_1},\sigma^2 \nonumber \\
\sim & N\big((\Phi_0+{X}_{\T_1}^T{X}_{\T_1}/\sigma^2)^{-1}\{\Phi_0\beta_0+{X}_{\T_1}^T({Y}_{\T_1}-\hat Z_{\T_1})/\sigma^2\}, (\Phi_0+{X}_{\T_1}^T{X}_{\T_1}/\sigma^2)^{-1}\big) .\end{align}
We then sample $\sigma^2$ from its full conditional. If the prior distribution for $\sigma^2$ is inverse gamma $\sigma^2 \sim \mr{IG}(a_0, b_0)$, then the posterior of $\sigma^2$ is also inverse gamma:
\begin{equation}\label{eq sigma post} \sigma^2|{Y}_{\T_1}, \hat Z_{\T_1},\beta \sim \mr{IG}\left(a_0+n/2, b_0+\|{Y}_{\T_1}-{X}_{\T_1}\beta-\hat Z_{\T_1}\|_2^2/2\right).\end{equation}
Denote the precision matrix of $\hat{Z}_{\T_1}|\theta$ from the RadGP prior as $\hat{\Phi}$. We sample the spatial random effects on $\T_1$ with conjugate gradients as illustrated in \cite{nishimura2022prior},
using the full conditional distribution:
$$\hat Z_{\T_1}|{Y}_{\T_1},\beta,\sigma^2,\theta \sim N({\xi},(\hat{\Phi}+{I}_n/\sigma^2)^{-1}), $$
where ${\xi} = (\hat{{\Phi}}+{I}_n/\sigma^2)^{-1}({Y}_{\T_1}-{X}_{\T_1}\beta)/\sigma^2$.

Finally, for the parameter $\theta$ in the covariance function $K_{\theta}$, there is no available conjugate prior in general. We sample $\theta$ using the Metropolis-Hastings algorithm relying on the full conditional $p(\theta|\hat Z_{\T_1}) \propto p(\hat Z_{\T_1}|\theta)p(\theta)$. 
A detailed algorithm for posterior sampling can be found in Section S4.1 of the Supplementary Material.

If prediction on a test set is of interest, given MCMC samples of random effects on training set $\T_1$ and parameter $\theta$, $Z_{\T_2}$ can be sampled from each unidimensional Gaussian distribution as in equation (\ref{eq ZT}). Unlike most previous related work, this prediction accounts for dependence among test locations, which is important in many applied contexts.

In practical applications where the sample sizes of training set and test set are similar and much larger than the cardinality of parent sets, the computational complexity is $O( M l n_{\text{cg}} n )$, where $l$ is the number of MCMC iterations, $n_{\text{cg}}$ is the average number of conjugate gradient iterations, $n$ is the sample size of training set, and $M$ is the average number of nonzero elements per column in the RadGP precision matrix $\hat{\Phi}$. The derivation of computational complexity is in Section 
S4 of the Supplementary Material; here we discuss implications. The variable $n_{\text{cg}}$ is determined by the distribution of eigenvalues of matrix $\hat{\Phi}+I_n/\sigma^2$ and numerical tolerance. The latter is fixed as $10^{-6}$ in our R package while the former is determined by the data set if the RadGP precision matrix $\hat{\Phi}$  approximates the original Gaussian process precision matrix $\Phi$ decently well. 
Therefore, the computational complexity is roughly proportional to $M$, the average number of nonzero elements per column in the precision matrix.

\section{Experimental Studies}\label{sec:experiments}
\subsection{Approximation of a Gaussian Process Prior}\label{subsec:exp prior}
\begin{figure}[h]
    \centering
    \includegraphics[width=\textwidth]{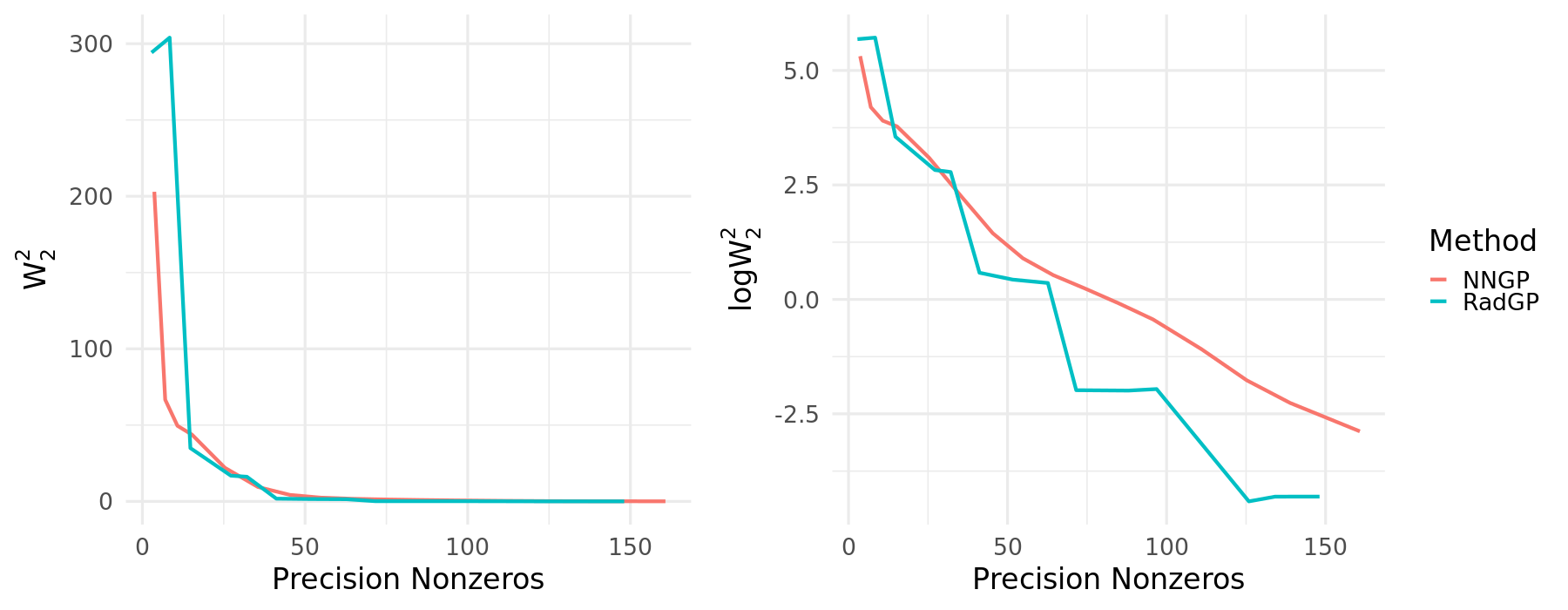}
    \caption{Accuracy of approximating a Gaussian process prior. X-axis is the average number of nonzero elements per column in the precision matrix. Y-axis is the squared Wasserstein-2 distance and its logarithm.}\label{W2}
\end{figure}

Throughout Section 5, we consider the specific RadGP described in Section 2.
We compare the RadGP prior with the nearest-neighbor Gaussian process (NNGP, \citealt{Datetal16a}) implemented on a maximin ordering of the data \citep{guinness2018permutation}, in terms of their approximation accuracy. 
The true Gaussian process has mean zero and an isotropic Mat\'ern covariance function
\begin{equation}\label{eq:exp}
K_0(x) = \frac{\tau^2 2^{1-\nu}}{\Gamma(\nu)} (\phi x)^{\nu} \mathcal{K}_{\nu} (\phi x) ,  
\end{equation} 
with the true parameters $(\phi, \tau^2, \nu) = (31.63, 1, 1.5)$. The value $\phi=31.63$ is obtained by setting $K_0(0.15)=0.05$, a choice suggested by \cite{Katetal20}.
The experiment is carried out on a $40\times 40$ grid in $[0,1]^2$. 
According to the analysis of computational complexity in Section \ref{sec:regression}, we use the average number of nonzero elements per column in the precision matrix as a measure of complexity for both methods. The quality of approximation is measured by squared Wasserstein-2 distance between each approximate process and the true Gaussian process, which can be computed efficiently using the R package \texttt{transport} available from \texttt{CRAN} (\citealt{transport24}). As shown in Figure \ref{W2}, the results are generally in favor of RadGP, as it achieves smaller approximation error under various complexity of the directed graphs.

\subsection{Approximation of a Gaussian Process Regression Posterior}\label{subsec:posterior}
We study the performance of RadGP for the spatial regression model (\ref{eq model}) with no covariates, where the true covariance function for spatial effects $Z(s)$ is the same as that in Section \ref{subsec:exp prior}. The white noise $\epsilon(s)$ follows $N(0,\sigma^2)$ with $\sigma^2=0.01$. The parameters $\phi, \tau^2, \sigma^2$ are unknown and need to be estimated from training data, while the Mat\'{e}rn smoothness parameter $\nu=3/2$ is assumed to be fixed and known.
The domain is set as $\Omega = [0,1]^2$, with training data consisting of $1600$ grid locations in $\Omega$. We replicate the analysis on 50 datasets. For each dataset, the test set is generated as $1000$ random samples from the uniform distribution on $\Omega$. 

Three methods are compared: RadGP proposed in this paper, nearest-neighbor Gaussian process  with maximin ordering, and Vecchia Gaussian process predictions (V-Pred) \citep{Katetal20} with maximin ordering. 
The last two methods share the same training procedure but differ in how they make predictions: NNGP 
assumes that testing locations are independent given the training locations whereas V-Pred 
allows dependencies among testing locations.
All methods are fitted via the MCMC algorithm described in Section \ref{sec:regression} and thus only differ in terms of their assumed graphs.

\begin{figure}[h]
    \centering
    \includegraphics[width=0.75\textwidth]{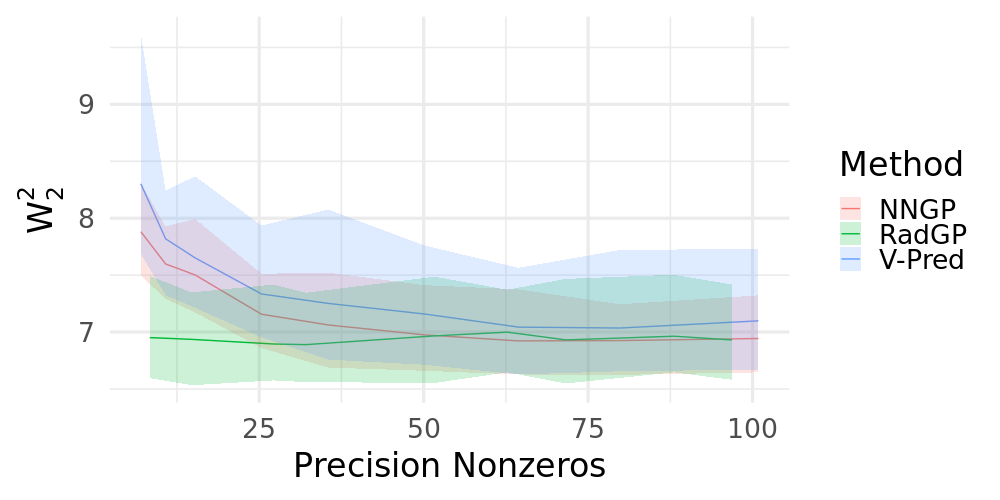}
    \caption{Accuracy of approximating the posterior predictive distribution in Gaussian process regression. X-axis is the average number of nonzero elements per column in the precision matrix, and Y-axis is the squared Wasserstein-2 distance. Solid lines are mean values and shaded regions are $90\%$ Bayesian credible intervals.
    }\label{fig:posterior}
\end{figure}

To measure the accuracy of approximations for prediction in the spatial regression model (\ref{eq model}), we compute the Wasserstein-2 distance between the posterior predictive distributions from each of the three approximation methods and the true Gaussian process regression model. The Wasserstein-2 distance is computed using 
MCMC samples.
The results are shown in Figure \ref{fig:posterior}. When the complexity of the directed graph is sufficiently large, NNGP and RadGP tend to have the same approximation error while V-Pred is slightly worse. However, RadGP is able to achieve the limiting approximation accuracy with a much smaller complexity compared to the other two methods, indicating a faster approximation by RadGP. Unlike prior approximation in Section \ref{subsec:exp prior}, the posterior approximation error does not go to zero when the complexity of directed graphs increases. This nonzero error comes from two parts: first, the Gaussian process regression model (\ref{eq model}) has unknown parameters that need to be estimated from a finite amount of data; second, the Wasserstein-2 distance computed based on a finite number of MCMC samples contains Monte Carlo error.

\subsection{Surface Temperature Data}
\begin{figure}[h]
    \centering
    \includegraphics[width=\textwidth]{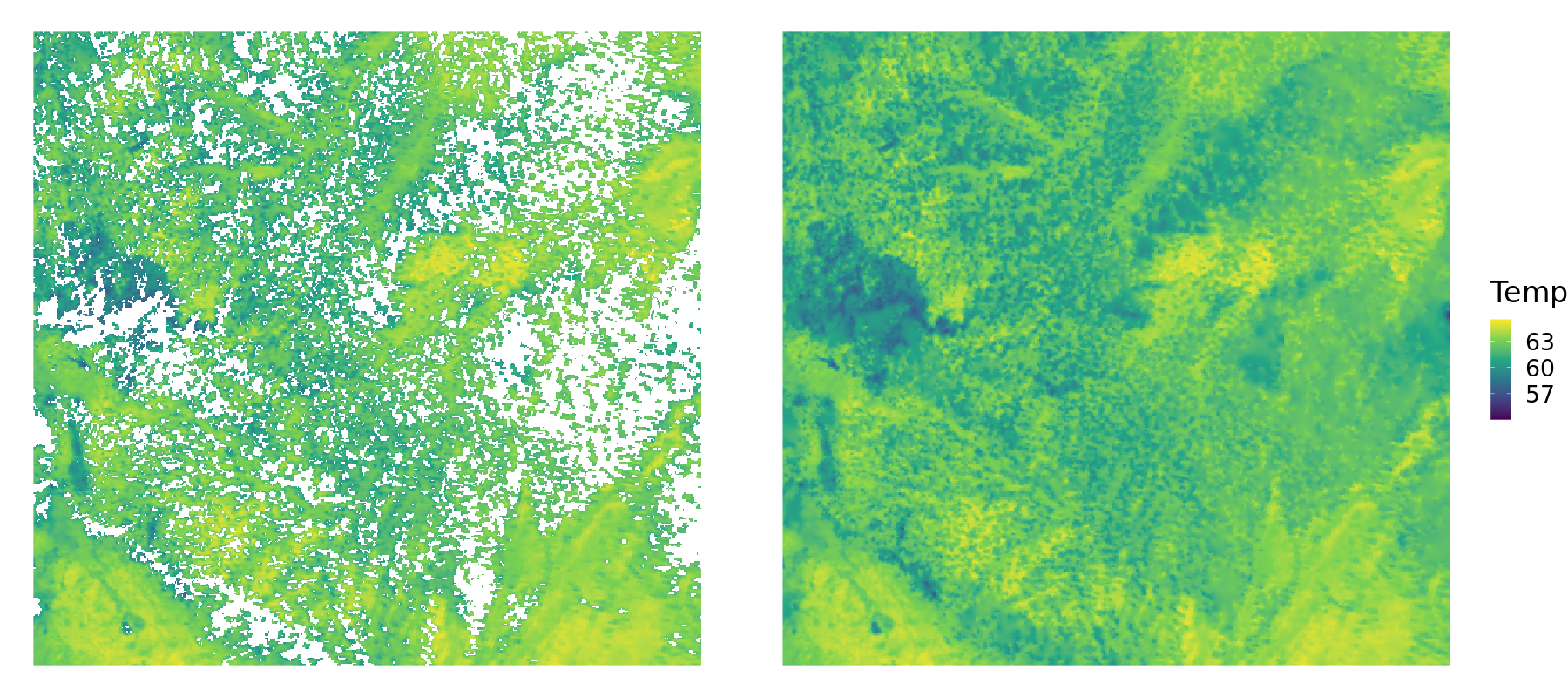}
    \caption{Ethiopian land surface temperature on April 1st, 2020. Left: observed data; Right: predicted values from RadGP.
    }\label{fig LST}
\end{figure}

We consider the land surface temperature data on a $499\mr{km}\times 581\mr{km}$ region southeast of Addis Ababa in Ethiopia at April 1st, 2020. The data are a product of the MODIS (Moderate Resolution Imaging Spectroradiometer) instrument on the Aqua satellite. The original data set comes as a $499\times 581$ image with each pixel corresponding to a $1\mr{km}\times1\mr{km}$ spatial region, as shown in Fig. \ref{fig LST} left. However, atmosphere conditions like clouds can block earth surfaces, resulting in missing data for many pixels. There are $211,683$ spatial locations with observed temperature values and $78,236$ locations with missing values. We randomly split the training data into $3$ folds, using two folds for training and one fold for out-of-sample validation. This results in a total of $141,122$ training locations and $148,797$ testing locations.

We consider the same Gaussian process regression model as Section 5.1 with no covariates, zero mean function and exponential covariance function $K_0(x) = \tau^2\exp(-\phi x)$. The observed temperature values are centered to compensate for the lack of intercept term. We evaluate the performances of RadGP with radius $4.01\mr{km}$ and nearest neighbor Gaussian process with neighbor size $15$. The priors for both methods are set as $p(\phi)\propto 1$, $\tau^2\sim \mr{IG}(2,1)$ and $\sigma^2\sim\mr{IG}(1,0.1)$.  Estimates for covariance parameters, out-of-sample mean squared errors and coverage for $95\%$ credible intervals are shown in Table \ref{tab 2}. The results are quite similar between the two methods.

The observed temperature and predicted temperature by RadGP are visualized in Fig. \ref{fig LST}. In general, the land surface temperature changes smoothly with respect to geological locations. The east and south part of the figure, which is covered by less vegetation, features slightly higher surface temperature. Small variations of temperature also occur in many local regions, indicating the geographical complexity of Ethiopia. There is a small blue plate in the west border of the figure that has significantly lower surface temperature than its surroundings, which corresponds to several lakes in the Great Rift Valley.

We further compare the performance of the two methods in terms of joint predictions at multiple held out locations. As we do not know the ground truth, we cannot calculate Wasserstein-2 error and instead focus on the joint likelihood ratio for held out data. Specifically, we select $567$ pairs of locations in the set reserved for out-of-sample validation. Each pair consists of two locations with $1$km distance from each other. Any two locations from different pairs are at least $10$km apart. For each pair of locations, we compute the logarithm of the likelihood ratio between RadGP and NNGP. The distribution of the log likelihood ratio can provide important information about joint inference. Specifically, the mean of the log likelihood ratio is $2.152$ and there are $68\%$ pairs with positive log likelihood ratios. Since the distribution of the log likelihood ratio has heavy tails, we remove the bottom and top $5\%$ of the values and visualize the rest in Fig. \ref{fig loglike}, where the red vertical line is $\mr{ratio}=0$. We can see in the majority of cases, RadGP performs better than NNGP in terms of characterizing dependencies between temperature at two close locations.

\begin{table}[h]
\begin{minipage}[b]{0.54\linewidth}
\small
    \centering
    \begin{tabular}{c|cc}
     & RadGP & NNGP \\ \hline
    $\phi\times 10^2$ & 5.761 (5.432, 6.133) & 5.689 (5.020, 6.321) \\
    $\tau^2$ & 1.706 (1.605, 1.807) & 1.744 (1.565, 1.972) \\
    $\sigma^2\times10^3$ & 1.222 (1.104, 1.351) & 1.105 (0.994, 1.222) \\
    MSE & 0.0761 & 0.0770  \\
    coverage & 0.954 & 0.953\\ \hline
    time & 1929 sec & 1811 sec\\ \hline
    \end{tabular}
    \caption{Performance of RadGP and NNGP on \\land surface temperature data.}\label{tab 2}
\end{minipage}
\begin{minipage}[b]{0.44\linewidth}
\includegraphics[width=\linewidth]{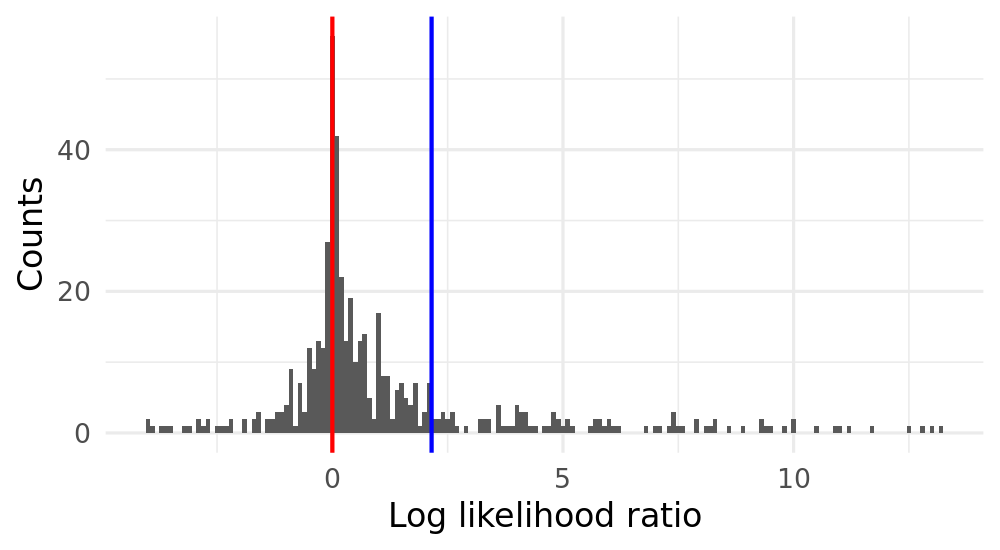}
\captionof{figure}{Histogram of log likelihood ratio between RadGP and NNGP. The blue and red lines correspond to the mean log likelihood ratio value, and zero, respectively.}
\label{fig loglike}
\end{minipage}
\end{table}

\section{Discussion}

Vecchia approximations are arguably one of the most popular approaches for scaling up Gaussian process models to large datasets in spatial statistics. The major contribution of this article is to propose the first Vecchia-like approximation method that has strong theoretical guarantees in quantifying the accuracy of approximation to the original process. Our numerical experiments have shown competitive performance from the proposed RadGP compared to other Vecchia approximation methods. We have demonstrated that our method leads to improvements in characterizing joint predictive distributions while matching the performance of state-of-the-art methods in parameter estimation and marginal predictions.

Our current theory leads to several potential future research directions. First, we have established upper bounds for the approximation error to the original Gaussian process in Wasserstein-2 distance. We can further study  lower bounds of this approximation error within specific families of approximations. It will be helpful to examine when the norm-controlled inversion bounds we have used in the proofs are tight, and whether the derived conditions on the approximation radius $\rho$ for various covariance functions in Corollary \ref{corollary} can be made optimal for such families. 

Second, existing Vecchia approximations that use maximin ordering \citep{Katetal20} systematically generate remote locations across the domain. While our theory shows close locations are effective in quantifying spatial dependence, it is still unclear whether and how remote locations help with the approximation of covariance function. We would like to expand our research into general Vecchia approximations. The objective is to study the approximation error for different approximation approaches and discover the optimal strategy that minimizes the approximation error under certain sparsity constraints.

Third, although our main focus has been on approximation of Gaussian processes, we can further consider how the approximate process will affect posterior estimation of model parameters, including spatial fixed effects and covariance parameters. Under our current theoretical framework, since the minimal separation distance $q$ needs to be lower bounded, we can study  Bayesian parameter estimation under increasing domain \citep{MarMar84} or mixed domain \citep{LahZhu06} asymptotic regimes. It remains an open problem whether our scalable approximation can maintain posterior consistency and optimal rates for all model parameters. We leave these directions for future exploration.

\section*{Acknowledgements}
The authors thank Dr. Karlheinz Gr\"{o}chenig for helpful discussion regarding the theory on norm-controlled inversion of Banach algebras, and two anonymous reviewers for detailed feedback and suggestion that greatly improved this article.
Yichen Zhu has received funding from the European Research Council (ERC) under grant No. 101041064. 
Yichen Zhu, Michele Peruzzi and David Dunson have received funding from ERC under grant No. 856506 and United States National Institutes of Health under grant No. R01ES028804.
Cheng Li has received funding from Singapore Ministry of Education Academic Research Funds Tier 1 Grant A-0004822-00-00. 
Part of this research was performed while Michele Peruzzi was visiting the Institute for Mathematical and Statistical Innovation (IMSI), which is supported by the National Science Foundation (grant No. DMS-1929348).

\vspace{-2mm}

\section*{Supplementary Material}
Proofs of lemmas, theorems and corollaries, detailed algorithms, computational complexity analysis and additional experimental studies are available at {\it Supplementary Material for ``Radial Neighbors for Provably
Accurate Scalable Approximations of Gaussian Processes.''} An R package for regression models using radial neighbors Gaussian processes is available at \url{https://github.com/mkln/radgp}.

\bibliographystyle{biometrika}
\bibliography{AGP}

\end{document}


\nolinenumbers
\jname{Submitted to Biometrika}


\markboth{Y. Zhu et~al.}{Radial Neighbors Gaussian Process}

\title{Radial Neighbors for Provably Accurate Scalable Approximations of Gaussian Processes}

\author{Yichen Zhu}
\affil{Bocconi Institute for Data Science and Analytics, Bocconi University, Milan, MI, Italy.
\email{yichen.zhu@unibocconi.it}}

\author{Michele Peruzzi}
\affil{Department of Biostatistics, University of Michigan, Ann Arbor, Michigan, U.S.A. \email{peruzzi@umich.edu}}

\author{Cheng Li}
\affil{Department of Statistics and Data Science, National University of Singapore, Sinagpore.
\email{stalic@nus.edu.sg}}

\author{\and David B. Dunson}
\affil{Department of Statistical Science \& Mathematics, Duke University, North Carolina, U.S.A. \email{dunson@duke.edu}}

\maketitle

This supplementary material is organized as follows. Section \ref{sec:proof.lemma.sec2} includes the proof of lemma \ref{valid process} in the main paper. Section \ref{sec:aux.lemma} provides several auxiliary lemmas for our main theorems. The subsequent Section \ref{sec:theory proof} states the proofs of Theorems 1 and 2 and Corollary 1 in the main paper.
Section \ref{sec:algo.radgp} provides the posterior sampling algorithms and computational complexity analysis for the RadGP regression model described in Section \ref{sec:regression} of the main paper. Section \ref{sec:add exp} provides experimental studies that are not included in the main paper due to length limit. Section \ref{sec:taper} provides further discussion on the relation between our RadGP theory and the existing tapering-
based Gaussian process theory.

We define some notation that will be used throughout this supplementary material. The notation that is only used for specific lemmas will be defined before or in those lemmas. Let $\Omega$ be the spatial domain that is a connected subspace of $\mathbb{R}^d$. Let $Z = (Z_s: s\in\Omega)$ be a real valued Gaussian process on $\Omega$ with zero mean and covariance function $K:\Omega\times\Omega\to\mathbb{R}$, and let $\hat Z$ be the radial neighborhood Gaussian process (RadGP) on $\Omega$. 
Let $\D\subset\Omega$ be an arbitrary finite subset of $\Omega$. Without loss of generality, we index the elements of $\D$ as $\D = \{w_1, w_2, \ldots, w_n\}$, such that for any $w_i\in\D$, all locations in $\mr{pa}(w_i)$ precede $w_i$ in the order prescribed by the radial neighbors graph in Section \ref{subsec:dag practice} of the main text.
For any index $i=1,\ldots,n$, let $\nu(i) = \{j\in\mathbb{N}: j<i\}$ be the collection of indices that are smaller than $i$. Since $\Omega$ can be an uncountable set, we use $s$ to denote a generic spatial location of $\Omega$. In the context where a directed graph exists, we use $\mr{pa}(s)$ to denote the parent set of $s$. 

For a finite set $A\subset\Omega$, let $Z_A$ be the finite-dimensional random variable of the process $Z$ on $A$. For two finite sets $A,B\subset\Omega$, we denote the covariance matrix between $Z_A$ and $Z_B$ as $\Sigma_{AB}$. For a generic square matrix $C$, let $\lambda_{\min}(C)$ and $\lambda_{\max}(C)$ be its smallest and largest eigenvalues. The matrix $l_1$, $l_2$ and trace norms are denoted by $\|\cdot\|_1$, $\|\cdot\|_2$ and $\|\cdot\|_{\mr{tr}}$. The vector $l_1$ and $l_2$ norms are denoted using the same notation as the matrix norms. Finally, for simplicity, we denote $(\hat{\Sigma}_{\D\D})^{-1}$ as $\hat\Phi$.

\section{Proof of Lemma \ref{valid process}}
\label{sec:proof.lemma.sec2}

\begin{proof}
By the Kolmogorov extension theorem, $\hat{\bm{Z}}(\cdot)$ defines a stochastic process if its finite density $p(\hat{\bm{Z}}_A)$ on any finite set $A\subset \Omega$ defined in Equation (\ref{eq alt-dag 4}) of the main paper satisfies the following two conditions:
\vspace{2mm}

\noindent (C1) For all finite set $A\subset\Omega$, let $\tilde A$ be a permutation of elements in $A$, then $p(\hat{\bm{Z}}_{\tilde{A}}) = p(\hat{\bm{Z}}_{A})$; 
\vspace{2mm}

\noindent (C2) For all finite sets $A, A'\subset\Omega$ with $A\cap A'=\emptyset$, 
$\int p(\hat{\bm{Z}}_{A\cup A'}) \prod_{s\in A'} ds = p(\hat{\bm{Z}}_{A})$.
\vspace{2mm}

We verify these two conditions for the radial neighbors Gaussian process, defined in Equations (\ref{eq alt-dag 2}), (\ref{eq alt-dag 3}) and (\ref{eq alt-dag 4}) of the main text.
\vspace{2mm}

\underline{\textsc{Condition C1}} This essentially requires that for all $A\subset \Omega$, Equation (\ref{eq alt-dag 4}) in the main paper defines a valid finite dimensional distribution. In the following, we prove a stronger condition that for any finite set $A\subseteq \Omega$, the distribution on $A$ defined by Equation (\ref{eq alt-dag 4}) in the main paper is a Gaussian distribution. 

{\it Let $\bm{u}_1\in\mathbb{R}^{d_1}, \bm{u}_2\in\mathbb{R}^{d_2}$ be two random vectors. If $\bm{u}_1\in N(\bm\mu_1,\bm\Sigma_1)$, $\bm{u}_2|\bm{u}_1\sim N(\bm\mu_2 + \bm B(\bm{u}_2-\bm\mu_1), \bm\Sigma_2)$, then $(\bm{u}_2,\bm{u}_1)$ follows a joint Gaussian distribution.
}

This claim can be verified by direct calculation and its proof is omitted. We now show $\hat{\bm{Z}}_A$ follows a Gaussian distribution. Recall the notations in the main paper that $A = U\cup V$ with $U\subset\Omega\backslash\D$ and $V\subset\D$.
By the claim, $p(\hat{Z}_\D) = \prod_{i=1}^n p(\hat{Z}_{w_i}|\hat{Z}_{\mr{pa}(w_i)})$ follows a joint Gaussian distribution. Similarly, since for all $s\in U$, $\mr{pa}(s)\subset\D$ and $\hat{Z}_{\D\cup U} = p(\hat{Z}_\D) \prod_{s\in U} p(\hat{Z}_s|\mr{pa}(s))$,
we have $\hat{Z}_{\D\cup U}$ also follows a Gaussian distribution. Noticing that $\hat{Z}_A$ is can be obtained from $\hat{Z}_{\D\cup U}$ by integrating out $Z_{w_i}, \forall i \in \D\backslash A$,
we conclude $\hat{Z}_A$ also follows a Gaussian distribution.

\vspace{2mm}

\underline{\textsc{Condition C2}}
Since $A'$ is a finite set, it suffices to prove condition C2 when $A' = \{x\}$ for arbitrary finite set $A\subset \Omega$ and $x\in \Omega\backslash A$. 
We split the proof into two cases: $x\in \D$ and $x\in \Omega\backslash\D$.

For $x\in \D$, by definition, 
\begin{align*}
p(\hat{\bm{Z}}_{A}) 
= & \int p(\hat{\bm{Z}}_U|\hat{\bm{Z}}_{\D}) \prod_{i=1}^n p(\hat{Z}_{w_i}|\hat{Z}_{\mr{pa}(w_i)}) \prod_{w_i\in \D\backslash V} dw_i \\ 
= & \int\left\{\int p(\hat{\bm{Z}}_U|\hat{\bm{Z}}_{\D}) \prod_{i=1}^n p(\hat{Z}_{w_i}|\hat{Z}_{\mr{pa}(w_i)}) \prod_{w_i\in \D\backslash (V\cup\{x\})} dw_i\right\} dx \\ 
= & \int p(\hat{\bm{Z}}_{A\cup\{x\}}) dx.
\end{align*}

For $x\in \Omega\backslash\D$, since $x$ is independent of all other locations when conditional on $\D$, we have
\begin{align*}
p(\hat{\bm{Z}}_{A}) 
= & \int p(\hat{\bm{Z}}_U|\hat{\bm{Z}}_{\D}) \prod_{i=1}^n p(\hat{Z}_{w_i}|\hat{Z}_{\mr{pa}(w_i)}) \prod_{w_i\in \D\backslash V} dw_i \\ 
= & \int\left\{\int p(\hat{\bm{Z}}_U|\hat{\bm{Z}}_{\D}) \prod_{i=1}^n p(\hat{Z}_{w_i}|\hat{Z}_{\mr{pa}(w_i)}) \prod_{w_i\in \D\backslash V} dw_i\right\} p(\hat{Z}_x|\hat{Z}_\D) dx \\ 
= & \int p(\hat{\bm{Z}}_{A\cup\{x\}}) dx.
\end{align*}

\underline{\textsc{Conclusion}}
Using the Kolmogorov extension theorem, we conclude that $\hat{\bm{Z}}$ is a well-defined stochastic process. Since we have also verified each finite-dimensional distribution of $\hat{\bm{Z}}$ is Gaussian, we further conclude that $\hat{\bm{Z}}$ is a valid Gaussian process.
\end{proof}

\section{Auxiliary Lemmas for Main Theorems}\label{sec:aux.lemma}

\subsection{Upper Bound for $W_2$ Distance}
We begin our theoretical analysis by studying the upper bound for the quantity of interest, the Wasserstein-2 distance.
Let a decomposition of $\Sigma_{\D\D}^{-1}$ be $\Sigma_{\D\D}^{-1}=LL^T$ with $L=(l_1,\ldots,l_n)$; similarly let $\hat{\Phi}=\hat L\hat L^T$ with $\hat L=(\hat l_1,\ldots,\hat l_n)$. 
For a generic $m\times m$ symmetric matrix $A$, let $\lambda_{\max}(A)$ be its largest eigenvalue. For a generic $k\times m$ matrix $B$, let $\|B\|_2 = \{\lambda_{\max}(B^{\Tr} B)\}^{1/2}$ and $\|B\|_{\tr} =\tr\{(B^{\Tr} B)^{1/2}\}$.
\begin{lemma}\label{W2 decom}
For any decomposition $\Sigma_{\D\D}^{-1}= L  L^\Tr$ and $\hat{ {\Phi}} = \hat{ L} \hat{ L}^\Tr$, if $\|\hat{L} - L\|_2 \le \|\Sigma_{\D\D}\|_2^{-1/2} / 2$, then we have
\begin{align*}
W_2^2(Z_{\D}, \hat{Z}_{\D}) 
& \le \|\Sigma_{\D\D} - \hat{\Sigma}_{\D\D}\|_{\tr}   \\
& \le 8n\|\Sigma_{\D\D}\|_2^2  \big(2\max_i \|{l}_i\|_2 \max_i \|{l}_i-\hat{{l}}_i\|_2 + \max_i \|{l}_i-\hat{{l}}_i\|_2^2\big).
\end{align*}
\end{lemma}

\begin{proof}
By Proposition 1 of \cite{quang2021convergence}, we have
\begin{equation}\label{eq W2 tr}
W_2^2(\bm{Z}_{\D}, \hat{\bm{Z}}_{\D}) \le 
    \|\bm{\Sigma}_{\D\D} - \hat{\bm{\Sigma}}_{\D\D}\|_{\tr}.    
\end{equation}
Thus, it suffices to derive the upper bound for $\|\bm{\Sigma}_{\D\D} - \hat{\bm{\Sigma}}_{\D\D}\|_{\tr}$.
Plugging in the decomposition of $\bm{\Sigma}_{\D\D}^{-1}$ and $\hat{\bm{\Phi}}$, we have that
\begin{align}\label{eq trace decom}
\bm{\Sigma}_{\D\D}-\hat{\bm{\Sigma}}_{\D\D} 
= & \bm{\Sigma}_{\D\D} (\hat{\bm\Phi}-\bm{\Sigma}_{\D\D}^{-1}) \hat{\bm\Phi}^{-1} \nonumber\\
= & \bm{\Sigma}_{\D\D} (\hat{\bm L}\hat{\bm L}^{\Tr}-\bm{L}\bm{L}^{\Tr}) \hat{\bm\Phi}^{-1} \nonumber\\
= &  \sum_{i=1}^n\left\{\bm{\Sigma}_{\D\D} \hat{\bm l}_i (\hat{\bm\Phi}^{-1}\hat{\bm l}_i)^{\Tr} - \bm{\Sigma}_{\D\D} \bm{l}_i(\hat{\bm\Phi}^{-1}\bm{l}_i)^{\Tr} \right\}.
\end{align}

We need an auxiliary result regarding the trace norm of matrix operations.
\begin{techlemma}\label{lem trace op}
For any $\bm{u}_1,\bm{u}_2,\bm{v}_1,\bm{v_2}\in\mathbb{R}^n$, we have 
$$\|\bm{u}_2\bm{v}_2^{\Tr}-\bm{u}_1\bm{v}_1^{\Tr}\|_{\tr}\le 2(\|\bm u_2- \bm u_1\|_2\|\bm v_1\|_2 + \|\bm v_2- \bm v_1\|_2\|\bm u_1\|_2 + \|\bm u_2- \bm u_1\|_2\|\bm v_2-\bm v_1\|_2).$$
\end{techlemma}
\begin{proof}
We first decompose the $l_2$ norm of $\bm{u}_2\bm{v}_2^{\Tr}-\bm{u}_1\bm{v}_1^{\Tr}$ as
\begin{align*}
\|\bm{u}_2\bm{v}_2^{\Tr}-\bm{u}_1\bm{v}_1^{\Tr}\|
& \le \|\bm u_2- \bm u_1\|_2\|\bm v_1\|_2 + \|\bm v_2- \bm v_1\|_2\|\bm u_1\|_2 + \|\bm u_2- \bm u_1\|_2\|\bm v_2-\bm v_1\|_2.
\end{align*}
Since the right eigenvectors corresponding to nonzero singular values of matrix $\bm{u}_2\bm{v}_2^{\Tr}-\bm{u}_1\bm{v}_1^{\Tr}$ must lie in a rank two space $\mr{span}\{\bm{v}_1,\bm v_2\}$, there are at most two nonzero singular values for this matrix. Therefore
\begin{align*}
& \|\bm{u}_2\bm{v}_2^{\Tr}-\bm{u}_1\bm{v}_1^{\Tr}\|_{\tr}\le 2\|\bm{u}_2\bm{v}_2^{\Tr}-\bm{u}_1\bm{v}_1^{\Tr}\|_2 
\end{align*}
Combining the above two equations completes the proof. \end{proof}

We now use the result of Technical Lemma \ref{lem trace op} and plug Equation (\ref{eq trace decom}) into Equation (\ref{eq W2 tr}) to derive that
\begin{align*}
&\quad~ \|\bm{\Sigma}_{\D\D}-\hat{\bm{\Sigma}}_{\D\D}\|_{\tr} 
\le  \sum_{i=1}^n \|\bm{\Sigma}_{\D\D}\hat{\bm{l}}_i (\hat{\bm{\Phi}}^{-1}\hat{\bm{l}}_i)^{\Tr}- \bm{\Sigma}_{\D\D}\bm{l}_i (\hat{\bm{\Phi}}^{-1}\bm{l}_i)^{\Tr}\|_{\tr} \\
&\le \sum_{i=1}^n 2\Big\{ \|\bm{\Sigma}_{\D\D}(\hat{\bm{l}}_i-\bm{l}_i)\|_2 \|\hat{\bm{\Phi}}^{-1}\bm{l}_i\|_2 + \|\bm{\Sigma}_{\D\D} \bm{l}_i\|_2 \|\hat{\bm{\Phi}}^{-1}(\hat{\bm{l}}_i-\bm{l}_i)\|_2 \\
&\qquad + \|\bm{\Sigma}_{\D\D}(\hat{\bm{l}}_i-\bm{l}_i)\|_2 \|\hat{\bm{\Phi}}^{-1}(\hat{\bm{l}}_i-\bm{l}_i)\|_2\Big\} \\
&\le \sum_{i=1}^n 2\|\bm{\Sigma}_{\D\D}\|_2 \|\hat{\bm{\Phi}}^{-1}\|_2 \big(2\|\hat{\bm l}_i-\bm l_i\|_2\|\bm l_i\|_2 + \|\hat{\bm l}_i-\bm l_i\|_2^2\big) \\
&\le 2n\|\bm{\Sigma}_{\D\D}\|_2 \|\hat{\bm{\Phi}}^{-1}\|_2 \big(2\max_i\|\hat{\bm l}_i-\bm l_i\|_2\max_i\|\bm l_i\|_2 + \max_i\|\hat{\bm l}_i-\bm l_i\|_2^2\big).
\end{align*}
If $\|\hat{L} - L\|_2 \le \|\Sigma_{\D\D}\|_2^{-1/2} / 2$, then we have that $\lambda_{\min} (\hat L) \ge \lambda_{\min} (L)/2$ and hence $\lambda_{\min}(\hat\Phi) \ge \lambda_{\min}\{(\Sigma_{\D\D})^{-1}\}/4$. This gives the upper bound $4\|\hat\Phi^{-1}\|_2 \le \|\Sigma_{\D\D}\|_2$. Thus
$$\|\bm{\Sigma}_{\D\D}-\hat{\bm{\Sigma}}_{\D\D}\|_{\tr} \le 8n\|\bm{\Sigma}_{\D\D}\|_2^2 \big(2\max_i\|\hat{\bm l}_i-\bm l_i\|_2\max_i\|\bm l_i\|_2 + \max_i\|\hat{\bm l}_i-\bm l_i\|_2^2\big). $$
\mbox{}
\end{proof}

\subsection{Matrix Norms}
Following the upper bound of Wasserstein distance in Lemma \ref{W2 decom}, we provide estimates for various matrix norms that are related to this upper bound.

\begin{lemma}\label{D bounds}
The following bounds regarding the set $\D$ hold:
\begin{enumerate}[(1)]
\item If $\bm{Z}\in\mathscr{Z}_{v_r}$ for some $r>0$, then 
$$\|\bm{\Sigma}_{\DD}\|_2 \le \|\bm{\Sigma}_{\DD}\|_1 \le 1+d 2^d q^{-d} 
\int_{q/2}^{1} x^{d-1} \frac{1}{v_r(x-q/2)\{1+(x-q/2)^2\}}dx.$$ 
\item If $\bm{Z}\in\mathscr{Z}_{c_r}$ for some $r\ge d+1$, then 
$$\|\bm{\Sigma}_{\DD}\|_2 \le \|\bm{\Sigma}_{\DD}\|_1 \le 1 + d2^{d}q^{-d}\max\{1,(q/2)^{d-1}\}.$$
\item Let $c_2 = 12\left\{\pi \Gamma^2(d/2+1)/9 \right\}^{1/(d+1)}$, $c_1 = 2\Gamma(d/2+1)(2^{3/2}/c_2)^d$ and $\phi_0(c_2/q) = \inf_{\|\omega\|_2\le 2c_2/q} \hat{K}_0(\omega)$ with $\hat{K}_0$ being the Fourier transform of $d$-dimensional function $K_0$, then
$$\|\bm{\Sigma}_{\DD}^{-1}\|_2 \le c_1 \{\phi_0(c_2/q')\}^{-1} q'^d, \quad \text{ for all } q' \le q.$$
\item If $\bm{Z}\in\mathscr{Z}_{v_r}$ for some $r>0$, then $\|\bm{\Sigma}_{\DD}\|_{v_r,l} \le 2^dq^{-d}(1+q/2)^{d-1} (1+ q/2 + d)$ for all $l\in \mathbb{N}$.
\item If $\bm{Z}\in\mathscr{Z}_{c_r}$ for some $r\ge d+1$, then
$\|\bm{\Sigma}_{\DD}\|_{c_{r-(d+1)/2}} \le d 2^{d-1} \pi/ \sqrt{6}$.
\end{enumerate}
\end{lemma}
\begin{proof}
(1)
Since the matrix $l_2$ norm is bounded by the matrix $l_1$ norm, we have
\begin{equation}\label{eq Sigma l2 1}
\|\bm{ \Sigma}_{\D\D}\|_2 
\le \|\bm{ \Sigma}_{\D\D}\|_1 = \sup_i \sum_{j=1}^n K(w_i,w_j).
\end{equation}
Now consider the term $\sum_{j=1}^n K(w_i,w_j)$. Recall that $q=\min_{1\le i<j \le |\D|}\|w_i-w_j\|$ is the minimal separation distance among all locations in $\D$. Define an auxiliary function $\varphi: \Omega \to \Omega$ such that
(1) $\|\varphi(s)-s\|_2< q/2$;
(2) If there exists some $w_i\in \D$ such that $\|w_i-s\|<q/2$, then $\varphi(s) = w_i$.

The function $\varphi$ maps the unit ball $U(w_i, q/2), 1\le i\le n$ into a singleton $\{w_i\}$. We extend the definition of $v_r(\cdot)$ to $(-\infty,0)$ by letting $v_r(x)=1$ for all $x<0$. Then we have that
\begin{align}\label{eq Sigma l2 2}
\sum_{j=1}^n K(w_i,w_j) 
&= \frac{1}{\pi^{d/2} (q/2)^d /\Gamma(d/2+1)}\int_{s\in\bigcup_{j=1}^n U(w_j, q/2)} K(w_i,w_j) ds \nonumber\\
&= \frac{1}{\pi^{d/2} (q/2)^d /\Gamma(d/2+1)}\int_{s\in\bigcup_{j=1}^n U(w_j, q/2)} K(w_i,\varphi(s)) ds \nonumber\\
&\le \frac{\Gamma(d/2+1) 2^d}{\pi^{d/2} q^d} \int_{s\in\Omega} K(w_i,\varphi(s)) ds \nonumber\\
&\le \frac{\Gamma(d/2+1) 2^d}{\pi^{d/2} q^d} \int_{s\in\Omega} \frac{1}{v_r(\max\{0,\|s\|_2-q/2\}) [1+\{\max(0,\|s\|_2-q/2)\}^{d+1}]} ds \nonumber\\
&\overset{\Omega\subset\mathbb{R}^d}{\le} \frac{\Gamma(d/2+1) 2^d}{\pi^{d/2} q^d} \int_0^{+\infty} \frac{2\pi^{d/2} x^{d-1}}{\Gamma(d/2)} \nonumber \\
&\quad \times \frac{1}{v_r(\max\{0,x-q/2\}) [1+\{\max(0,x-q/2)\}^{d+1}]} dx \nonumber\\
&= 1 + d2^dq^{-d} \int_{q/2}^{+\infty} \frac{x^{d-1}}{v_r(x-q/2) \{1+(x-q/2)^{d+1}\}} dx,
\end{align}
Combining equations (\ref{eq Sigma l2 1}) and (\ref{eq Sigma l2 2}), we finish the proof.
\vspace{2mm}

{ (2)}
Using the same notation as the proof of (1) above, we have
\begin{align*}
\sum_{j=1}^n K(w_i,w_j) 
&\le \frac{\Gamma(d/2+1) 2^d}{\pi^{d/2} q^d} \int_{s\in\Omega} K(w_i,\varphi(s)) ds \nonumber\\
&\le \frac{\Gamma(d/2+1) 2^d}{\pi^{d/2} q^d} \int_{s\in\Omega} \frac{1}{ (1+\max\{0,\|s\|_2-q/2\})^{r}} ds \nonumber\\
&\overset{\Omega\subset\mathbb{R}^d}{\le} \frac{\Gamma(d/2+1) 2^d}{\pi^{d/2} q^d} \int_0^{+\infty} \frac{2\pi^{d/2} x^{d-1}}{\Gamma(d/2)} \frac{1}{(1+\max\{0,\|x\|_2-q/2\})^{r}} dx \nonumber\\
&= 1 + d2^dq^{-d} \int_{q/2}^{+\infty} \frac{x^{d-1}}{(1+x-q/2)^{r}} dx \\
&\overset{r\ge d+1}{\le}  1 +  d2^{d}q^{-d}\max\{1,(q/2)^{d-1}\}.
\end{align*}

{ (3)} The conclusion follows directly from Theorem 12.3 of \cite{wendland2004scattered}.

{ (4)} For a vector $w\in \mathbb{R}^d$, we use $w[j]$ to denote its $j$th component. We have that for any $l \in \mathbb{N}$, 
\begin{align*}
\|\bm{\Sigma}_{\D\D}\|_{v_r,l}
= & \sum_{k=0}^{+\infty} \frac{1}{(k!)^r} \|\nabla_l^k(\bm{\Sigma})\|_2 \\
\le & \sum_{k=0}^{+\infty} \frac{1}{(k!)^r} \|\nabla_l^k(\bm{\Sigma})\|_1 \\
= & \sum_{k=0}^{+\infty} \frac{1}{(k!)^r} \sup_i\sum_{j=1}^n |w_i[l]-w_j[l]|^k K(w_i,w_j) \\
\le & \sum_{k=0}^{+\infty} \frac{1}{(k!)^r} \sup_i\sum_{j=1}^n |w_i[l]-w_j[l]|^k  \frac{1}{v_r(\|w_i-w_j\|_2)(\|w_i-w_j\|_2^2+1)} \\
= & \sup_i\sum_{j=1}^n \frac{v_r(|w_i[l]-w_j[l]|)}{v_r(\|w_i-w_j\|_2)(\|w_i-w_j\|_2^2+1)} \\
\le & \sup_i \sum_{j=1}^n \frac{1}{\|w_i-w_j\|^{d+1}+1}.
\end{align*}
Using the same trick as in the proof of (1) to turn the summation into an integration, we have that
\begin{align*}
& \sup_i \sum_{j=1}^n \frac{1}{\|w_i-w_j\|^{d+1}+1} \\
= {} & \frac{1}{\pi^{d/2} (q/2)^d /\Gamma(d/2+1)}\int_{s\in\bigcup_{j=1}^n U(w_j, q/2)} \frac{1}{\|w_i-\varphi(s)\|^{d+1}+1} ds \nonumber\\
\le {} & \frac{\Gamma (d/2+1)2^d}{\pi^{d/2} q^d} \int_{s\in\Omega} \frac{1}{\|w_i-\varphi(s)\|^{d+1}+1}ds \\
\le {} & d2^dq^{-d}\int_0^{+\infty} \frac{x^{d-1}}{(\max\{0,x-q/2\})^2+1} dx \\
\le {} & d2^dq^{-d} \left\{\int_0^{1+q/2} x^{d-1} dx + \int_{1+q/2}^{+\infty} \frac{x^{d-1}}{(x-q/2)^{d+1}} dx\right\}\\
\le {} & 2^dq^{-d}(1+q/2)^d + d2^dq^{-d} (1+q/2)^{d-1}\int_{1+q/2}^{+\infty}\frac{1}{(x-q/2)^2} dx \\
= {} & 2^dq^{-d}(1+q/2)^d + d2^dq^{-d} (1+q/2)^{d-1} \\
= {} & 2^dq^{-d}(1+q/2)^{d-1} (1+ q/2+d).
\end{align*}

{(5)}
\begin{align*}
\|\bm{\Sigma}_{\D\D}\|_{c_{r-(d+1)/2}} 
&\le \Bigg\{\sum_{k\in\mathbb{Z}^d}\sup_{w_i,w_j\in\D} (1+\|w_i-w_j\|_\infty)^{-2r}  \nonumber \\
&\qquad\times (1+\|w_i-w_j\|_\infty)^{2(r-(d+1)/2)} \mathbbm{1}_{\{w_i-w_j\in [0,1)^d+k\}}  \Bigg\}^{\frac{1}{2}} \\
&\overset{\Omega\subset\mathbb{R}^d}{\le} \left\{\sum_{k\in\mathbb{Z}^d}\sup_{s\in [0,1)^d+k} (1+\|s\|_\infty)^{-(d+1)} \right\}^{\frac{1}{2}} \\
&\le \left\{\sum_{k=1}^{\infty} [(2k)^d-\{2(k-1)\}^d] k^{-(d+1)} \right\}^{\frac{1}{2}} \\
&\le \left\{\sum_{k=1}^{\infty} d(2k)^{d-1} k^{-(d+1)} \right\}^{\frac{1}{2}} \\
&= d2^{d-1} \Big(\sum_{k=1}^{\infty} k^{-2}  \Big)^{\frac{1}{2}}
= d2^{d-1} \frac{\pi}{\sqrt{6}}.
\end{align*}
\mbox{}
\end{proof}

\subsection{Bounds regarding Matrix Inversion}
To utilize the upper bound proved in Lemma \ref{W2 decom}, it is necessary to bound the quantity $\max_i\|l_i-\hat{l}_i\|_2$, which involves the Cholesky decomposition of an inverse matrix. In this section, we present some general tools that describe the properties of matrix inversion in some norms.

We first discuss the Cholesky decomposition that yields matrix $L$ and $\hat{L}$.
We introduce the notation of quoting a submatrix by its indices that will only be used in this subsection.
For a generic matrix $B$, let $B[i,j]$ be its entry at the $i$th row and $j$th column. Recall $\nu(i)=\{j\in\mathbb{N}: j<i\}$, further define $N(i) = \{j\in\nu(i): \|w_i-w_j\|_2<\rho\}$. Our directed acyclic graphs constructed in Section \ref{subsec:construction} of the main text induce the following decomposition of the precision matrices $\Sigma_{\D\D}^{-1}$ from the original Gaussian processes and $\hat{\Phi}$ from the radial neighbors Gaussian processes.

\begin{lemma}\label{lem chol pre}
The precision matrices $\Sigma_{\D\D}^{-1}$ and $\hat\Phi$ satisfy
$$\Sigma_{\D\D}^{-1} = (I_n-B^\Tr)D^{-1}(I_n-B), \quad \hat{\Phi} = (I_n-\hat B^\Tr)\hat D^{-1}(I_n-\hat B),$$
where $B$ and $\hat B$ are lower triangular matrices such that 
$$B[\nu(i),i]= (\Sigma[\nu(i),\nu(i)])^{-1}\Sigma[\nu(i),i]\;\; \mr{and} \;\;B[j,i]=0,\;\forall\; j\not\in\nu(i),$$
$$\hat B[N(i),i]= (\Sigma[N(i),N(i)])^{-1}\Sigma[N(i),i] \;\;\mr{and}\;\; \hat B[j,i]=0,\;\forall \; j\not\in N(i).$$
$D$ and $\hat D$ are diagonal matrices, such that $$D[i,i]=K(w_i,w_i)-\Sigma[\nu(i),i]^\Tr (\Sigma[\nu(i),\nu(i)])^{-1}\Sigma[\nu(i),i],$$ 
$$\hat D[i,i]=K(w_i,w_i) - \Sigma[N(i),i]^\Tr (\Sigma[N(i),N(i)])^{-1} \Sigma[N(i),i].$$
\end{lemma}
\vspace{2mm}

\begin{proof}
Radial neighbors Gaussian process implies an algorithm to find a sparse approximation for the Cholesky of the inverse of a matrix. 
We begin by constructing the exact decomposition $\bm{\Sigma}_{\DD}^{-1}=\bm L \bm L^{\Tr}$ with $\bm L = (\bm l_1, \cdots \bm l_n)$. By Bayes' rule, the joint density of $\bm{Z}_{\D}$ can be decomposed as
\begin{equation}\label{eq full DAG}
p(\bm{Z}_{\D}) = p(\bm{Z}_{w_1})\prod_{i=2}^N p(\bm{Z}_{w_i}|\bm{Z}_{w_1}, \ldots, \bm{Z}_{w_{i-1}}).
\end{equation}
The decomposition of density induces a decomposition on the precision matrix of Gaussian distribution $p(\bm{Z}_{\D})$. Specifically, writing (\ref{eq full DAG}) in the form of conditional regression, we have
\begin{align*}
\bm{Z}_{w_1} &= \eta_1, \qquad
\bm{Z}_{w_i} =  \sum_{j=1}^{i-1} b_{i, j} \bm{Z}_{w_j} + \eta_i, \text{ for all } 2\le i\le n ,
\end{align*}
where $b_{i, j},1\le j \le i-1$ are conditional regression coefficients satisfying
$$\begin{bmatrix} b_{i, 1} & b_{i, 2} & \cdots & b_{i, i-1}\end{bmatrix}^{\Tr}
= (\bm{\Sigma}[\nu(i),\nu(i)])^{-1} \bm{\Sigma}[\nu(i),i]
$$
and
$\eta_i, 1\le i\le n$ are independent mean zero Gaussian random variables with variance 
$$\var(\eta_i) = \var(\bm{Z}_{w_i}|\bm{Z}_{w_1},\ldots, \bm{Z}_{w_{i-1}}) = K(w_{i},w_{i}) - (\bm{\Sigma}[\nu(i),i])^{\Tr} (\bm{\Sigma}[\nu(i),\nu(i)])^{-1} \bm{\Sigma}[\nu(i),i].
$$
Define the coefficient matrix $\bm B$ such that $\bm B[i,j] = b_{i, j}$ if and only if $i>j$ and $\bm B[i,j] = 0$ otherwise. That is, $B[i,\nu(i)]^T= (\bm{\Sigma}[\nu(i),\nu(i)])^{-1} \bm{\Sigma}[\nu(i),i]$.
Also define a diagonal matrix $\bm D$ such that $\bm D[i,i] = \var(\eta_{i})$. Then by the equality $\bm w = \bm B \bm w + \bm \eta$, we have that
$$
\bm{\Sigma}^{-1}_{\D\D} = (\bm I_n - \bm{B}^{\Tr}) \bm{D}^{-1} (\bm I_n - \bm B).$$
To obtain the decomposition $\bm{\Sigma}_{\D\D}^{-1}=\bm L \bm L^{\Tr}$, we let $\bm L$ be
\begin{equation*}
\bm L = (\bm I_n - \bm{B}^{\Tr}) \bm{D}^{-1/2}.
\end{equation*}
This has proved the first decomposition in Lemma \ref{lem chol pre}.

For the decomposition of the precision matrix $\hat{\bm{\Phi}}$ of $\hat{\bm{Z}}$, we now use a similar way to derive the decomposition of $\hat{\bm{\Phi}}$. Noticing the definition of ordering $(w_i)_{1\le i\le n}$ implies $\mr{pa}(w_i)\subset \nu(i)$, we have
$$ 
p(\hat{\bm{Z}}_{\D}) = p(\hat{\bm{Z}}_{w_1}) \prod_{i=2}^n  p(\hat{\bm{Z}}_{w_i}|\hat{\bm{Z}}_{\mr{pa}(w_i)}). $$
Denote $N(i)=\{j\in\nu(i):\|w_i-w_j\|_2<\rho\}$. The above equation similarly induces a decomposition of $\hat{\bm{\Phi}}_{\D\D}$ as
\begin{equation*}
\hat{\bm{\Phi}}_{\D\D} = (\bm I_n - \hat{\bm B}^{\Tr})\hat{\bm{D}}^{-1}(\bm I_n - \hat{\bm B}),    
\end{equation*}
where $\hat{\bm B}$ is an $n\times n$ matrix such that $\hat{\bm B}[i,j] = 0$ for all $j\not\in N(i)$. For the nonzero elements of $\hat{\bm B}$, we have for all $1\le i\le n$
\begin{equation*}
\hat{\bm{B}}[i,N(i)]^T = \bm{\Sigma}_{N(i),N(i)}^{-1} \bm{\Sigma}_{N(i),i}.
\end{equation*}
The matrix $\hat{\bm{D}}$ is a diagonal matrix with entries
\begin{equation*}
\hat{\bm D}[i,i] = K(w_i,w_i) - \bm{\Sigma}_{N(i),i}^{\Tr} \bm{\Sigma}_{N(i),i}^{-1} \bm{\Sigma}_{N(i),i}.    
\end{equation*}
The matrix decomposition $\hat{\bm{\Phi}}=\hat{\bm L}\hat{\bm L}^{\Tr}$ is defined by 
\begin{equation*}
\hat{\bm L} = (\bm I_n - \hat{\bm{B}}^{\Tr}) \hat{\bm{D}}^{-1/2}.
\end{equation*}
\mbox{}
\end{proof}

We next introduce a result from literature of norm-controlled inversion of Banach Algebra that helps us transfer the spatial decay properties from covariance matrices to precision matrices and the Cholesky of precision matrices.

Specifically, let the matrix $A\in\mathbb{R}^{n\times n}$ be associated with $n$ spatial locations $w_1,\ldots, w_n$ such that the $(i,j)$-entry of $A$, denoted by $A[i,j]$, is a function of the difference $w_i-w_j$. Let $w_i[l]$ be the $l$th coordinate of $w_i$. For all $1\le l\le d$, define a linear matrix operator $\nabla_l$ such that 
$$\nabla_l({A})[i,j]= \big(w_i[l]-w_j[l]\big) {A}[i,j].$$
For such a square matrix $A$, we define its $l$th order $v_r$ norm and $c_r$ norm as 
$$\|{A}\|_{v_r,l} = \sum_{k=0}^{+\infty} \|\nabla_l^k({A})\|_2/\{(k!)^r\},$$
$$\|{A}\|_{c_r} = \Big\{\sum_{k\in\mathbb{Z}^d}\sup_{w_i,w_j\in\D} ({A}[i,j])^2(1+\|w_i-w_j\|_\infty)^{2r} \mathbbm{1}_{\{w_i-w_j\in [0,1)^d+k\}}  \Big\}^{1/2},$$
which are related to the functions $v_r$ and $c_r$ defined above.
These two matrix norms describe the spatial decaying properties in the sense that if a matrix $A$ has a finite $\|\cdot\|_{v_r,l}$ norm (or $\|\cdot\|_{c_r}$ norm), then its $(i,j)$-entry decays at the rate $1/v_r(\|w_i-w_j\|_\infty)$ (or $1/c_r(\|w_i-w_j\|_\infty)$).
The following lemma is based on the work of \cite{grochenig2014norm} and \cite{fang2020norm}.

\begin{lemma}\label{norm-controlled}
Let ${A}\in n\times n$ be an invertible matrix. 
If $\|A\|_{v_r,l}<\infty$ for some $r>1$ and all $1\le l\le d$, then
$$\|{A}^{-1}\|_{v_r,l} \le \|{A}^{-1}\|_2 v_{r-1}(\|{A}^{-1}\|_2  \|{A}\|_{v_r,l}) .$$
If $\|\cdot\|_{c_r}<\infty$ for some $r\ge d+1$, then there exist positive constants $c_4, c_5$ only dependent on the dimension $d$, such that
$$\|{A}^{-1}\|_{c_r} \le c_4 c_5^r \{\max(1, 1/q)\}^{6d} \|{A}^{-1}\|_2 \big(\|{A}^{-1}\|_2  \|{A}\|_{c_r}\big)^{r+d/2}.$$
\end{lemma}

\begin{proof}
\underline{\textsc{Part 1}} We prove the inequality for $\|\cdot\|_{v_r,l}$ norm.
The idea is to first prove the bound holds if $A$ is an infinite matrix using the theory of norm-controlled inversion, and then embed a finite-dimensional matrix into an infinite matrix.

Let $A\in\mathbb{R}^{\mathbb{Z}\times\mathbb{Z}}$ be an infinite matrix with finite $\|\cdot \|_{v_r,l}$ norm. We have two basic facts:
\vspace{2mm}

\noindent (1) the collection of matrices of $\mathbb{R}^{\mathbb{Z}\times\mathbb{Z}}$ with finite matrix $L_2$ norm forms an algebra with matrix addition and multiplication as algebra addition and multiplication; 
\vspace{2mm}

\noindent (2) $\nabla_{l}$ is a differential operator on the algebra of matrices with finite $L_2$ norm. 
\vspace{2mm}

Thus we can apply the norm-controlled inversion theory in Section 2.4 of \cite{grochenig2014norm} to bound the inversion $\bm{A}^{-1}$. Specifically, let $\bm{B} = \bm{A}/\|\bm{A}\|_{v_{r,l}}$ for some $r \ge 1$ and $l\in \mathbb{N}$, then applying equation (2.26) of \cite{grochenig2014norm} yields
\begin{equation*}
\|\bm{B}^{-1}\|_{v_r,l} \le \|\bm{B}^{-1}\|_2 v_{r-1}(\|\bm{B}^{-1}\|_2).
\end{equation*}
Replacing $\bm{B}^{-1}=\|\bm{A}\|_{v_r,l}\bm{A}^{-1}$, we have
\begin{equation*}
\|\bm{A}^{-1}\|_{v_r,l} \le \|\bm{A}^{-1}\|_2 v_{r-1}(\|\bm{A}^{-1}\|_2\|\bm{A}\|_{v_r,l}).
\end{equation*}

Now consider the case when $A$ is a finite matrix. Without loss of generality, we assume $\mathcal{I}=\{1,2,\ldots, n\}$ for $n\in\mathbb{N}$. Now expand $\mathcal{I}$ into $\mathbb{Z}$ by assigning each $i\in\mathbb{Z}\backslash \mathcal{I}$ a location $w_i\in\Omega$ such that for any $ i,j\in\mathbb{Z}$ with $i\ne j$, $w_i\ne w_j$. We embed $\bm{A}$ into a matrix $\tilde{\bm{A}}\in\mathbb{R}^{\mathbb{Z}\times\mathbb{Z}}$ such that $$\tilde{\bm{A}}[i,j] = \bigg\{ 
\begin{aligned}
& \bm{A}[i,j], & i\in\mathcal{I} \text{ and } j\in\mathcal{I}, \\
& \mathbbm{1}_{\{i=j\}}, & i\not\in\mathcal{I} \text{ or } j\not\in\mathcal{I}.
\end{aligned}$$
Then again applying equation (2.26) of \cite{grochenig2014norm} on $\tilde{\bm{A}}/\|\tilde{\bm{A}}\|_{v_r,l}$, we similarly obtain
\begin{equation}\label{eq tilde A ori}
\|\tilde{\bm{A}}^{-1}\|_{v_r,l} \le \|\tilde{\bm{A}}^{-1}\|_2 v_{r-1}(\|\tilde{\bm{A}}^{-1}\|_2\|\tilde{\bm{A}}\|_{v_r,l}).
\end{equation}
Since the matrix $\tilde{\bm{A}}$ can be represented as a block diagonal matrix
\begin{equation}\label{eq tilde A bound}
\tilde{\bm{A}} = \begin{bmatrix} \bm{A} & \bm{0} \\ \bm{0} & \bm{I}_{\mathbb{Z}\backslash\mathcal{I}} \end{bmatrix},   
\end{equation}
we have
$$\nabla_l^k(\tilde{\bm{A}}) = \begin{bmatrix} \nabla_l^k(\tilde{\bm{A}}) & \bm{0} \\ \bm{0} & \bm{I}_{\mathbb{Z}\backslash\mathcal{I}} \end{bmatrix},$$
which implies that
\begin{equation}\label{eq tilde A norm}
\|\tilde{\bm{A}}^{-1}\|_2 = \|\bm{A}^{-1}\|_2,\;\; \|\tilde{\bm{A}}^{-1}\|_{v_r,l} =  \|\bm{A}^{-1}\|_{v_r,l}.
\end{equation}
Combining equation (\ref{eq tilde A norm}) with (\ref{eq tilde A bound}), we obtain that
$$
\|\bm{A}^{-1}\|_{v_r,l} \le \|\bm{A}^{-1}\|_2 v_{r-1}(\|\bm{A}^{-1}\|_2\|\bm{A}\|_{v_r,l}).
$$

\underline{\textsc{Part 2}} We prove the inequality for $\|\cdot\|_{c_r,l}$ norm.
Based on the same trick of part 1 that extends a finite matrix to an infinite matrix, Theorem 2 of \cite{fang2020norm} directly yields
\begin{align} \label{eq:fang2020norm}
\|\bm{A}^{-1}\|_{c_r} \le C \|\bm{A}^{-1}\|_2 \big(\|\bm{A}^{-1}\|_2  \|\bm{A}\|_{c_r}\big)^{r+\frac{d}{2}} ,
\end{align}
for some positive constant $C$ dependent on $r$, $d$ and $R(\Lambda)$ in their paper. We only need to verify we can find $C', C''$ dependent on $d$, such that $C\le C' {C''}^r (\max\{1, 1/q\})^6$. In the rest of the proof, we use $O(\cdot)$ to denote a quantity is no greater than some constants only dependent on $d$ and $r$.

We first notice their $R(\Lambda)$ is the maximal number of locations in a cube $k+[0,1)^d$. Thus we have
$$R(\Lambda) \le \{\max(1, \lceil 1/q\rceil)\}^d = O(\{\max(1,1/q)\}^d).$$
By the proof of their Proposition 1 (4), and our condition $r\geq d+1$ in Lemma \ref{norm-controlled}, their constant $C_1$ satisfies
$$C_1 = 2^{r+1} R(\Lambda) \left(\frac{3^d r}{r-d/2}\right)^{1/2}
= O\left(2^r \{\max(1,1/q)\}^d \right).$$
By the definition of $D_{d,p,r}$ in their proof of Theorem 2, we have
\begin{align*}
D_{d,p,r} = 2^{3r+3} 11^{d+r} &\{\max(1,1/q)\}^d \left\{\left(\frac{d}{2r-d}\right)^{1/2} + \left(\frac{d+2r}{2r-d-2}\right)^{1/2} \right\} \\&= O \big(88^r \{\max(1,1/q)\}^d \big).\end{align*}
Finally, by their equations (42), (34) and (35), the constant $C$ in \eqref{eq:fang2020norm} satisfies that
\begin{align*}
C \le{} & 10^{d/2+r+d+1} N_0^{d/p+r} \\
={} & O\Big(10^r (C_1 D_{d,p,r})^{\sup_{r\ge d+1} \frac{r+d/2}{\min\{1,r-d/2\}} }   \Big) \\
={} & O\Big( 10^r (C_1 D_{d,p,r})^3   \Big) \\
={} & O\Big( 10^r \big[(2\times88)^r \{\max(1,1/q)\}^{2d}\big]^3 \Big) \\
={} & O\Big( (10\times 176^3)^r \{\max(1,1/q)\}^{6d} \Big).
\end{align*}
The conclusion of Lemma \ref{norm-controlled} follows from this upper bound on $C$ and \eqref{eq:fang2020norm}.  
The constants computed above using the results in \cite{fang2020norm} are merely for the purpose of proving the existence of the constant $C$. In the literature of norm-controlled inversion, these constants are not carefully tuned and not tight in general.
\end{proof}

\subsection{Bounds for $\max_i\|l_i-\hat{l}_i\|_2$ under Decaying Covariance Functions}

We are finally ready to present upper bounds for $\max_i\|l_i-l\|_2$. The first result applies to the situation where the decay rate of the covariance function is no slower than $1/v_r$.
\begin{lemma}\label{l bounds}
Suppose that $\bm{Z}\in\mathscr{Z}_{v_r}$ for some $r>1$. 
If $0<q<1$ and
\begin{align*} 
\frac{ n^{1/2} }{v_r( \rho d^{-1/2} )}\{\phi_0(c_2/q)\}^{-9/2} v_{r-1}(c_3\{\phi_0(c_2/q)\}^{-1})  \le c_6
\end{align*}
for some constant $c_6$ only dependent on $d$, then
\begin{align}
\max_i \|\bm{l}_i - \hat{\bm l}_i\|_2 & \lesssim \frac{1}{v_r( \rho d^{-1/2} )}\{\phi_0(c_2/q)\}^{-\frac{9}{2}}q^{\frac{1}{2}d} v_{r-1}(c_3\{\phi_0(c_2/q)\}^{-1}), \label{eq:lemS3.l.diff1} \\
\mr{and}\;\; \|\hat{\bm{\Phi}} - \bm{\Sigma}^{-1}_{\D\D}\|_2 & \le  \frac{1}{2} \|\bm{\Sigma}_{\D\D}^{-1}\|_2. \label{eq phi prelim}
\end{align}
Else if $q\ge 1$ and
$ n^{1/2} \{v_r( \rho d^{-1/2})\}^{-1} \le c'_6$
for some constant $c'_6$ only dependent on $d$, then equation (\ref{eq phi prelim}) still holds, and
\begin{align}
\max_i\|\bm{l}_i-\hat{\bm{l}}_i\|_2 \lesssim\frac{1}{v_r( \rho d^{-1/2} )} . \label{eq:lemS3.l.diff2}
\end{align}
\end{lemma}

\begin{proof}\mbox{}
\underline{\textsc{Part 1}}
We first prove the upper bounds for $\max_i\|\bm{l}-\hat{\bm{l}}\|_2$ for both the cases of $0<q<1$ and $q\ge 1$.
Fix an arbitrary $i\in \{1,2,\ldots, n\}$. For simplicity of notation, define the sets $N$ and $O$ as
$$N=\{w_j\in\D: j<i, \|w_i-w_j\|_2<\rho\},$$ 
$$O= \{w_j\in\D:j<i,\|w_i-w_j\|_2\ge\rho\}.$$
With a little abuse of notation, we denote the $i$th column of matrix $\bm{I}_n-\bm{B}^{\Tr}$ and $\bm{I}_n-\hat{\bm{B}}^{\Tr}$ as $\bm{b}_i$ and $\hat{\bm{b}}_i$; denote $\bm{\Phi}_{NN}$ as the submatrix of $(\bm{\Sigma}[\nu(i),\nu(i)])^{-1}$ whose rows and columns correspond to set $N$; similarly define $\bm{\Phi}_{ON}$, $\bm{\Phi}_{NO}$ and $\bm{\Phi}_{OO}$. In the rest of the proof, we reorder the indices of elements in the sets $N$ and $O$ such that the indices of elements in $N$ are always smaller than those in $O$. In this way, we are able to formulate various computations as block matrix computations. By the definition of $\bm{\Phi}_{NN}$ and $\bm{\Phi}_{ON}$, we have
$$\bm{\Sigma}_{NN} \bm{\Phi}_{NN} + \bm{\Sigma}_{NO}\bm{\Phi}_{ON} = \bm{I}_N.$$
Thus
\begin{equation*}
\bm{\Phi}_{NN} = \bm{\Sigma}_{NN}^{-1} - \bm{\Sigma}_{NN}^{-1}\bm{\Sigma}_{NO}\bm{\Phi}_{ON}.
\end{equation*}
We can formulate $\bm{b}_i$ and $\hat{\bm{b}}_i$ as
$$\bm{b}_i = 
\begin{bmatrix}
1 \\ \bm{\Phi}_{NN}\bm{\Sigma}_{N,w_i}+\bm{\Phi}_{NO}\bm{\Sigma}_{O,w_i} \\
\bm{\Phi}_{ON}\bm{\Sigma}_{N,w_i}+\bm{\Phi}_{OO}\bm{\Sigma}_{O,w_i}
\end{bmatrix},
\quad
\hat{\bm{b}}_i = 
\begin{bmatrix}
1 \\ \bm{\Sigma}_{NN}^{-1}\bm{\Sigma}_{N,w_i}\\
0
\end{bmatrix}.
$$
Therefore we have
\begin{align}\label{eq b diff}
\|\bm{b}_i-\hat{\bm{b}}_i\|_2 
\le &  \|(\bm{\Phi}_{NN}-\bm{\Sigma}_{NN}^{-1})\bm{\Sigma}_{N,w_i}+\bm{\Phi}_{NO}\bm{\Sigma}_{O,w_i}\|_2 + \|\bm{\Phi}_{ON}\bm{\Sigma}_{N,w_i}+\bm{\Phi}_{OO}\bm{\Sigma}_{O,w_i}\|_2 \nonumber \\
\le & \|-\bm{\Sigma}_{NN}^{-1}\bm{\Sigma}_{NO}\bm{\Phi}_{ON}\bm{\Sigma}_{N,w_i}\|_2+\|\bm{\Phi}_{NO}\bm{\Sigma}_{O,w_i}\|_2 + \|\bm{\Phi}_{ON}\bm{\Sigma}_{N,w_i}\|_2+\|\bm{\Phi}_{OO}\bm{\Sigma}_{O,w_i}\|_2.
\end{align}
Similarly, for $\bm{D}[i,i]$ and $\hat{\bm{D}}[i,i]$, we have
$$\bm{D}[i,i] = \bm{\Sigma}_{w_i,w_i} - \bm{\Sigma}_{N,w_i}^{\Tr}\bm{\Phi}_{NN}\bm{\Sigma}_{N,w_i} - 2\bm{\Sigma}_{N,w_i}^{\Tr}\bm{\Phi}_{NO}\bm{\Sigma}_{O,w_i} - \bm{\Sigma}_{O,w_i}^{\Tr}\bm{\Phi}_{OO}\bm{\Sigma}_{O,w_i},
$$
$$\hat{\bm{D}}[i,i] = \bm{\Sigma}_{w_i,w_i}-\bm{\Sigma}_{N,w_i}^{\Tr}\bm{\Sigma}_{NN}^{-1}\bm{\Sigma}_{N,w_i}.
$$
Thus
\begin{align}\label{eq d diff}
|\bm{D}[i,i]-\hat{\bm{D}}[i,i]|
\le &  - \bm{\Sigma}_{N,w_i}^{\Tr}(\bm{\Phi}_{NN}-\bm{\Sigma}_{NN}^{-1})\bm{\Sigma}_{N,w_i} - 2\bm{\Sigma}_{N,w_i}^{\Tr}\bm{\Phi}_{NO}\bm{\Sigma}_{O,w_i} - \bm{\Sigma}_{O,w_i}^{\Tr}\bm{\Phi}_{OO}\bm{\Sigma}_{O,w_i} \nonumber\\
\le & |\bm{\Sigma}_{N,w_i}^{\Tr}\bm{\Sigma}_{NN}^{-1}\bm{\Sigma}_{NO}^{-1}\bm{\Phi}_{ON}\bm{\Sigma}_{N,w_i}| + 2|\bm{\Sigma}_{N,w_i}^{\Tr}\bm{\Phi}_{NO}\bm{\Sigma}_{O,w_i}| + |\bm{\Sigma}_{O,w_i}^{\Tr}\bm{\Phi}_{OO}\bm{\Sigma}_{O,w_i}|.
\end{align}
The term $\bm{\Phi}_{ON}\bm{\Sigma}_{N,w_i}$ appears multiple times in bounds (\ref{eq b diff}) and (\ref{eq d diff}). The next technical lemma shows it can be controlled by approximation radius $\rho$.
\begin{techlemma}\label{sigma nw}
For the submatrices $\bm{\Phi}_{ON}$ and $\bm{\Sigma}_{N,w_i}$ defined above, we have 
$$\|\bm{\Phi}_{ON}\bm{\Sigma}_{N,w_i}\|_2\le \frac{1}{v_r( \rho d^{-1/2} )} c_3 \{\phi_0(c_2/q)\}^{-1} v_{r-1}(c_3 \{\phi_0(c_2/q)\}^{-1} ).$$
\end{techlemma}
\begin{proof}
Applying Lemma \ref{norm-controlled} to the matrix $\bm{\Sigma}_{\nu(i),\nu(i)}$, we have that for any $l\in \mathbb{N}$,
\begin{align} \label{eq:Phi.ON}
\|\bm{\Phi}_{ON}\|_{v_r,l}
\le & \|(\bm{\Sigma}[\nu(i),\nu(i)])^{-1}\|_{v_r,l} \nonumber \\
\le & \|(\bm{\Sigma}[\nu(i),\nu(i)])^{-1}\|_2 v_{r-1}\left(\|(\bm{\Sigma}[\nu(i),\nu(i)])^{-1}\|_2 \|(\bm{\Sigma}[\nu(i),\nu(i)])\|_{v_r,l} \right) \nonumber \\
\le & \|\bm{\Sigma}_{\D\D}^{-1}\|_2 v_{r-1}\left(\|\bm{\Sigma}_{\D\D}^{-1}\|_2 \|\bm{\Sigma}_{\D\D}\|_{v_r,l} \right) \nonumber \\
\le & c_1 \{\phi_0(c_2/q)\}^{-1} q^d v_{r-1}\left(c_1 \{\phi_0(c_2/q)\}^{-1}
2^d(1+q/2)^{d-1} (1+ q/2 + d)
 \right),  
\end{align}
where the last step follows from Lemma \ref{D bounds}.

Now define a matrix operator $\tilde\nabla$ such that for all $ \bm{A}\in\mathbb{R}^{\mathcal{I}\times\mathcal{I}}$, $\tilde\nabla(\bm{A}) \in \mathbb{R}^{\mathcal{I}\times\mathcal{I}}$ has its $(i,j)$-entry defined as
$$\tilde\nabla(\bm{A})[i,j] = v_r(\|w_i-w_j\|_2 /d^{1/2}) \bm{A}[i,j].$$
Since for all $w_i,w_j\in\D$, $\sup_{1\le l\le d} |w_i[l] -w_j[l]|\ge \|w_i-w_j\|_2 /d^{1/2} $, we have that for the column of $\bm{\Phi}_{ON}$ corresponding to $w_k\in N$,
\begin{align*}
\|\tilde\nabla\bm{\Phi}_{O,w_k}\|_2
&= \left[ \sum_{w_j\in O} \{v_r(\|w_i-w_j\|_2 /d^{1/2} ) \bm{\Phi}_{w_j,w_k}\}^2 \right]^{1/2} \\
& \le \left[ \sum_{w_j\in O} \Big\{\sup_{1\le l\le d} v_r(|w_i[l] -w_j[l]|) \bm{\Phi}_{w_j,w_k}\Big\}^2 \right]^{1/2} \\
& \le \left[\sum_{w_j\in O} \Big\{\sum_{l=1}^d v_r(|w_i[l] -w_j[l]|) \bm{\Phi}_{w_j,w_k}\Big\}^2 \right]^{1/2} \\
& \le \sum_{l=1}^d \left[\sum_{w_j\in O} \{v_r(|w_i[l] -w_j[l]|) \bm{\Phi}_{w_j,w_k}\}^2 \right]^{1/2}  \\
& \le \sum_{1\le l\le d} \|\bm{\Phi}_{O,w_k}\|_{v_r,l} \le \sum_{1\le l\le d} \|\bm{\Phi}_{ON}\|_{v_r,l} \\
& \le d c_1 \{\phi_0(c_2/q)\}^{-1} q^d v_{r-1}\left(c_1 \{\phi_0(c_2/q)\}^{-1} 2^d(1+q/2)^{d-1}(1+ d+q/2) \right),
\end{align*}
where the last inequality follows from \eqref{eq:Phi.ON}.

Therefore, when $Z \in\mathscr{Z}_{v_r}$, we can derive that
\begin{align*}
\|\bm{\Phi}_{ON}\bm{\Sigma}_{N,w_i}\|_2
& \le \sum_{w_k\in N} \|\bm{\Phi}_{O,w_k} \|_2 K_0(w_k,w_i) \\
& \le \sum_{w_k\in N} \|\tilde\nabla\bm{\Phi}_{O,w_k}\|_2 \frac{1}{\inf_{w_j\in O} v_r(\|w_j-w_k\|_2 d^{-1/2} )} \\
&\quad \times \frac{1}{v_r(\|w_k-w_i\|_2) (1+\|w_k-w_i\|_2^{d+1})}  \\
& \le \sum_{w_k\in N} \|\tilde\nabla\bm{\Phi}_{O,w_k}\|_2 \frac{1}{\inf_{w_j\in O}v_r(\|w_j-w_i\|_2 d^{-1/2} )} \frac{1}{1+\|w_k-w_i\|_2^{d+1}} \\
& \le \sum_{w_k\in N} \|\tilde\nabla\bm{\Phi}_{O,w_k}\|_2 \frac{1}{v_r( \rho d^{-1/2} )} \frac{1}{1+\|w_k-w_i\|_2^{d+1}} \\
& \le  \frac{1}{v_r( \rho d^{-1/2} )} d2^dq^{-d} (1+q/2)^{d-1}(1+d+q/2) \cdot d c_1 \{\phi_0(c_2/q)\}^{-1} q^d \\
&\quad \times v_{r-1}\left(c_1 \{\phi_0(c_2/q)\}^{-1} 2^d(1+q/2)^{d-1}(1+ d+q/2) \right) \\
& \le \frac{1}{v_r( \rho d^{-1/2} )} c_3 \{\phi_0(c_2/q)\}^{-1} v_{r-1}(c_3 \{\phi_0(c_2/q)\}^{-1} ).
\end{align*}
where $c_3 = c_1 d^2 2^d(1+ d+q/2) (1+q/2)^{d-1}$, and the third inequality is due to the fact that $v_r$ is submultiplicative. When $q<1$, $c_3$ can be regarded as a constant independent of $q$.
\end{proof}

We now come back to the proof of the main Lemma \ref{l bounds} and first consider the situation where $q<1$. The elements of the column vector $\bm{\Sigma}_{O,w_i}$ are the covariances between locations that are at least $\rho$ distance apart. Therefore, we have
\begin{align}\label{eq sigma ow}
\|\bm{\Sigma}_{O,w_i}\|_2
\le & \sum_{j=1}^n K(w_i,w_j) \mathbbm{1}_{\{\|w_i-w_j\|_2\ge \rho\}} \nonumber\\
\le & \sum_{j=1}^n \frac{1}{v_r(\rho)}\frac{1}{1+\|w_i-w_j\|_2^{d+1}} \mathbbm{1}_{\{\|w_i-w_j\|_2\ge \rho\}} \nonumber \\
\le & \frac{\Gamma(d/2+1)2^d}{\pi^{d/2}q^d} \frac{1}{v_r(\rho)} \int_{\rho-q/2}^{+\infty} \frac{2\pi^{d/2}x^{d-1}}{\Gamma(d/2)} \frac{1}{1+(x-q/2)^{d+1}} dx \nonumber \\
\le & \frac{\pi d2^{d-1}q^{-d}}{v_r(\rho)} .
\end{align}

Using Technical Lemma \ref{sigma nw} and equation (\ref{eq sigma ow}) while controlling all other terms in $\|\bm{b}_i-\hat{\bm{b}}_i\|$ and $|\bm{D}[i,i]-\hat{\bm{D}}[i,i]|$ with the matrix $l_2$ norm, we have that
\begin{align}\label{eq b diff bound}
\|\bm{b_i}-\hat{\bm{b}}_i\|_2 
\le & (\|\bm{\Sigma}_{\D\D}\|_2\|\bm{\Sigma}_{\D\D}^{-1}\|_2+1)\|\bm{\Phi}_{ON}\bm{\Sigma}_{N,w_i}\|_2 + 2\|\bm{\Sigma}_{\D\D}^{-1}\|_2\|\bm{\Sigma}_{O,w_i}\|_2 \nonumber \\
\lesssim & 
\frac{1}{v_r( \rho d^{-1/2} )}q^{-d} \{\phi_0(c_2/q)\}^{-1} q^d c_3 \{\phi_0(c_2/q)\}^{-1} v_{r-1}(c_3 \{\phi_0(c_2/q)\}^{-1} ) \nonumber \\
& + \{\phi_0(c_2/q)\}^{-1} q^d q^{-d} v_r(\rho)^{-1}
\nonumber \\
\lesssim & \frac{1}{v_r( \rho d^{-1/2} )} \{\phi_0(c_2/q)\}^{-2} v_{r-1}(c_3\{\phi_0(c_2/q)\}^{-1}) ,
\end{align}
\begin{align}\label{eq d diff bound}
|\bm{D}[i,i]-\hat{\bm{D}}[i,i]|
\le & \|\bm{\Sigma}_{\D\D}\|_2^2\|\bm{\Sigma}_{\D\D}^{-1}\|_2 \|\bm{\Phi}_{ON}\bm{\Sigma}_{N,w_i}\|_2 + (\|\bm{\Sigma}_{\D\D}\|_2+1)\|\bm{\Sigma}_{\D\D}^{-1}\|_2\|\bm{\Sigma}_{O,w_i}\|_2 \nonumber \\
\lesssim & \frac{1}{v_r( \rho d^{-1/2} )}q^{-2d} \{\phi_0(c_2/q)\}^{-1} q^d c_3 \{\phi_0(c_2/q)\}^{-1} v_{r-1}(c_3 \{\phi_0(c_2/q)\}^{-1})\nonumber\\
& + q^{-d}\{\phi_0(c_2/q)\}^{-1} q^d v_r(\rho)^{-1} \nonumber \\
\lesssim & \frac{1}{v_r( \rho d^{-1/2} )} \{\phi_0(c_2/q)\}^{-2} q^{-d} v_{r-1}(c_3\{\phi_0(c_2/q)\}^{-1}) ,
\end{align}
where the multiplicative constants under the $\lesssim$ relations only depend on $d$. We also have the following bounds for $\bm{b}_i$,  $(\bm{D}[i,i])^{-1}$ and $(\hat{\bm{D}}[i,i])^{-1}$:
\begin{align}
& \|\bm{b}_i\|_2 = \|(\bm{\Sigma}[\nu(i),\nu(i)])^{-1}\bm{\Sigma}[\nu(i),i]\|_2
\lesssim \{\phi_0(c_2/q)\}^{-1} , \nonumber \\
& \min\big\{(\bm{D}[i,i])^{-1}, (\hat{\bm{D}}[i,i])^{-1}\big\} \ge K(w_i,w_i)^{-1}
= 1/K_0(0),  \label{eq d bound}  \\
& (\bm{D}[i,i])^{-1} \le \|(\bm{\Sigma}[\nu(i),\nu(i)])^{-1}\|_2
\lesssim \{\phi_0(c_2/q)\}^{-1} q^d ,  \nonumber \\
& (\hat{\bm{D}}[i,i])^{-1} \le \|(\bm{\Sigma}_{N\cup\{w_i\},N\cup\{w_i\}})^{-1}\|_2
\lesssim \{\phi_0(c_2/q)\}^{-1} q^d .  \nonumber  
\end{align}
Combining equations (\ref{eq b diff bound}) to (\ref{eq d bound}), when $q<1$, 
we can bound $\|\bm{l}_i-\hat{\bm{l}}_i\|_2$ as
\begin{align*}
 & \|\bm{l}_i-\hat{\bm{l}}_i\|_2 \\
={} & \|\bm{b}_i(\bm{D}[i,i])^{-1/2}-\hat{\bm{b}}(\hat{\bm{D}}[i,i])^{-1/2}\|_2 \\
\le {} & \|\bm{b}_i\|_2|(\bm{D}[i,i])^{-1/2}-(\hat{\bm{D}}[i,i])^{-1/2}| + \|\bm{b}_i-\hat{\bm{b}}_i\|_2 (\bm{D}[i,i])^{-1/2} \\ 
& + \|\bm{b}_i-\hat{\bm{b}}_i\|_2|(\bm{D}[i,i])^{-1/2}-(\hat{\bm{D}}[i,i])^{-1/2}| \\
\lesssim {} & \|\bm{b}_i\|_2 \Big( [\{\phi_0(c_2/q)\}^{-1} q^d]^{3/2}|\bm D[i,i]-\hat{\bm{D}}[i,i]| \\
&\quad +[\{\phi_0(c_2/q)\}^{-1} q^d]^{5/2} O\big((\bm D[i,i]-\hat{\bm{D}}[i,i])^2\big) \Big) \\
& + \|\bm{b}_i-\hat{\bm{b}}_i\|_2 (\bm{D}[i,i])^{-1/2}  \\
& + \|\bm{b}_i-\hat{\bm{b}}_i\|_2 \Big( [\{\phi_0(c_2/q)\}^{-1} q^d]^{3/2}|\bm D[i,i]-\hat{\bm{D}}[i,i]| \\
&\quad +[\{\phi_0(c_2/q)\}^{-1} q^d]^{5/2} O\big((\bm D[i,i]-\hat{\bm{D}}[i,i])^2\big) \Big) \\ 
\lesssim {} & \frac{1}{v_r(\rho d^{-1/2} )}\{\phi_0(c_2/q)\}^{-\frac{9}{2}}q^{\frac{1}{2}d} v_{r-1}(c_3\{\phi_0(c_2/q)\}^{-1})  \\
& + \frac{1}{v_r(\rho d^{-1/2} )}\{\phi_0(c_2/q)\}^{-\frac{5}{2}}q^{\frac{1}{2}d} v_{r-1}(c_3\{\phi_0(c_2/q)\}^{-1})  \\
\lesssim {} & \frac{1}{v_r(\rho d^{-1/2} )}\{\phi_0(c_2/q)\}^{-\frac{9}{2}}q^{\frac{1}{2}d} v_{r-1}(c_3\{\phi_0(c_2/q)\}^{-1}).
\end{align*}
This completes the proof of \eqref{eq:lemS3.l.diff1}.

If $q\ge 1$, both the minimal and maximal eigenvalues of $\bm{\Sigma}_{\D\D}$ are bounded by constants by \eqref{eq d bound}. Therefore, all terms involving $q$ in \eqref{eq:lemS3.l.diff1} become constant and we have
$$\|\bm{l}_i-\hat{\bm{l}}_i\|_2 \lesssim \frac{1}{v_r(\rho d^{-1/2} )} .$$ 
Since the above arguments hold for arbitrary $i$, we finish the proof of \eqref{eq:lemS3.l.diff2}.
\vspace{2mm}

\underline{\textsc{Part 2}}
We now derive a sufficient condition for $\|\hat{L}-L\|_2\le \|\bm{\Sigma}_{\D\D}\|_2^{-1/2}/2$. The left hand side can be bounded as
\begin{align*}
\|\hat{L}-L\|_2 \le n^{1/2} \max_i\|\hat l_i - l_i\|_2.
\end{align*}
We first consider the case $0<q<1$. 
A sufficient condition for $\|\hat{L}-L\|_2\le \|\bm{\Sigma}_{\D\D}\|_2^{-1/2}/2$ is
\begin{equation}\label{eq sufficient 1}
2 \|\bm{\Sigma}_{\D\D}\|^{1/2}_2 n^{1/2} \max_i \|\hat{\bm{l}}_i-\bm{l}_i\|_2 \le 1.
\end{equation} 
By applying Lemma \ref{D bounds} and part 1 of this proof to equation (\ref{eq sufficient 1}), we get the following sufficient condition 
\begin{equation}\label{eq rho q cond 2}
\frac{n^{1/2}}{v_r(\rho d^{-1/2})}\{\phi_0(c_2/q)\}^{-9/2} v_{r-1}(c_3\{\phi_0(c_2/q)\}^{-1})  \le c_6 ,
\end{equation}
for some constant $c_6$ only dependent on $d$.

Now for the case $q>1$, all terms involving $q$ can be considered as constants, leaving $n$ and $\rho$ as the only variables. Therefore, a sufficient condition for $\|\hat{L}-L\|_2\le \|\bm{\Sigma}_{\D\D}\|_2^{-1/2}/2$ is
$n^{1/2}\{v_r(\rho d^{-1/2})\}^{-1}\le c'_6$
for some constant $c'_6$ only dependent on $d$.
\end{proof}
\vspace{6mm}

For the polynomial decaying class $\mathscr{Z}_{c_r}$, we have a similar result.
\vspace{2mm}

\begin{lemma}\label{l bounds poly}
Suppose that $\bm{Z}\in\mathscr{Z}_{c_r}$ for some $r\ge d+1$. If $0<q<1$ and
\begin{equation*}
\frac{ n^{1/2} }{ (1+ \rho d^{-1/2} )^{r-d-1}} q^{(r-7)d} \{\phi_0(c_2/q)\}^{-(r+4)} (c_1c_5d2^{d-1}\pi/\sqrt{6})^r \le c_7
\end{equation*}
for some constant $c_7$ only dependent on $d$, then equation (\ref{eq phi prelim}) holds and
$$
\|\bm{l}_i-\hat{\bm{l}}_i\|_2 \lesssim \frac{1}{ (1+ \rho d^{-1/2} )^{r-d-1}} q^{(r-13/2)d} \{\phi_0(c_2/q)\}^{-(r+4)} (c_1c_5d2^{d-1}\pi/\sqrt{6})^r .$$
Else if $q\ge 1$ and $ n^{1/2}(1+ \rho d^{-1/2} )^{-(r-d-1)} \{\phi_0(c_2/q)\}^{-r} (c_1c_5d2^{d-1}\pi/\sqrt{6})^r \le c'_7$ for some constant $c'_7$ only dependent on $d$, then equation (\ref{eq phi prelim}) still holds and
$$\|\bm{l}_i-\hat{\bm{l}}_i\|_2\lesssim \frac{1}{ (1+ \rho d^{-1/2} )^{r-d-1}} \{\phi_0(c_2/q)\}^{-r} \{c_1c_5d2^{d-1}\pi/\sqrt{6}\}^r.$$
\end{lemma}

\begin{proof}\mbox{}
\underline{\textsc{Part 1}}
We first prove the upper bounds for $\max_i\|\bm{l}-\hat{\bm{l}}\|_2$. The proof is the same as that of Lemma \ref{l bounds} till equation (\ref{eq d diff}). We begin with the bounds on $\|\bm{\Sigma}_{O,w_i}\|_2$ and $\|\bm{\Phi}_{ON}\bm{\Sigma}_{N,w_i}\|_2$.
\begin{align*}
\|\bm{\Sigma}_{O,w_i}\|_2
\le & \left[\sum_{j=1}^n \{K(w_i,w_j)\}^2 \mathbbm{1}_{\{\|w_i-w_j\|_2\ge \rho\}} \right]^{1/2} \nonumber\\
\le & \left[\sum_{j=1}^n \frac{1}{(1+\|w_i-w_j\|_2)^{2r}} \mathbbm{1}_{\{\|w_i-w_j\|_2\ge \rho\}} \right]^{1/2}\nonumber\\
\le & \frac{1}{(1+\rho)^{r-(d+1)/2}} \left[\sum_{j=1}^n \frac{1}{(1+\|w_i-w_j\|_2)^{d+1}}  \mathbbm{1}_{\{\|w_i-w_j\|_2\ge \rho\}} \right]^{1/2}\nonumber\\
\le & \frac{1}{(1+\rho)^{r-(d+1)/2}} \left[\frac{\Gamma(d/2+1)2^d}{\pi^{d/2}q^d} \int_{\rho-q/2}^{+\infty} \frac{2\pi^{d/2}x^{d-1}}{\Gamma(d/2)} \frac{1}{1+(x-q/2)^{d+1}} dx \right]^{1/2} \nonumber \\
\le & \frac{ (\pi d2^{d-1}q^{-d})^{1/2}} {(1+\rho)^{r-(d+1)/2}}.
\end{align*}
To bound $\|\bm{\Phi}_{ON}\bm{\Sigma}_{N,w_i}\|_2$, we first consider the norm for columns of $\|\bm{\Phi}_{ON}\|_2$. Define an operator $c_r\odot$ such that for any matrix $\bm{A}\in\mathbb{R}^{\mathcal{I}\times\mathcal{I}}$, $c\odot\bm{A} \in \mathbb{R}^{\mathcal{I}\times\mathcal{I}}$ has its $(i,j)$-entry defined as
$$c\odot\bm{A}[i,j] = (1+\|w_i-w_j\|_2 d^{-1/2})^{r-(d+1)/2}\bm{A}[i,j].$$
Since for all $w_i,w_j\in\D$, $\|w_i-w_j\|_\infty\ge \|w_i-w_j\|_2 d^{-1/2}$, we have for all $w_k\in N$, for the column of $\bm{\Phi}_{ON}$ corresponding to $w_k\in N$,
\begin{align*}
\|c\odot\bm{\Phi}_{O,w_k}\|_2 
\le & [\max\{1,\lceil 1/q\rceil\}]^{d/2} 
\|\bm{\Phi}_{O,w_k}\|_{r-(d+1)/2}\\
\le & [\max\{1,\lceil 1/q\rceil\}]^{d/2}  \|(\bm{\Sigma}[\nu(i),\nu(i)])^{-1}\|_{r-(d+1)/2} \\
\overset{(i)}{\le} & [\max\{1,\lceil 1/q\rceil\}]^{\frac{13}{2}d} c_4c_5^{r-(d+1)/2} \|(\bm{\Sigma}[\nu(i),\nu(i)])^{-1}\|_2 \\
& \{\|(\bm{\Sigma}[\nu(i),\nu(i)])^{-1}\|_2 \|(\bm{\Sigma}[\nu(i),\nu(i)])\|_{c_{r-(d+1)/2}}\}^{r-1/2}\\
\le & [\max\{1,\lceil 1/q\rceil\}]^{\frac{13}{2}d} c_4c_5^{r-(d+1)/2} \|\bm{\Sigma}_{\D\D}^{-1}\|_2 (\|\bm{\Sigma}_{\D\D}^{-1}\|_2 \|\bm{\Sigma}_{\D\D}\|_{c_{r-(d+1)/2}})^{r-1/2}\\
\overset{(ii)}{\le} & c_4 c_5^{r-(d+1)/2} [\max\{1,\lceil 1/q\rceil\}]^{\frac{13}{2}d} [c_1\{\phi_0(c_2/q)\}^{-1}q^d]^{r+1/2} (d2^{d-1}\pi/\sqrt{6})^{r-1/2},
\end{align*} 
where $(i)$ follows from Lemma \ref{norm-controlled}, and $(ii)$ follows from (5) of Lemma \ref{D bounds}.

We first consider the case $q<1$. We have
$$\|c\odot\bm{\Phi}_{O,w_k}\|_2 \lesssim q^{(r-6)d} \{\phi_0(c_2/q)\}^{-(r+1/2)} (c_1c_5d2^{d-1}\pi/\sqrt{6})^r ,$$
where the constant in $\lesssim$ is only dependent on $d$. Therefore we have
\begin{align*}
& \|\bm{\Phi}_{ON}\bm{\Sigma}_{N,w_i}\|_2 \\
\le{} & \sum_{w_k\in N} \|\bm{\Phi}_{O,w_k} \|_2 K_0(w_k,w_i) \\
\le{} & \sum_{w_k\in N} \|c\odot\bm{\Phi}_{O,w_k}\|_2 \frac{1}{\inf_{w_j\in O} (1+\|w_j-w_k\|_2 d^{-1/2} )^{r-(d+1)/2}} \frac{1}{(1+\|w_i-w_k\|_2)^{r}}  \\
\le & \sum_{w_k\in N} \|c\odot\bm{\Phi}_{O,w_k}\|_2 \frac{1}{\inf_{w_j\in O} (1+\|w_j-w_k\|_2 d^{-1/2} +\|w_i-w_k\|_2)^{r-(d+1)/2}} \\
&\quad \times \frac{1}{(1+\|w_i-w_k\|_2)^{(d+1)/2}} \\
\le{} & \sum_{w_k\in N} \|c\odot\bm{\Phi}_{O,w_k}\|_2 \frac{1}{ (1+ \rho d^{-1/2} +\|w_i-w_k\|_2)^{r-(d+1)/2}} \frac{1}{(1+\|w_i-w_k\|_2)^{(d+1)/2}}  \\
\le{} & \sum_{w_k\in N} \|c\odot\bm{\Phi}_{O,w_k}\|_2 \frac{1}{ (1+ \rho d^{-1/2} )^{r-d-1}} \frac{1}{(1+\|w_i-w_k\|_2)^{d+1}}  \\
\le{} & \pi d 2^{d-1} q^{-d} \sup_{w_k\in N} \|c\odot\bm{\Phi}_{O,w_k}\|_2 \frac{1}{ (1+ \rho d^{-1/2} )^{r-d-1}} \\
\lesssim{} & \frac{1}{ (1+ \rho d^{-1/2} )^{r-d-1}} q^{-d}
q^{(r-6)d} \{\phi_0(c_2/q)\}^{-(r+1/2)} (c_1c_5d2^{d-1}\pi/\sqrt{6})^r  \\
\lesssim{} & \frac{1}{ (1+ \rho d^{-1/2} )^{r-d-1}} q^{(r-7)d} \{\phi_0(c_2/q)\}^{-(r+1/2)} (c_1c_5d2^{d-1}\pi/\sqrt{6})^r .
\end{align*}
The rest of the proof follows a similar strategy as that of Lemma \ref{l bounds}. We now only list the key steps.

\begin{align*}
\|\bm{b_i}-\hat{\bm{b}}_i\|_2 
\le & (\|\bm{\Sigma}_{\D\D}\|_2\|\bm{\Sigma}_{\D\D}^{-1}\|_2+1)\|\bm{\Phi}_{ON}\bm{\Sigma}_{N,w_i}\|_2 + 2\|\bm{\Sigma}_{\D\D}^{-1}\|_2\|\bm{\Sigma}_{O,w_i}\|_2 \\
\lesssim & \frac{1}{ (1+ \rho d^{-1/2} )^{r-d-1}} q^{(r-7)d} \{\phi_0(c_2/q)\}^{-(r+3/2)} (c_1c_5d2^{d-1}\pi/\sqrt{6})^r ,
\end{align*}
\begin{align}\label{eq d diff poly}
|\bm{D}[i,i]-\hat{\bm{D}}[i,i]|
\le & \|\bm{\Sigma}_{\D\D}\|_2^2\|\bm{\Sigma}_{\D\D}^{-1}\|_2 \|\bm{\Phi}_{ON}\bm{\Sigma}_{N,w_i}\|_2 + (\|\bm{\Sigma}_{\D\D}\|_2+1)\|\bm{\Sigma}_{\D\D}^{-1}\|_2\|\bm{\Sigma}_{O,w_i}\|_2 \nonumber\\
\lesssim & \frac{1}{ (1+ \rho d^{-1/2} )^{r-d-1}} q^{(r-8)d} \{\phi_0(c_2/q)\}^{-(r+3/2)} (c_1c_5d2^{d-1}\pi/\sqrt{6})^r\nonumber .
\end{align}
\begin{align*}
\|\bm{l}_i-\hat{\bm{l}}_i\|_2
\le & \|\bm{b}\|_2|(\bm{D}[i,i])^{-1/2}-(\hat{\bm{D}}[i,i])^{-1/2}| + \|\bm{b}_2-\hat{\bm{b}}_2\|_2 (\bm{D}[i,i])^{-1/2} \\
& + \|\bm{b}_2-\hat{\bm{b}}_2\|_2|(\bm{D}[i,i])^{-1/2}-(\hat{\bm{D}}[i,i])^{-1/2}| \\
\lesssim & \frac{1}{ (1+ \rho d^{-1/2} )^{r-d-1}} q^{(r-13/2)d} \{\phi_0(c_2/q)\}^{-(r+4)} (c_1c_5d2^{d-1}\pi/\sqrt{6})^r .
\end{align*}
If $q\ge 1$, then we have a simplified bound as
$$\|\bm{l}_i-\hat{\bm{l}}_i\|_2 \lesssim \frac{1}{ (1+ \rho d^{-1/2} )^{r-d-1}} \{\phi_0(c_2/q)\}^{-r} (c_1c_5d2^{d-1}\pi/\sqrt{6})^r .$$

\underline{\textsc{Part 2}}
A sufficient condition for $\|\hat{L}-L\|_2 \le \|\bm{\Sigma}_{\D\D}\|_2^{-1/2}/2$ can be derived by satisfying equation (\ref{eq sufficient 1}) in the proof of Lemma \ref{l bounds}. When $q<1$, this yields
$$ \frac{ n^{1/2} }{ (1+ \rho d^{-1/2} )^{r-d-1}} q^{(r-7)d} \{\phi_0(c_2/q)\}^{-(r+4)} (c_1c_5d2^{d-1}\pi/\sqrt{6})^r \le c_7$$
for some constant $c_7$ only dependent on $d$.
When $q\ge 1$, this yields
$$ \frac{n^{1/2}}{ (1+ \rho d^{-1/2} )^{r-d-1}} \{\phi_0(c_2/q)\}^{-r} (c_1c_5d2^{d-1}\pi/\sqrt{6})^r \le c'_7$$
for some constant $c'_7$ only dependent on $d$.
\end{proof}

\section{Proof of Theorems and Corollary in Section 3 of the Main Text}
\label{sec:theory proof}

\subsection{Proof of Theorems \ref{W22 rate} and \ref{W22 rate poly}}
\begin{proof}[of Theorem \ref{W22 rate}]
We first consider the case $0<q<1$. Since the second term in Lemma \ref{W2 decom} of the main text is dominated by the first term as shown in the proof of Lemma \ref{l bounds}, we plug the results of Lemmas \ref{D bounds} and \ref{l bounds} into Lemma \ref{W2 decom} to obtain that 
\begin{align*}
&\quad~ W_2^2(Z_{\D},\hat{Z}_{\D}) \\
& \lesssim  n\|\bm{\Sigma}_{\D\D}\|_2^2 (\|\bm{\Sigma}_{\D\D}^{-1}\|_2)^{-1/2} \max_i\|\bm{l}_i-\hat{\bm{l}}_i\|_2  \\ 
& \lesssim  q^{-2d} \{\phi_0(c_2/q)\}^{-\frac{1}{2}}q^{\frac{d}{2}} \cdot  \frac{1}{v_r( \rho d^{-1/2} )}\{\phi_0(c_2/q)\}^{-\frac{9}{2}}q^{\frac{1}{2}d} v_{r-1}(c_3\{\phi_0(c_2/q)\}^{-1})  \\
& \lesssim \frac{n}{v_r( \rho d^{-1/2} )}\{\phi_0(c_2/q)\}^{-5} q^{-d} v_{r-1}(c_3\{\phi_0(c_2/q)\}^{-1}) .
\end{align*}
For the case when $q\ge 1$, all terms involving $q$ can be considered as constants, thus we have
$$W_2^2(Z_{\D},\hat{Z}_{\D}) \lesssim \frac{n}{v_r( \rho d^{-1/2} )}.$$
Since the above results hold for all $Z\in\mathscr{Z}_{v_r}$, we finish the proof of Theorem \ref{W22 rate}.
\end{proof}
The proof of Theorem \ref{W22 rate poly} follows the same steps as the proof of Theorem \ref{W22 rate}, with Lemma \ref{l bounds} replaced by Lemma \ref{l bounds poly}.

\subsection{Proof of Corollary \ref{corollary}}
We first prove a technical lemma.
\vspace{2mm}
\begin{techlemma}\label{vr connect}
For all $x>0$, we have 
$$v_2(x^3) \ge \{v_1(x)\}^3.$$
\end{techlemma}
\begin{proof}
\begin{align*}
v_2(x^3)  & =\sum_{k=0}^{\infty} \frac{x^{3k}}{(k!)^{3/2}} \frac{1}{(k!)^{1/2}} \\
& = \sum_{j=0}^{\infty} \frac{1}{(j!)^{1/2}} \sum_{k=0}^{\infty}  \left\{\frac{x^k}{(k!)^{1/2}}\right\}^3 \frac{\frac{1}{(k!)^{1/2}}}{\sum_{j=0}^{\infty} \frac{1}{(j!)^{1/2}}} \\
& \overset{\mr{Jensen}}{\ge} \sum_{j=0}^{\infty} \frac{1}{(j!)^{1/2}} \left\{ \sum_{k=0}^{\infty} \frac{x^k}{(k!)^{1/2}} \frac{1}{(k!)^{1/2}} \right\}^3\\
& \ge [v_1(x)]^3.
\end{align*}
\mbox{}
\end{proof}

\vspace{3mm}

\begin{proof}
We now prove the results for the three covariance functions. 

\underline{\textsc{Mat\'{e}rn}} Let $r=2$. Since the Mat\'{e}rn covariance function decays at the rate $K_0(\|s_1-s_2\|_2) = O\big(\|s_1-s_2\|_2^{\nu-1/2}\exp(-\phi\|s_1-s_2\|) \big) \lesssim  v_2(\phi\|s_1-s_2\|_2)$, we have $\bm{Z}\in\mathscr{Z}_{v_2}$ up to a scale parameter $\phi$. Thus when $q<1/\phi$, Theorem \ref{W22 rate} in the main paper gives
\begin{equation}\label{eq thm1 sim}
W_2^2(\bm{Z}_{\D},\hat{\bm{Z}}_{\D}) \lesssim \frac{1}{v_2(\phi \rho d^{-1/2} )}\exp\left(c_3 \{\phi_0(c_2/(\phi q))\}^{-1} + \ln\Big[n q^{-d} \{\phi_0(c_2/(\phi q))\}^{-5}\Big] \right) ,
\end{equation}
where $\phi_0(\cdot)$ is the Fourier transform of Mat\'{e}rn covariance function with scale parameter $\phi=1$, which has the closed form solution
$$\hat{K}_0(\xi) = \frac{\tau^2 2^d \pi^{d/2} \Gamma(\nu+d/2) }{\Gamma(\nu)} (1 + \|\xi\|_2^2)^{-(\nu+d/2)}.$$
Thus, the function $[\phi_0(c_2/(\phi q))]^{-1}$ has closed form as
\begin{align}\label{eq phi0}
[\phi_0(c_2/(\phi q))]^{-1} & = \frac{1}{\inf_{\|\xi\|_2\le 2c_2/q}\hat{K}_0(\xi)} \nonumber\\
& = \frac{\Gamma(\nu) \{1+4c_2^2/(\phi^2q^2)\}^{\nu+d/2}}{\tau^2 2^d \pi^{d/2} \Gamma(\nu+d/2) } = c_{m,1} \{1+4c_2^2/(\phi^2q^2)\}^{\nu+d/2},
\end{align}
where $c_{m,1}$ is a constant dependent only on $d,\nu,\tau^2$.
Therefore, let $\rho = \frac{d^{1/2}}{\phi}\left[c_3c_{m,1} \Big(1 + \frac{ 4c_2^2}{\phi^2 q^2}\Big)^{\nu+\frac{d}{2}} + \ln\Big\{ c_{m,1}n q^{-d} \Big(1 + \frac{ 4c_2^2}{\phi^2q^2}\Big)^{5(\nu+\frac{d}{2})} \Big\} \right]^3$. We plug the results of equation (\ref{eq phi0}) and Technical Lemma \ref{vr connect} into equation (\ref{eq thm1 sim}) to obtain that
\begin{align*}
 W_2^2(\bm{Z}_{\D},\hat{\bm{Z}}_{\D})
\lesssim & \frac{\exp\left[c_3 c_{m,1} \Big(1 + \frac{ c_2^2}{\phi^2q^2}\Big)^{\nu+\frac{d}{2}} + \ln\Big\{ c_{m,1}n q^{-d} \Big(1 + \frac{ c_2^2}{\phi^2q^2}\Big)^{5(\nu+\frac{d}{2})} \Big\} \right]}
{\left(\exp\left[c_3 c_{m,1} \Big(1 + \frac{ c_2^2}{\phi^2q^2}\Big)^{\nu+\frac{d}{2}} + \ln\Big\{ c_{m,1}n q^{-d} \Big(1 + \frac{ c_2^2}{\phi^2q^2}\Big)^{5(\nu+\frac{d}{2})} \Big\} \right]\right)^3} \\
= & \frac{1}
{\left(\exp\left[c_3 c_{m,1} \Big(1 + \frac{ c_2^2}{\phi^2q^2}\Big)^{\nu+\frac{d}{2}} + \ln\Big\{ c_{m,1}n q^{-d} \Big(1 + \frac{ c_2^2}{\phi^2q^2}\Big)^{5(\nu+\frac{d}{2})} \Big\} \right]\right)^2} \\
= & o(1)
\end{align*}
\vspace{2mm}

\underline{\textsc{Gaussian}}
The Gaussian covariance function decays even faster than Mat\'ern as the spatial distance increases. Therefore, we can let $r=2$ and we have $\bm{Z}\in\mathscr{Z}_{v_2}$ up to a scale parameter $a^{1/2}$. The Fourier transform of the $d$-dimensional Gaussian function $K_0$ under the case $a=1$ is
$$\hat{K}_0(\xi) = \tau^2 \exp(-\|\xi\|_2^2/4).$$
Thus the function $\{\phi_0(c_2/(a^{1/2}q))\}^{-1}$ has the closed form solution as
\begin{equation}\label{eq phi0 exp}
\{\phi_0(c_2/(a^{1/2}q))\}^{-1} =  \tau^{-2} \exp\left(\frac{c_2^2}{aq^2}\right).
\end{equation}
Similar to the case of Mat\'ern, we combine equation (\ref{eq phi0 exp}) and equation (\ref{eq thm1 sim}) to derive the condition on $\rho$ as in Corollary \ref{corollary}.
\vspace{3mm}

\underline{\textsc{Generalized Cauchy}}
By definition, $K_0(\|s_1-s_2\|_2) \lesssim 1/(1+\|s_1-s_2\|_2)^\lambda$. Thus $\bm{Z}\in\mathscr{Z}_{c_r}$ and Theorem 2 can be directly applied here. By Theorem 1 of \cite{bevilacqua2019estimation}, we have
$$ \{\phi_0(c_2\phi /q) \}^{-1} \lesssim (2c_2\phi/q)^{d+\delta} ,$$
where the multiplicative constant in the $\lesssim$ relation only depends on $d,\tau^2, \lambda,\delta,\phi$. Therefore
\begin{align*}
W_2^2(\bm{Z}_{\D},\hat{\bm{Z}}_{\D})
\lesssim &  \frac{1}{(1+\rho/(\phi d^{1/2}))^{\lambda-d-1}} nq^{(\lambda-8)d} (2c_2\phi/q)^{-(d+\delta)(\lambda+9/2)} (c_1c_5d2^{d-1}\pi/\sqrt{6})^{\lambda}  \\
&  \lesssim \frac{nq^{-d(\lambda+9/2+8-\lambda)-\delta(\lambda+9/2)}}{\{1+\rho/(\phi d^{1/2})\}^{\lambda-d-1}}  \\
&  \lesssim \frac{n q^{-\frac{25}{2} d - \delta(\lambda+9/2)}}{\{1+\rho/(\phi d^{1/2})\}^{\lambda-d-1}} .
\end{align*}
Setting the right-hand side to be $o(1)$ and reversely solving for $\rho$ gives the condition on $\rho$ in Corollary \ref{corollary}.
\end{proof}

\section{Posterior Sampling Algorithms for RadGP Regression}
\label{sec:algo.radgp}
\subsection{Algorithm and Computational Complexity for Latent Effects Model}
We provide the algorithm to perform the posterior sampling on the latent effects model described in Section 4 of the main paper. 
\begin{algorithm}[ht]
\SetAlgoLined
Input training locations $\mathcal{T}_1\subset \Omega$, test locations $\mathcal{T}_2\subset\Omega$, covariates $\bm{X}_{\T_1}$, $\bm{X}_{\T_2}$ and response $\bm{Y}_{\T_1}$. \\
Set an approximation radius $\rho$, the total number of MCMC steps $l_1$ and the number of burn-in steps $l_2$ ($l_2<l_1$). Set initial values for $\beta$, $\theta$, $\sigma$ and $Z_{\T_1}$.\\
\For{$1\le l \le l_1$}{
Sample $\beta$ from equation (\ref{eq beta post}) in the main paper; \\
Sample $\sigma$ from equation (\ref{eq sigma post}) in the main paper;\\
Generate random vector ${W}\sim N(\sigma^{-2}({Y}_{\T_1}-{X}_{\T_1}\beta), \hat{\Phi}+{I}_n/\sigma^2)$ by 
${W} = \sigma^{-2}({Y}_{\T_1}-{X}_{\T_1}\beta) + \hat{{L}}{W}_1 + \sigma^{-1} {W}_2,$   
where ${W}_1, {W}_2$ are independent $N(0,{I}_n)$ random vectors and $\hat{L}$ is the Cholesky factor of the precision matrix $\hat{\Phi}$; \\
Compute a sample of $\hat Z_{\T_1}$ by solving the linear system $(\hat{\Phi}+{I}_n/\sigma^2)\hat Z_{\T_1} = {W}$, e.g., using conjugate gradient;\\
Update $\theta$ using a Metropolis Hastings sampling step; \\ 
\If{$l\ge l_2$}{Sample $\hat{Z}_{\T_2}$ according to equation (\ref{eq ZT}).}
}
\KwResult{Output the posterior samples of $\beta, \theta, \delta, \hat{Z}_{\mathcal{T}_1}$ and $\hat{Z}_{\mathcal{T}_2}$.}
\caption{Posterior Sampling for Latent RadGP Regression}\label{alg cg}
\end{algorithm}

We proceed to analyze the computational complexity of Algorithm \ref{alg cg}. Since the number of MCMC steps $l_1$ is typical of order $10^3\sim 10^5$, the cost of operations performed outside MCMC iterations (such as the construction of radial neighbors graph) is negligible. Hereafter we only analyze the operations inside the MCMC loop. 

Recall that $n$ denotes the sample size of training set $\T_1$. Further let $n_{\mr{test}}$ be the sample size of test set, and $d_\beta$ be the dimension of coefficients $\beta$. Sampling $d_\beta$ from equation (\ref{eq beta post}) in the main paper involves matrix inversion of dimension $d_\beta$ and multiplication between $d_\beta\times n$ and $n\times d_\beta$ matrices, which takes $O(d_\beta^3+nd_\beta^2)$ time. Once the mean and variance of $\beta$ are computed, sampling $\beta$ itself only takes $O(\beta)$ time.
Assuming $d_\beta<n$, the time complexity of sampling $\beta$ is $O(nd_\beta^2)$.

Sampling $\sigma^2$ from equation (\ref{eq sigma post}) in the main paper involves the computation of linear regression residuals for $n$ training samples, which takes $O(n)$ time. 

Next, we consider the computational complexity of sampling the spatial random effects $\hat Z_{\T_1}$ at the training locations. Let $M_1$ be the maximal number of parents in the radial neighbors graph on $\T_1$. Computing the vector $W$ involves the computation of Cholesky factor $\hat{L}$, which takes no more than $O(nM_1^3)$ time, and the matrix-vector multiplication between $\hat{L}$ and $W_1$, which takes $O(nM_1)$ time. The overall complexity of sampling $W$ is no greater than $O(n M_1^3)$.

Let $M_2$ be the average number of nonzero elements per column in the precision matrix $\hat{\Phi}$. For RadGP, $M_2$ is equivalent to the number of locations within $2\rho$ radius, averaged at all training locations. If we solve the linear system $(\hat{\Phi}+{I}_n/\sigma^2)\hat Z_{\T_1} = {W}$ with conjugate gradients, the most costly operation in each iteration is the multiplication between the matrix $\hat{\Phi}+I_n/\sigma^2$ and a vector, which takes $O(n M_2)$ time. Let $n_{\text{cg}}$ be the average number of conjugate gradient steps. In theory, $n_{\text{cg}}$ can be as large as $n$; in practice, under a fixed tolerance ($10^{-6}$ in our package),  $n_{\text{cg}}$ can be much smaller than $n$. See Section \ref{subsec:MCMC} for some simulation studies.
The total complexity of solving this linear system is $O(n_{\text{cg}} n M_2)$.

Let $M_3$ be the maximal number of parents in the radial neighbors graph on $\T_1\cup \T_2$, the union of training and test sets. Sampling random effects on test locations is equivalent to sampling $n_{\mr{test}}$ unidimensional Gaussian distributions, which has the cost $O(n_{\mr{test}} M_3^3)$.

Combining all the sampling steps, the overall time complexity of Algorithm \ref{alg cg} is 
\begin{align} \label{eq:complexity1}
& O\left(l_1 n (d_\beta^2+M_1^3+n_{\text{cg}}M_2) + l_2 n_{\mr{test}} M_3^3 \right).
\end{align}
In practice, $n$ is usually much larger than $M_1, M_2$ and $M_3$ to the point that $M_i, 1\le i \le 3$ can be regarded as constants compared to $n$. Assuming $n_{\mr{test}}$ is in the same order or smaller than $n$, then the most computationally expensive step is to sample the spatial random effects at the training locations because $n_{cg}$ is the only variable (other than $n$ itself) that increases with the sample size $n$. Therefore, the overall computational complexity in \eqref{eq:complexity1} can be simplified to $O(M_2 l_1 n n_{\text{cg}})$.

If the spatial locations are distributed on a grid in $\mathbb{R}^d$ with minimal separation distance $q$, then $M_2 \approx (\sqrt{\pi}\rho/q)^d \Gamma^{-1}(d/2+1) $. The overall computational complexity can be further simplified to $O( l_1 n n_{\text{cg}}(\rho/q)^d )$. We can see the computational complexity is proportional to the $d$th power of approximation radius $\rho$, with smaller $\rho$ leading to faster computation.

\subsection{Posterior Sampling for Response Model}
Algorithm \ref{alg cg} outputs posterior samples of all spatial random effects along with posterior samples of parameters $\beta, \theta, \delta$. Because high dimensionality of the spatial random effects may negatively impact the mixing of MCMC chains, we directly approximate the marginal covariance using radial neighbors without estimating these latent effects. The resulting model of the response \citep{Finetal19} only involves a low-dimensional parameter vector to be estimated via MCMC. Let $\tilde{\bm{Z}}_{\T_1} = \bm{Z}_{\T_1} +\epsilon_{\T_1}$ and $\tilde{\bm{Z}}_{\T_2} = \bm{Z}_{\T_2} +\epsilon_{\T_2}$. We use $\tilde{\bm{\Phi}}$ to denote the RadGP precision that approximates $(\bm{\Sigma}_{\T_1\T_1}+\sigma^2 \bm{I})^{-1}$. The joint posterior density now becomes
$$ \mr{det}(\tilde{\bm{\Phi}})^{1/2}\exp\left\{-\frac{1}{2} (\bm{Y}_{\T_1}-\bm{X}_{\T_1}\beta)^{\Tr} \tilde{\bm{\Phi}} (\bm{Y}_{\T_1}-\bm{X}_{\T_1}\beta)\right\}p(\theta)p(\beta)p(\sigma^2), 
 $$
where, using $\hat{\bm{B}}$ and $\hat{\bm{D}}$ as defined in Lemma \ref{lem chol pre}, we have $\tilde{\bm{\Phi}} = (\bm{I}_n-\hat{\bm{B}}^{\Tr})\hat{\bm{D}}^{-1}(\bm{I}_n-\hat{\bm{B}})$,
which is a sparse matrix. Let $\tilde M = \sup_{1\leq i\leq n} |N(i)|$ be the maximal number of points in a radius $\rho$ ball. By the derivation in Lemma \ref{lem chol pre}, each row of $\hat{\bm{B}}$ has at most $\tilde M$ nonzero elements, indicating that the computation of all rows of $\hat{\bm B}$ can proceed in parallel with a total computational complexity $O(n\tilde{M}^3)$. Computations of the quadratic form $\bm{X}_\D^{\Tr}\tilde{\bm{\Phi}}\bm{X}$ and the determinant $\mr{det}(\tilde{\bm{\Phi}})$ have $O(n\tilde{M}^2)$ and $O(n)$ time complexity, respectively.

Posterior sampling of unknown parameters proceeds as a hybrid Gibbs adaptive Metropolis-Hastings sampler. If the prior for $\beta$ is normal, $\beta\sim N(\beta_0,\Phi_0^{-1})$, then the full conditional posterior distribution is also normal: 
\begin{equation}\label{eq beta response}
 \beta|\bm{Y}_{\T_1},\bm{Z}_{\T_1},\sigma^2 \sim N\big((\Phi_0+\bm{X}_{\T_1}^{\Tr} \tilde{\bm{\Phi}}\bm{X}_{\T_1})^{-1}(\Phi_0\beta_0+\bm{X}_{\T_1}^{\Tr}\tilde{\bm{\Phi}}\bm{Y}_{\T_1}), (\Phi_0+\bm{X}_{\T_1}^{\Tr}\tilde{\bm{\Phi}}\bm{X}_{\T_1})^{-1}\big) .
\end{equation}
We use robust adaptive Metropolis-Hastings steps \citep{vihola2012} to update $\sigma^2$ and $\theta$ targeting an acceptance probability $\approx 24\%$.
The sampling of $\bm{Y}_\T$ given $\beta, \theta, \sigma^2$ relies on:
\begin{align}
\tilde{\bm{Z}}_{\mathcal{T}_2} & \sim 
N(\bm{\Sigma}_{\T_2\T_1} \hat{\bm\Phi} (\bm{Y}_{\T_1}-\bm{X}_{\T_1}\beta), \bm{\Sigma}_{\T_2\T_2} - \bm{\Sigma}_{\T_2\T_1} \hat{\bm\Phi} \bm{\Sigma}_{\T_1\T_2}), \label{eq ZT response 1}  \\
\tilde{\bm{Y}}_{\T_2} & = \bm{X}_{\T_2}\beta + \tilde{\bm{Z}}_{\T_2}. \label{eq ZT response 2}
\end{align}

Algorithm \ref{alg AltGP response} summarizes the MCMC for Bayesian inference of the RadGP response model.

\begin{algorithm}[H]
\SetAlgoLined
Input training locations $\mathcal{T}_1\subset \Omega$, test locations $\mathcal{T}_2\subset\Omega$, covariates $\bm{X}_{\T_1}$, $\bm{X}_{\T_2}$ and response $\bm{Y}_{\T_1}$. \\
Set an approximation radius $\rho$, iteration constants $L_1, L_2$ and initial values for $\beta$, $\theta$, $\sigma$ and $Z_{\T_1}$.\\
Compute the initial decomposition $\tilde{\bm{\Phi}} = (\bm{I}_n-\hat{\bm{B}}^{\Tr})\hat{\bm{D}}^{-1}(\bm{I}_n-\hat{\bm{B}})$. \\
\For{$1\le l \le l_1$}{
Sample $\beta$ from equation (\ref{eq beta response});\\
Sample $\sigma$ and $\theta$ from using Metropolis Hastings updates;  \\
Compute the decomposition $\tilde{\bm{\Phi}} = (\bm{I}_n-\hat{\bm{B}}^{\Tr})\hat{\bm{D}}^{-1}(\bm{I}_n-\hat{\bm{B}})$. \\
\If{$l\ge l_2$}{Sample $\tilde{\bm{Z}}_{\T_2}$ according to equation (\ref{eq ZT response 1}) and compute $\tilde{\bm{Y}}_{\T_2}$ from equation (\ref{eq ZT response 2}).}
}
\KwResult{Output the posterior samples of $\beta, \theta, \delta$ and $\tilde{Z}_{\mathcal{T}_2}$.}
\caption{Posterior Sampling for Response RadGP Process}\label{alg AltGP response}
\end{algorithm}

Algorithm \ref{alg AltGP response} does not require the computation of spatial random effects. Hence by a similar analysis for the latent effects model, the overall computational complexity is $O(l_1 n (d_\beta^2+M_1^3) + l_2 n_{\mr{test}} M_3^3)$.

\section{Additional Experimental Studies}\label{sec:add exp}
\subsection{Prior Approximations Under Increasing Domain Asymptotics}
In this section, we employ an increasing domain asymptotic setting in a unidimensional space where the minimal separation distance is fixed at $q=1/50$. As the sample size $n$ increases, the size of domain $\Omega$ expands accordingly. We choose $\rho$ as a function of $n$ according to Corollary \ref{corollary} and study how the prior approximation error changes as $n$ increases. Specifically, we consider two covariance functions for the original Gaussian process:
Mat\'{e}rn with smoothness $\nu=3/2$ and generalized Cauchy with power $\lambda=5$. Using the notations in Corollary \ref{corollary}, we choose parameters $(\phi, \tau^2, \nu) = (17.95, 1, 3/2)$ for Mat\'{e}rn and $(\phi, \tau^2, \lambda, \delta)=(2.13, 1, 5, 1)$ for generalized Cauchy. The $\phi$ parameter in both cases are chosen such that $K_0(0.15)=0.25$. The approximation radius as a function of $n$ is chosen to be $\rho(n) = (\log n)^3/200$ for Mat\'{e}rn and $\rho(n) = n^{1/2}/20$ for generalized Cauchy. 

Because our main theorems only specify the upper bound of Wasserstein distance up to an absolute constant, it is necessary to numerically determine a ``tight'' constant. Specifically, for Mat\'{e}rn, we have set the theoretical upper bound to be $c_{out}n\exp(-c_{in}\rho)$. The constant $c_{in}$ is selected so that $\exp(-c_{in}\rho)$ decays at about the same speed as the Mat\'{e}rn covariance function in the domain of our simulation. The constant $c_{out}$ is numerically adjusted to be the smallest constant such that the theoretical upper bound of the squared Wasserstein-2 distance $W_2^2$ is no smaller than simulated values for all $n$. This procedure yields $c_{in} =13.68$ and $c_{out} =5.19\times10^{-9}$. For the generalized Cauchy covariance function, we set the theoretical upper bound to be $c_{ca} n / (1+\sqrt{n})^3$, where $c_{ca}$ is determined numerically to be $7.71\times10^{-3}$.

\begin{figure}[h]
    \centering
    \includegraphics[width=\textwidth]{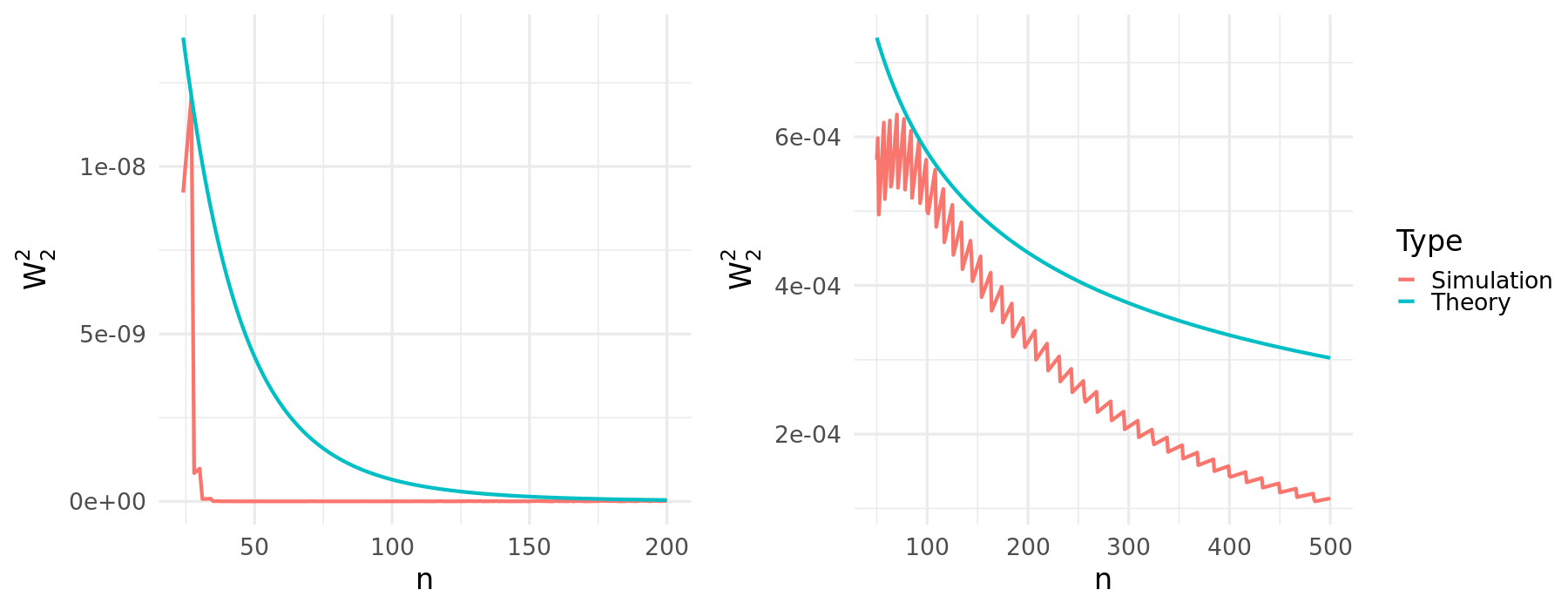}
    \caption{Comparison between simulation and theoretical bounds under increasing domain asymptotics. $X$-axis is the sample size. $Y$-axis is the squared Wasserstein-2 distance. Left: Mat\'{e}rn covariance function with smoothness $\nu=3/2$; Right: Generalized Cauchy covariance function with power $\lambda=5$.}
    \label{fig:asympt}
\end{figure}

The results for both simulated and theoretical upper bounds of $W_2^2$ are shown in Figure \ref{fig:asympt}. We can see for both Mat\'{e}rn and generalized Cauchy covariance functions, the simulation curves decay at a faster rate compared to the theoretical upper bound. For the Mat\'{e}rn covariance function, the theoretical curve decays at a faster-than-polynomial but slower-than-exponential rate, whereas the simulation curve exhibits a roughly exponential decay. For the generalized Cauchy covariance function, the simulation curve initially increases when the sample size is small, and then it decreases faster than the theoretical curve as the sample size becomes large. These findings indicate that our current theoretical upper bounds may not be sufficiently tight. The actual decay rate of the Wasserstein-2 distance between the RadGP and the true Gaussian process could be faster than what we have proved in Corollary \ref{corollary}.

\subsection{Parameter Estimation}\label{subsec:para}
\begin{figure}[h]
    \centering    \includegraphics[width=\textwidth]{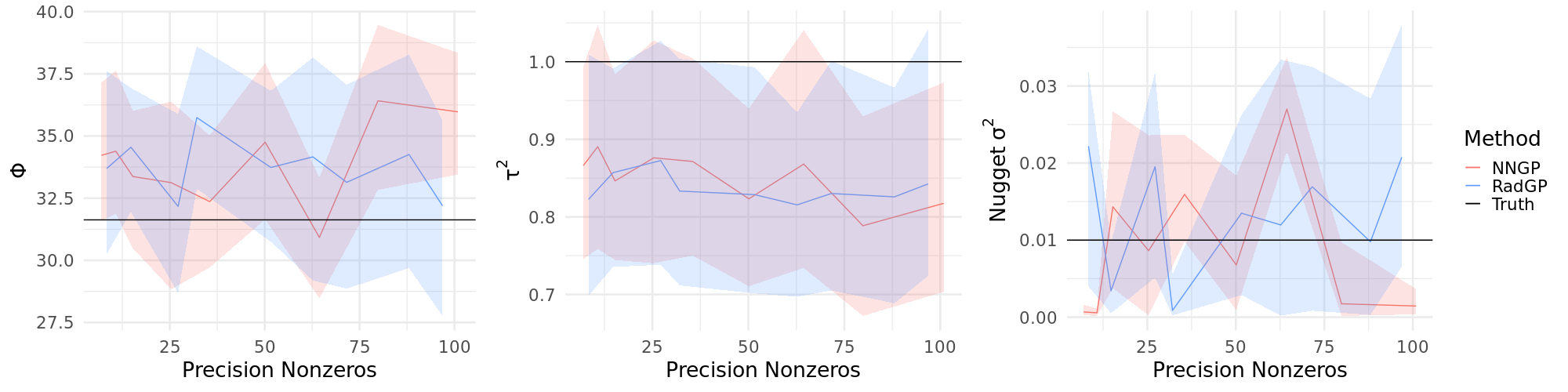}
    \caption{Posterior mean and Bayesian credible intervals for parameters $\phi, \tau^2$ and $\sigma^2$, where $X$-axis is the average number of nonzero elements per column in the precision matrix. Solid lines are the posterior means, and shaded regions are the $90\%$ Bayesian credible intervals.}
    \label{fig:para}
\end{figure}

Recall that the Mat\'{e}rn covariance function is
\begin{equation*}
K_0(x) = \frac{\tau^2 2^{1-\nu}}{\Gamma(\nu)} (\phi x)^{\nu} \mathcal{K}_{\nu} (\phi x).  
\end{equation*} 
In this section, we present the results of estimating the  parameters $\phi, \tau^2$ and $\sigma^2$ under the same settings described in Section \ref{subsec:posterior} of the main paper. The experiment is conducted once for each $x$-axis value, which represents the average number of nonzero elements per column in the precision matrix. We directly compute the posterior mean and $90\%$ credible intervals for all parameters from the posterior samples. The results are shown in Figure \ref{fig:para}. Generally, RadGP and NNGP exhibit similar performance in parameter estimation, with RadGP showing a slight advantage in estimating $\phi$ at large $x$-axis value. Neither method estimates $\tau^2$ and the nugget $\sigma^2$ very effectively, which is not surprising as distinguishing the Gaussian process and the nugget effect is known to be a difficult problem. According to the fixed domain asymptotic theory, when $\phi$ is known, the convergence rates of both frequentist estimators and Bayesian posteriors for $\tau^2$ will deteriorate from a parametric rate in the model without nugget effect to a much slower nonparametric rate in the model with nugget effect \citep{Chenetal00,Tangetal21,LiSunZhu23}. The presence of the nugget effect means that increasing the complexity of graphical structure does not significantly improve the accuracy of parameter estimation.

\subsection{MCMC Mixing and Conjugate Gradient Iterations}\label{subsec:MCMC}
\begin{figure}
    \centering
    \includegraphics[width=\textwidth]{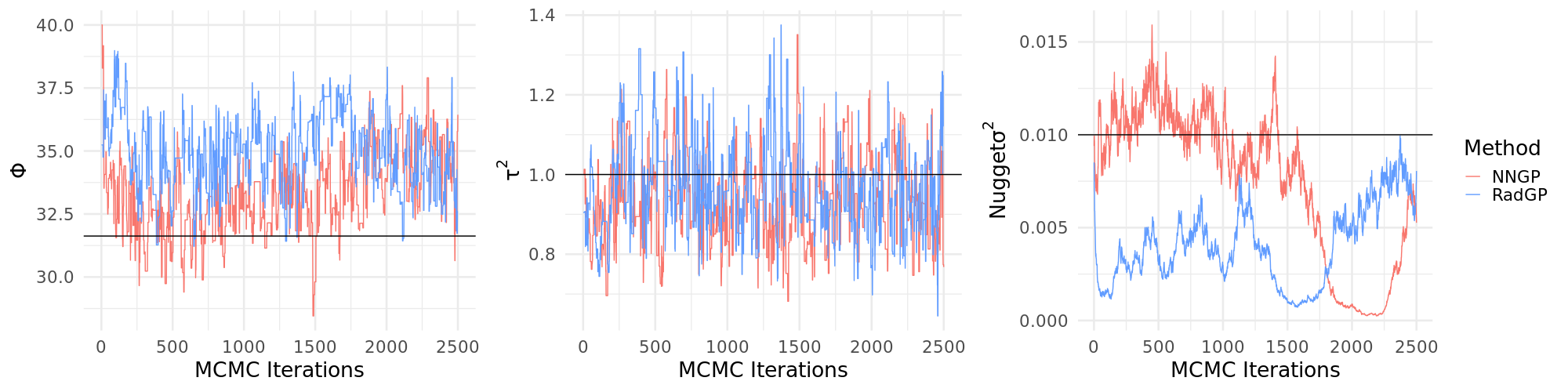}
    \caption{Trace plots of model parameters $\phi$, $\tau^2$ and $\sigma^2$ for both RadGP and NNGP methods.}
    \label{fig:trace}
\end{figure}

\begin{figure}
    \centering
    \includegraphics[width=\textwidth]{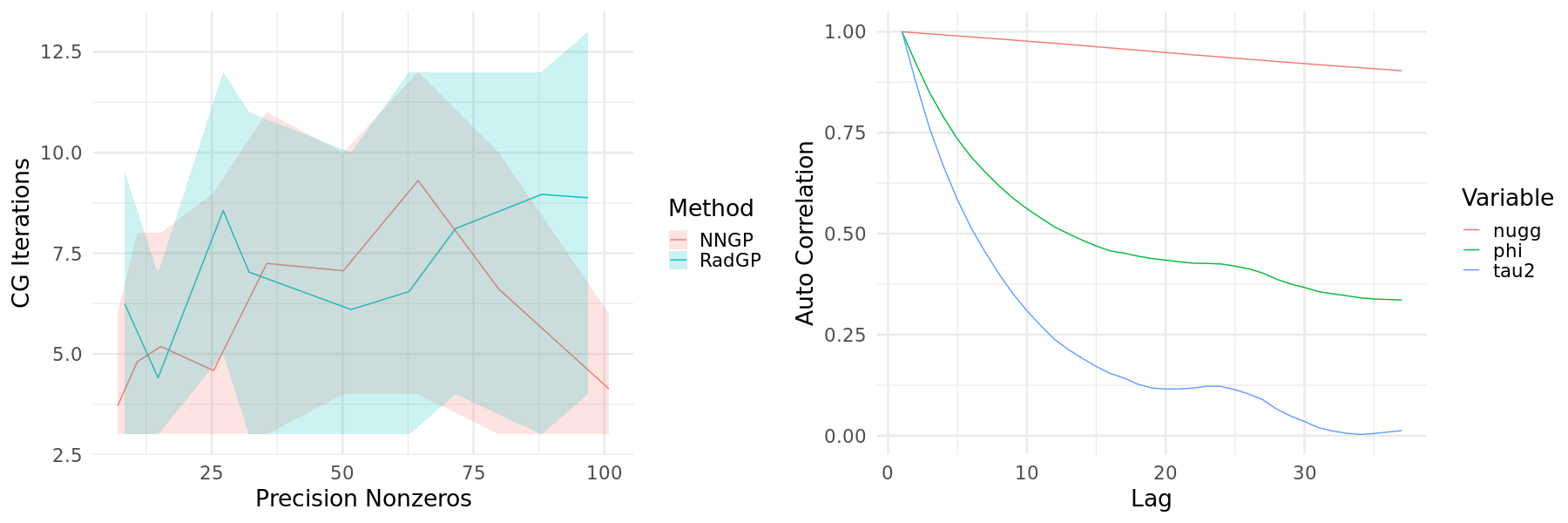}
    \caption{Left: Number of conjugate gradient iterations under different graph complexities, where $X$-axis is the average number of nonzero elements per column in the precision matrix. Solid lines are the mean values and shaded regions are the $90\%$ Bayesian credible intervals; Right: Autocorrelation functions for RadGP, where $X$-axis is the lag in the MCMC chains.}
    \label{fig:trace-cg}
\end{figure}

We study the posterior mixing of Algorithm \ref{alg cg} under the settings described in Section \ref{subsec:posterior}. Figure \ref{fig:trace} displays the trace plots of the Gaussian process parameters $\phi$, $\tau^2$ and $\sigma^2$ for both RadGP and NNGP methods. The radius $\rho$ for RadGP is set at $0.055$, and the number of parents for NNGP is set at $4$. With these parameter settings, the average number of nonzero elements per column in the precision matrix 
is $14.5$ for RadGP and $15$ for NNGP, which are very close.\footnote{Because the locations are distributed on a grid and the number of parents for NNGP can only take integer numbers, the average number of nonzero elements can only take discrete values, making it impossible to set exactly the same average number of nonzeros for RadGP and NNGP methods.} We also plot the autocorrelation function for RadGP in Figure \ref{fig:trace-cg} (right). We omit the trace plots under other graph complexities and the autocorrelation plot for NNGP, as they are very similar to those displayed here.

Despite the requirement for MCMC to sample all spatial random effects, our simulation results suggest that the mixing for the Gaussian process parameters $\phi$ and $\tau^2$ is quite good. The autocorrelation drops below $0.5$ when the lag exceeds $7$. Conversely, the posterior mixing for the nugget effect $\sigma^2$ is markedly slower. As discussed in Section \ref{subsec:para}, this phenomenon is attributed to the challenge of differentiating the nugget effect from the covariance function.

Figure \ref{fig:trace-cg} (left) displays the number of conjugate gradient iterations required under the settings described in Section \ref{subsec:posterior}. With a common tolerance of $10^{-6}$, the average number of conjugate gradient iterations is approximately $7.5$, with the $90\%$ credible interval $(2.6, 12.6)$, which is much smaller than the training sample size $n=1600$.

\subsection{Prior and Posterior Approximations under Higher Spatial Correlation}
\begin{figure}[h]
    \centering    \includegraphics[width=\textwidth]{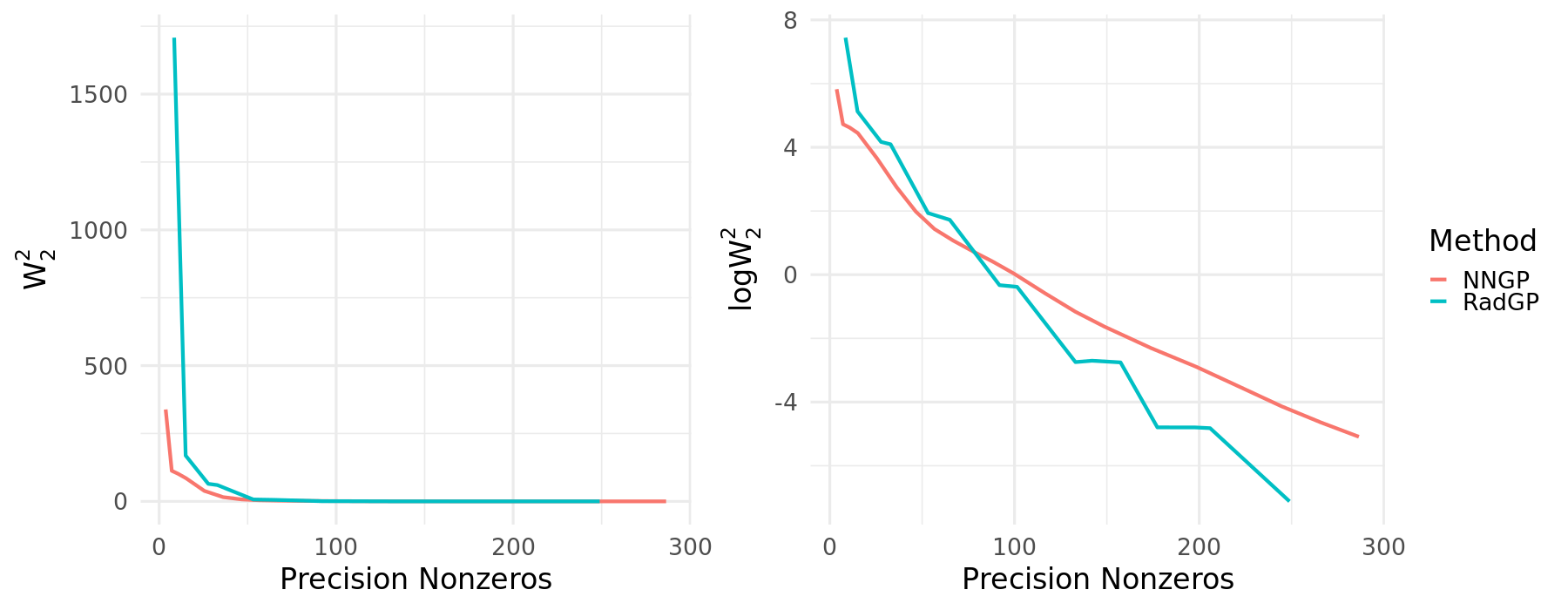}
    \caption{Prior comparison for Mat\'{e}rn covariance function under the high spatial correlation setting. X-axis is the average number of nonzero elements per column in the precision matrix. Y-axis is the squared Wasserstein-2 distance and its logarithm.}
    \label{fig:high-corr-prior}
\end{figure}

\begin{figure}[h]
    \centering    \includegraphics[width=0.8\textwidth]{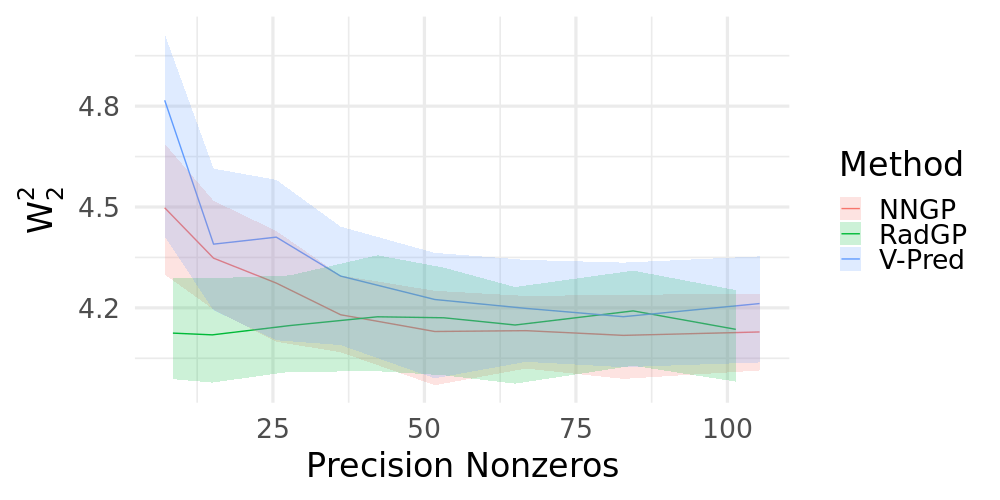}
    \caption{Posterior comparison under the high spatial correlation setting. X-axis is the average number of nonzero elements per column in the precision matrix. Y-axis is the squared Wasserstein-2 distance. Solid lines are mean values and shaded regions are $90\%$ Bayesian credible intervals.}
    \label{fig:high-corr-posterior}
\end{figure}

In this section, we perform numerical studies for the case with denser training sets and high spatial correlations. Specifically, the training locations are set on a $50\times 50$ equally spaced grid in $[0,1]^2$. The length-scale parameter $\phi$ is selected so that $K(0.15) = 0.15$. The above configuration generates much stronger spatial dependence and $50\%$ more training samples compared to the settings in Section \ref{subsec:posterior}. Other setups of the experiments are the same as described in Section \ref{subsec:posterior}. The prior and posterior comparisons are provided in Figures \ref{fig:high-corr-prior} and \ref{fig:high-corr-posterior}, respectively. For the prior comparison, we can see that RadGP performs better than NNGP when the average number of nonzero elements in the precision matrices is greater than $80$. For the posterior comparison, while the Bayesian credible bands for all three methods largely overlap when the graph structure is sufficiently complex, RadGP outperforms the other two methods when the graph complexity is lower. Based on these comparisons, we conclude that RadGP performs at least as well as the state-of-the-art methods like NNGP in settings characterzied by high spatial correlation.

\subsection{Prior and Posterior Approximations in 4-Dimensional Spaces}
\begin{figure}
    \centering
    \includegraphics[width=\textwidth]{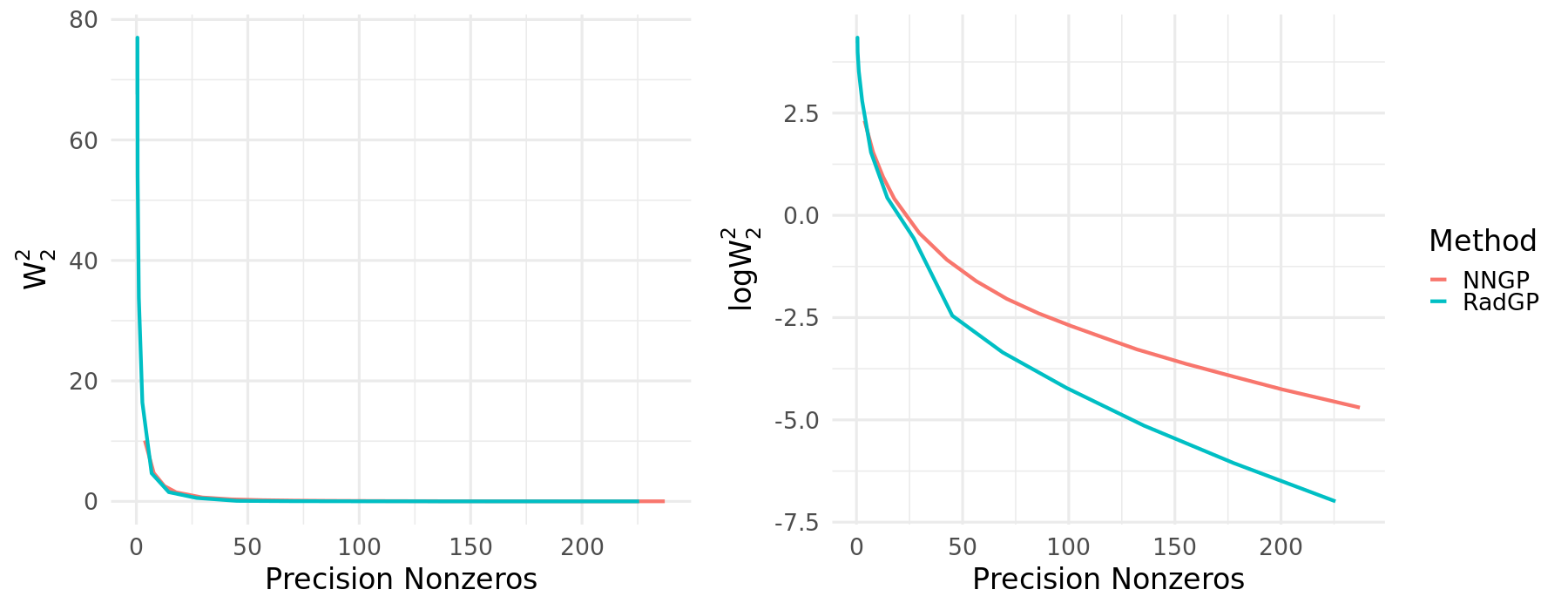}
    \caption{Prior comparison for Mat\'{e}rn covariance function in a 4-dimensional space. X-axis is the average number of nonzero elements per column in the precision matrix. Y-axis is the squared Wasserstein-2 distance and its logarithm.}
    \label{fig:prior-4d}
\end{figure}

\begin{figure}
    \centering
    \includegraphics[width=0.8\textwidth]{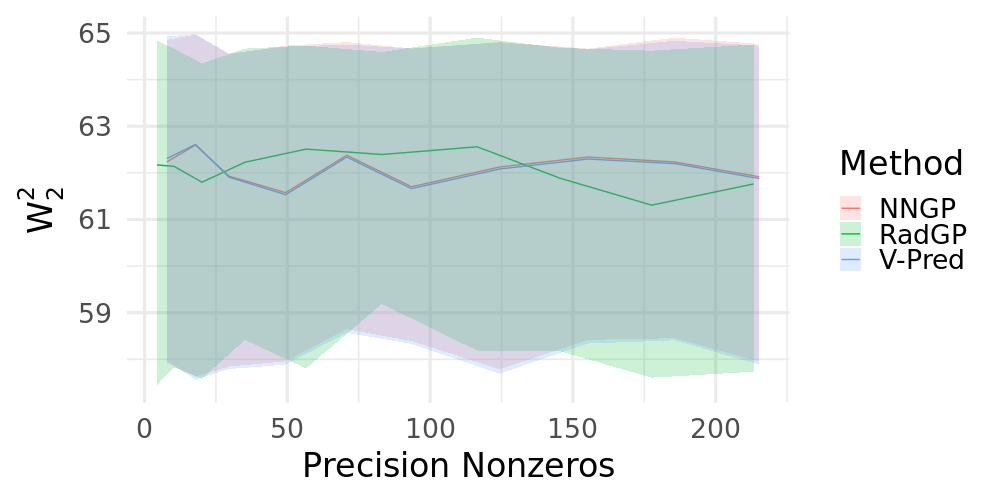}
    \caption{Posterior comparison for Mat\'{e}rn covariance function in a 4-dimensional space. X-axis is the average number of nonzero elements per column in the precision matrix. Y-axis is the squared Wasserstein-2 distance. Solid lines are mean values and shaded regions are $90\%$ Bayesian credible intervals.}
    \label{fig:posterior-4d}
\end{figure}

We examine the numerical performance of RadGP in the prior and posterior approximations with a 4-dimensional domain $\Omega$. Due to the curse of dimensionality, setting the training locations on a regular grid is no longer feasible. Therefore, we draw 2500 training samples from a random Latin Hypercube design in $[0,1]^4$ and fix them throughout the simulation studies in this section.\footnote{Changing training locations will also change the x-axis values (representing the average numbers of nonzeros in the precision matrix), which makes it harder to compare different methods. For simplicity, we fix the training locations. Based on our simulation studies, the key information conveyed by Figures \ref{fig:prior-4d} and \ref{fig:posterior-4d} is insensitive to the randomness in training locations.} The posterior comparison consists of $50$ repeated experiments. For each replication, the test data consist of another 2500 locations, complementing the training locations using random Latin hypercube design. The rest of the simulation settings are the same as described in Section \ref{subsec:posterior}.

The prior and posterior comparison are displayed in Figures \ref{fig:prior-4d} and \ref{fig:posterior-4d}. On the prior level, RadGP performs better than NNGP for the majority of x-axis values (the average numbers of nonzero elements in the precision matrix) considered in the comparison. However, all methods display equally large approximation errors for posterior prediction at the test locations, as measured in squared Wasserstein-2 distance. This indicates that the RadGP prior can still approximate the original Gaussian process prior accurately. However, due to the curse of dimensionality, all methods under investigation are unable to produce accurate posterior predictions of spatial random effects at unobserved locations with $n=2500$ training samples in a 4-dimensional space.

\subsection{Discussion on NNGP and V-Pred}

In the posterior comparison conducted in Section 5.2 of the main paper as well as Sections S5.4 and S5.5 above, NNGP performs similarly or sometimes even better than V-Pred in terms of the Wasserstein-2 distance between the posterior predictive distribution and the true predictive distribution. This phenomenon may appear confusing at first glance, since V-Pred allows dependence among the test locations and should ideally perform better. However, we must also take into account the impact of noise on prediction.
In our numerical implementation, both V-Pred and NNGP have a fixed number of parents for all spatial locations. Therefore, compared to NNGP, for each test location, V-Pred will replace some parents in the training
set with some parents in the test set that are closer. The presence of white noise makes the prediction of spatial random effects on test locations not so good as that on training locations,
which results in a worse prediction performance from V-Pred than NNGP on some occasions.

\section{Discussion on RadGP Theory and Tapering-based Gaussian Processes Theory} \label{sec:taper}

We would like to discuss the relations between our RadGP theory and the existing tapering-based Gaussian process theory (see for example, \citealt{Furetal06}, \citealt{Kauetal08}, \citealt{wang2011fixed}, and \citealt{ShaRup12}). Despite both theories working on the theoretical properties of Gaussian processes, these two theories exhibit the following four major differences.

First, the RadGP theory provided in this paper studies the finite sample approximation error bound of the radial neighbors Gaussian process, whereas the tapering-based Gaussion processes theory mainly studies the parameter estimation for the original Gaussian process. These two problems are different in the sense that an accurate finite sample approximation in Wassserstein distance does not imply good convergence rates in parameter estimation for the original process, and vice versa. To the best of our knowledge, we are not aware of any theory that provides the finite sample error bound for the approximation from either the covariance-tapered Gaussian processes or Vecchia approximated Gaussian processes to the
original Gaussian processes in Wasserstein distance under either increasing domain or fixed domain asymptotics.
Therefore, one major theoretical contribution in this paper is to introduce a new technique to spatial statistics, using the norm-controlled inversion of Banach algebras that can help quantify the approximation error for the radial neighbors Gaussian process. This new technique can potentially be further extended to other Vecchia approximation methods.

Second, the application scope of our RadGP theory and the previous tapering-based Gaussian processes theory are different. The theory of covariance tapering has mainly focused on parameter estimation for a parametric family of covariance functions, such as the isotropic Mat\'{e}rn covariance function, though in practice covariance tapering can be applied to any covariance function with good empirical performance for short range dependence. The setting of our RadGP theory is much more general and works for two large families of covariance functions ($\mathscr{Z}_{v_r}$ and $\mathscr{Z}_{c_r}$ defined in Equations (5) and (6) in the main paper) that only need to satisfy certain spatial decay conditions.

Third, since the Wasserstein-2 distance considered in our theory increases with the sample size (while decreases with the radius $\rho$), our current theorems and corollary work best for increasing domain and mixed domain asymptotics, but is not readily applicable to fixed domain (or infill) asymptotics. On the other hand, tapering-based Gaussian processes theory on parameter estimation, being well-developed in the literature, applies to both fixed domain asymptotics \citep{wang2011fixed} and increasing domain asymptotics \citep{ShaRup12}.

Finally, besides the difference in theory, RadGP and covariance tapering are different Gaussian process methods. Covariance tapering multiplies the original covariance function by another covariance function that is identically zero outside of a specific range. As such, covariance tapering induces sparsity only in the spatial covariance matrix, but not in the precision matrix. In contrast, RadGP induces sparsity in the precision matrix, aiming to accelerate the evaluation of Gaussian likelihoods, similar to other Vecchia approximation methods.

\bibliographystyle{chicago}
\bibliography{AGP}